\documentclass{amsart} 
\usepackage{amsmath}
\usepackage{amsthm}
\usepackage{hyperref}

\numberwithin{equation}{section}
\usepackage{amsfonts}
\usepackage{amssymb}
\usepackage{graphicx}
\usepackage{bbm}
\usepackage{tikz}
\usetikzlibrary{calc}
\usepackage[all]{xy}
\usepackage{blkarray}
\usepackage{indentfirst}
\usepackage{mathtools}
\usepackage{caption}
\usepackage{subcaption}
\newtheorem{theorem}{Theorem}[section]

\newtheorem{definition}{Definition}[section]

\newtheorem{proposition}{Proposition}[section]
\newtheorem{remark}{Remark}[section]

\newcommand{\Ad}{\mathrm{Ad}}
\newcommand{\ad}{\mathrm{ad}}

\newcommand{\SL}{\mathrm{SL}}

\newcommand{\CC}{\mathbb{C}}
\newcommand{\RR}{\mathbb{R}}

\newcommand{\ZZ}{\mathbb{Z}}

\newcommand{\Tr}{\mathrm{Tr}}

\newcommand{\frakk}{\mathfrak{k}}
\newcommand{\fraka}{\mathfrak{a}}
\newcommand{\frakg}{\mathfrak{g}}
\newcommand{\frakt}{\mathfrak{t}}
\newcommand{\frakh}{\mathfrak{h}}
\newcommand{\frakp}{\mathfrak{p}}

\newcommand{\fraku}{\mathfrak{u}}
\newcommand{\frakn}{\mathfrak{n}}
\newcommand{\frakm}{\mathfrak{m}}

\newcommand{\fraksl}{\mathfrak{sl}}

\newcommand{\Sym}{\mathrm{Sym}}

\newcommand{\fraksym}{\mathfrak{sym}}

\newcommand{\fraksp}{\mathfrak{sp}}
\newcommand{\ii}{\mathsf{i}}
\newcommand{\dr}{\mathsf{d}r}
\newcommand{\dd}{\mathsf{d}}
\newcommand{\dl}{\mathsf{d}l}

\newcommand{\ccc}{\mathsf{c}}

\newenvironment{spmatrix}{\left(\begin{smallmatrix}}{\end{smallmatrix}\right)}

\begin{document} 

\begin{abstract}
	Following the same method of \cite{MillerButtcane} and \cite{zz2019}, we calculate the $(\mathfrak{g},K)$ module structure for the principal series representation of $Sp(4,\RR)$. Furthermore, we introduced a hypergeometric generating function together with an inverse Mellin transform technique as an improvement to the method to calculate the intertwining operators. We have shown that the matrix entries of the simple intertwining operators for $Sp(4,\RR)$-principal series are Hahn polynomials, and the matrix entries of the long intertwining operator can be expressed as the constant term of the Laurent expansion of some hypergeometric function.
\end{abstract}
 
\title{Intertwining Operator for $Sp(4,\mathbb{R})$ and Orthogonal Polynomials } 
\author{Zhuohui Zhang} 
\email{zhuohui.zhang@weizmann.ac.il}
\maketitle

\section{An Introduction to the Results}
The intertwining operators for principal series representation was defined in \cite{KnappSteinIntertwining} and discussed in \cite{VoganWallachIntertwining}. \cite{KnappZuckermanIntertwining} and \cite{Muic} applied intertwining operators to study the decomposition of certain classes of principal series representations. In this paper, we have shown the following new result for long intertwining operator for the minimal principal series of $Sp(4,\RR)$:
\begin{theorem}\label{LongSp4R}
If the induction parameter $\delta$ of $Sp(4,\RR)$ principal series satisfies $\delta\in\{(0,0),(1,1)\}$ and $n- m_1\equiv \delta_1\text{ and } n+m_2\equiv \delta_2 \text{ mod }2$ for $i\in\{1,2\}$, the matrix entries $[A(\lambda)]^{j,n}_{m_1,m_2}$ for the long intertwining operator is the constant Laurent series coefficient of 
\begin{align}
	&[A(\lambda)]^{j,n}_{m_1,m_2}(t_1,t_2)=\nonumber\\
	&\frac{(-1)^{n}((2j)!)^2}{\left(\frac{\lambda_1+1}{2}\right)^{(\frac{j+n-\epsilon^{j,n}_{\delta}}{2})}\left(\frac{\lambda_1+1}{2}\right)^{(-\frac{j+n-\epsilon^{j,n}_{\delta}}{2})}\left(\frac{\lambda_1-\lambda_2+1}{2}\right)^{(j)}\left(\frac{\lambda_1+\lambda_2+1}{2}\right)^{(j)}}\times\nonumber\\
	&\frac{1}{c^j_{m_1}c^j_{m_2}}\frac{\left(\frac{\lambda_1-\lambda_2}{2}\right)^{(\frac{j+m_1-\epsilon^{j,n}_{\delta}}{2})}\left(\frac{\lambda_1+\lambda_2}{2}\right)^{(\frac{-j-m_2+\epsilon^{j,n}_{\delta}}{2})}}{\left(\frac{\lambda_2+1}{2}\right)^{\frac{m_2-n}{2}}\left(\frac{\lambda_2+1}{2}\right)^{-\frac{m_2-n}{2}}}\frac{\Gamma\left(\frac{1-\epsilon^{j,n}_{\delta}+j-m_1}{2}\right)}{\Gamma\left(\frac{1-\epsilon^{j,n}_{\delta}+j-m_2}{2}\right)}\times\nonumber\\
	&(1-t_1)^{\frac{-1+\epsilon^{j,n}_{\delta}-j+m_1}{2}}(1-t_2)^{\frac{-1-\epsilon^{j,n}_{\delta}+j-m_2}{2}}t_2^{\epsilon^{j,n}_{\delta}-2j}t_1^{-\epsilon^{j,n}_{\delta}}\times\nonumber\\
	&{}_2F_1\left(\begin{smallmatrix}-j+m_1,\frac{\lambda_1-\lambda_2-2j-1}{2}\\-2j\end{smallmatrix};t_1\right){}_2F_1\left(\begin{smallmatrix}-j-m_2,\frac{\lambda_1+\lambda_2-2j-1}{2}\\-2j\end{smallmatrix};t_2\right)\times\nonumber\\
	 &\,
   _5F_4\left(\begin{smallmatrix}-j+\epsilon^{j,n}_{\delta},1,\frac{-j+m_2+1+\epsilon^{j,n}_{\delta}}{2},\frac{-j-n-\lambda_1+1+\epsilon^{j,n}_{\delta}}{2},\frac{-j-m_2+\lambda_1+\lambda_2+\epsilon^{j,n}_{\delta}}{2}\\-j+\epsilon^{j,n}_{\delta},\frac{-j+m_1+1+\epsilon^{j,n}_{\delta}}{2},\frac{-j-n+\lambda_1+1+\epsilon^{j,n}_{\delta}}{2},\frac{-j-m_1-\lambda_1+\lambda_2+\epsilon^{j,n}_{\delta}}{2}+1\end{smallmatrix};\frac{t_2^2 \left(1-t_1\right)}{t_1^2\left(1-t_2\right)
   }\right).
\end{align}
where $
	\epsilon^{j,n}_{\delta} = \left\{\begin{smallmatrix}0&j-n\equiv\delta_i\text{ }\mathrm{mod}\text{ } 2\\1&j-n\not\equiv\delta_i\text{ }\mathrm{mod}\text{ } 2\end{smallmatrix}\right.
$.
\end{theorem}

The notations in the theorem above will be specified in the following sections.\\

This paper is the second part of my Ph.D thesis. I would like to thank my doctoral advisor Stephen D. Miller, Siddhartha Sahi, James Lepowsky and Doron Zeilberger for instructions and inspiring ideas. I also would like to thank Oleg Marichev from Wolfram Research for his support and generous sharing of his symbolic integration package.\\

\section{The Lie Group $Sp(4,\RR)$}
\subsection{The Structure of the Group $Sp(4,\RR)$}
The real symplectic group $G=Sp(4,\RR)$ is the subgroup of $SL(4,\RR)$ consisting of all elements $g\in SL(4,\RR)$ such that $g^t J g=J$, where
\[
	J = \begin{spmatrix}0&0&1&0\\0&0&0&1\\-1&0&0&0\\0&-1&0&0\end{spmatrix}.
\]
The real Lie algebra of the symplectic group is
\[
	\frakg=\fraksp(4,\RR) = \{X\in\fraksl(4,\RR)|X^t J+J X=0\}.
\]
A root space decomposition for the complexified Lie algebra $\frakg_\CC$ and its Chevalley basis are determined by the following choice of data:
\begin{enumerate}
	\item A Cartan subalgebra $\frakh_\CC$ generated by $H_1 = E_{1,1}-E_{3,3}, H_2 = E_{2,2}-E_{4,4}$;
	\item The simple roots $\alpha_1$ and $\alpha_2$ sending $t_1 H_1+t_2 H_2$ to
	\begin{align*}
	\begin{matrix}
		\alpha_1(t_1 H_1+t_2 H_2) = t_1-t_2, & \alpha_2(t_1 H_1+t_2 H_2) = 2 t_2
	\end{matrix}
	\end{align*}
	\item The simple coroots $\check{\alpha}_1$ and $\check{\alpha}_2$ sending $t_1 H_1+t_2 H_2$ to
	\begin{align*}
	\begin{matrix}
		\check{\alpha}_1(t_1 H_1+t_2 H_2) = t_1-t_2, & \check{\alpha}_2(t_1 H_1+t_2 H_2) =  t_2
	\end{matrix}
	\end{align*}
	\item The fundamental weights $\varpi_1$ and $\varpi_2$ in $\frakh_\CC^*$ satisfying $\langle\varpi_i, \check{\alpha}_j\rangle = \delta_{ij}$;
	\item The set of positive roots $\Delta^+(\frakg_\CC,\frakh_\CC) = \{\alpha_1,\alpha_2,\alpha_1+\alpha_2,2\alpha_1+\alpha_2\}$;
	\item $\rho_\CC=\frac{1}{2}\sum_{\alpha\in\Delta^+(\frakg_\CC,\frakh_\CC)}\alpha = \frac{4\alpha_1+3\alpha_2}{2}=\varpi_1+\varpi_2$;
	\item A basis for positive root spaces $\frakg_{\alpha}$ for $\alpha\in \Delta^+(\frakg_\CC,\frakh_\CC)$
	\begin{align*}
	\begin{matrix}X_{\alpha_1}=E_{1,2} - E_{4,3}&X_{\alpha_2}=E_{2,4}\\ X_{\alpha_1+\alpha_2}=E_{2,3}+E_{1,4}&X_{2\alpha_1+\alpha_2}=E_{1,3}
	\end{matrix}
	\end{align*}
	with the corresponding negative root vector given by $X_{-\alpha} = X_\alpha^t$; 
\end{enumerate}
We will also represent an element $\lambda\in \frakh_\CC^*$ by a pair of complex numbers $(\lambda_1,\lambda_2)$ where $\lambda_i = \lambda(H_i)$. In this coordinate, $\rho_\CC$ is represented by $(2,1)$.\\

On the real Lie group $Sp(4,\RR)$, there is a Cartan involution $\theta(g) = (g^t)^{-1}$ which determines a subgroup of fixed points $K = G^{\theta}$. The Lie algebra $\frakk\subset\frakg$ of the maximal compact subgroup $K\subset G$ is:
\[
	\frakk = \{\begin{spmatrix}A&B\\-B&A\end{spmatrix}|A\text{ antisymmetric}, B\in\fraksym(2)\}
\]
%Set
%\[
%	\ccc = \begin{spmatrix}\frac{\ii}{\sqrt{2}}&0&\frac{1}{\sqrt{2}}&0\\
%	0&\frac{\ii}{\sqrt{2}}&0&\frac{1}{\sqrt{2}}\\
%	-\frac{\ii}{\sqrt{2}}&0&\frac{1}{\sqrt{2}}&0\\
%	0&-\frac{\ii}{\sqrt{2}}&0&\frac{1}{\sqrt{2}}
%	\end{spmatrix}
%\]
%then
%\[
%\ccc\frakk\ccc^{-1} = \{\begin{spmatrix}U&0\\0&\bar{U}\end{spmatrix}|U\in\fraku(2)\}
%\]
The map $\begin{spmatrix}A&B\\-B&A\end{spmatrix}\mapsto A+\ii B$ from $\frakk$ to $2\times 2$ matrices identifies the Lie algebra $\frakk$ with $\fraku(2)$. The 4-dimensional Lie algebra $\frakk$ have generators $U_0,U_1, U_2, U_3$ corresponding to the 4 infinitesimal generators $\gamma_0,\gamma_1,\gamma_2,\gamma_3$ of $\fraku(2)$: 
\[
	\begin{smallmatrix}
		U_0 = \frac{1}{2}\begin{spmatrix}0&0&1&0\\0&0&0&1\\-1&0&0&0\\0&-1&0&0\end{spmatrix},&U_1=\frac{1}{2}\begin{spmatrix}0&0&0&1\\0&0&1&0\\0&-1&0&0\\-1&0&0&0\end{spmatrix}\\
		U_2 = \frac{1}{2}\begin{spmatrix}0&1&0&0\\-1&0&0&0\\0&0&0&1\\0&0&-1&0\end{spmatrix},&U_3= \frac{1}{2}\begin{spmatrix}0&0&1&0\\0&0&0&-1\\-1&0&0&0\\0&1&0&0\end{spmatrix}
	\end{smallmatrix}
\]
They satisfy a commutation relation:
\begin{align*}
	[U_0,U_i]=0, [U_i,U_j] = -\epsilon_{ijk}U_k
\end{align*}
where the \emph{Levi-Civita symbol} $\epsilon_{ijk}$ takes value $1$ if $(ijk)$ is an even permutation of $(123)$, $-1$ if $(ijk)$ is an odd permutation of $(123)$, and 0 if two or more elements in $\{i,j,k\}$ are equal.
\subsection{The Cartan Subgroups of $Sp(4,\RR)$}
\subsubsection{The Maximally Compact Cartan Subalgebra}
A Cartan subalgebra of $K\cong U(2)$ is $\frakt = \RR (U_0+U_3)\oplus \RR (U_0-U_3)$. The basis vectors $U_0\pm U_3$ are deliberately chosen so that they can be related to $H_{\alpha_1}$ and $H_{\alpha_2}$ by Cayley transforms. Let the the simple roots $\beta_1,\beta_2$ act on an element $t_1(U_0+U_3)+t_2(U_0-U_3)\in\frakt_\CC$ by
\begin{align*}
	\beta_1(t_1(U_0+U_3)+t_2(U_0-U_3)) &= \ii (t_1-t_2),\\\beta_2(t_1(U_0+U_3)+t_2(U_0-U_3)) &= 2\ii t_2.
\end{align*}
In the Vogan diagram (see \cite[Section~VI.8]{KnappRep}), the shorter root is named to be the compact root. In this situation, all roots are imaginary:
\begin{center}
\begin{tikzpicture}
	\draw
	(0,0) circle [radius=.1] node [above] {$\beta_1$};
	\draw[fill=black]
	(1,0) circle [radius=.1] node [above] {$\beta_2$};
	\draw
	(0.1,-.04) --++ (1,0)
	(0.1,+.04) --++ (1,0);
	\draw
	(.5,0) --++ (60:.2)
	(.5,0) --++ (-60:.2);
\end{tikzpicture}
\end{center}
We can decompose the set of positive roots $\Delta^+(\frakg_\CC,\frakt_\CC)$ into the union of the set of compact roots and noncompact roots:
\begin{align*}
	\Delta^+_c(\frakg_\CC,\frakt_\CC) &= \{\beta_1\}\\
	\Delta^+_{nc}(\frakg_\CC,\frakt_\CC) &= \{\beta_2,\beta_1+\beta_2,2\beta_1+\beta_2\};
\end{align*}
The compact root vectors are
\begin{align*}
	v_{\pm\beta_1} &= \frac{1}{\ii}(U_1\pm\ii U_2).
\end{align*}
These roots are displayed in the following picture, where the gray color stands for noncompact roots, and the light gray color stands for compact roots:
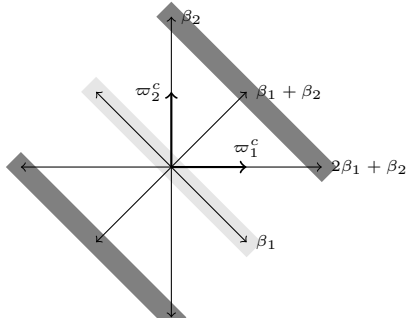
\begin{center}
\begin{tikzpicture}[scale=1]
	\clip (-3,-2) rectangle (3.2,2.2);
		\fill[fill=gray!100] (-0.2,2) -- (-0,2.2) -- (2.2,0) -- (2,-0.2);
		\fill[fill=gray!100] (-2.2,0) -- (-2,0.2) -- (0.2,-2) -- (0,-2.2);
		\fill[fill=gray!20] (-1.2,1) -- (-1,1.2) -- (1.2,-1) -- (1,-1.2);
		\draw[->] (0,0) -- (1,-1) node[right] {\tiny $\beta_1$};
		\draw[->]	(0,0)--(0,2) node[right] {\tiny $\beta_2$};
		\draw[->]	(0,0)--(1,1) node[right] {\tiny $\beta_1+\beta_2$};
		\draw[->] (0,0)--(2,0) node[right] {\tiny $2\beta_1+\beta_2$};
		\draw[->]	(0,0)--(-2,0);
		\draw[->]	(0,0)--(-1,-1);
		\draw[->]	(0,0)--(-1,1);
		\draw[->]	(0,0)--(0,-2);
		\draw[->, thick] (0,0) -- (1,0)  node[above] {\tiny $\varpi^c_1$};
		\draw[->, thick] (0,0) -- (0,1)  node[left] {\tiny $\varpi^c_2$};
%		\fill[fill=gray!20] (0.05,0.05) --  (1.782050,1.05) -- (0.05,2.05) -- (0.05,0.05);
%		\draw[dashed] (-1.732050,-1) -- (1.732050,1);
%		\draw[dashed] (-1.732050,1) -- (1.732050,-1);
%		\draw[dashed] (0,-2) -- (0,2);
	\end{tikzpicture}
\end{center}
The subspace $\frakt_\CC\oplus \frakg_{\beta_1}\oplus\frakg_{-\beta_1}$of $\frakg_\CC$ generated by the vectors $\{v_{\beta_1},v_{-\beta_1}, U_0+U_3, U_0-U_3\}$ is isomorphic to $\mathfrak{gl}(2,\CC)$. The embedding of this Lie subalgebra into $\frakg_\CC$ corresponds to the embedding of $U(2)$ into $G$ as $K$.\\
\indent The -1 eigenspace $\frakp_\CC$  for $\theta$ is generated by the following noncompact root vectors:
%\[
%	\frakp = \{\begin{spmatrix}t_1&s&u-v&r\\s&t_2&r&u+v\\u-v&r&-t_1&-s\\r&u+v&-s&-t_2\end{spmatrix}\}
%\]
\begin{align}
	v_{\pm(2\beta_1+\beta_2)}&=\frac{1}{2}(H_1\pm \ii(X_{2\alpha_1+\alpha_2}+X_{-2\alpha_1-\alpha_2}))\label{V1}\\
	v_{\pm\beta_2}&=\frac{1}{2}(H_2 \pm \ii(X_{\alpha_2}+X_{-\alpha_2}))\label{V2}\\
	v_{\pm(\beta_1+\beta_2)}&= \frac{1}{2}(X_{\alpha_1+\alpha_2}+X_{-\alpha_1-\alpha_2} \mp \ii(X_{\alpha_1}+X_{-\alpha_1}))\label{V3}.
\end{align}
\subsubsection{A Cayley transform and the maximal noncompact Cartan subalgebra}
There are 4 conjugacy classes of Cartan subgroups of $Sp(4,\RR)$ displayed in the diagram below, connected by Cayley transforms corresponding to the roots $\beta_2$ and $2\beta_1+\beta_2$:
	\begin{equation}\label{CayleyOnCartanSubgroups}
		\xymatrix{
		& \mathbb{S}^1\times\mathbb{S}^1 \ar[dl]_{\ccc_{\beta_2}} \ar[dr]^{\ccc_{2\beta_1+\beta_2}}&\\
		\RR^{\times}\times \mathbb{S}^1 \ar[dr]_{\ccc_{2\beta_1+\beta_2}} & & \mathbb{S}^1\times\RR^{\times} \ar[dl]^{\ccc_{\beta_2}}\\
		& \RR^{\times}\times \RR^{\times}
		}
	\end{equation}
Consider the \emph{Cayley} transforms $\ccc _{\beta} = \Ad\exp(\frac{\pi}{4}(\overline{v_{\beta}}-v_{\beta}))$ for a noncompact root $\beta$, when $\beta = \beta_2\text{ or }2\beta_1+\beta_2$, the Cayley transform is nontrivial:
\begin{align*}
	\ccc _{\beta_2} &=  \Ad\exp(\frac{\pi}{4}(\overline{v_{\beta_2}}-v_{\beta_2}))\\
	\ccc _{2\beta_1+\beta_2} &=  \Ad\exp(\frac{\pi}{4}(\overline{v_{2\beta_1+\beta_2}}-v_{2\beta_1+\beta_2})).
\end{align*}
The two Cayley transforms $c_{\beta_2}$ and $c_{2\beta_1+\beta_2}$ commute with each other. If we start from the real Cartan subalgebra $\frakt$, the other real Cartan subalgebras of $Sp(4,\RR)$ can be obtained by applying the Cayley transforms $\ccc _{\beta_2}$ and $\ccc _{2\beta_1+\beta_2}$:
\[
		\xymatrix{
		& \RR(U_0+U_3)\oplus \RR(U_0-U_3)\ar[dl]_{\ccc_{\beta_2}} \ar[dr]^{\ccc_{2\beta_1+\beta_2}}&\\
		\RR(U_0+U_3)\oplus \RR\ii H_2 \ar[dr]_{\ccc_{2\beta_1+\beta_2}} & & \RR\ii H_1\oplus \RR(U_0-U_3) \ar[dl]^{\ccc_{\beta_2}}\\
		& \RR\ii H_1\oplus \RR\ii H_2
		}
\]
This diagram shows how the generators of the Cartan subalgebras are mapped to each other by the Cayley transforms shown in the diagram (\ref{CayleyOnCartanSubgroups}). We define the \emph{maximally noncompact Cartan subalgebra} $\fraka = \RR H_1\oplus\RR H_2\subset\frakp$. The roots $\alpha_i$ acts on $\frakh_\CC = \fraka\otimes_\RR \CC$, and they can be related to the roots $\beta_i$ on the compact Cartan subgroup by \[\alpha_i\circ (\ccc _{\beta_2} \ccc _{2\beta_1+\beta_2}) = \beta_i.\] Moreover, we can apply the composite Cayley transform on the root vectors $v_{\beta}$:
\begin{align*}
	 \ccc _{\beta_2} \ccc _{2\beta_1+\beta_2}v_{\beta}=\ccc _{2\beta_1+\beta_2} \ccc _{\beta_2} v_{\beta} &= \begin{cases}\ii X_{\beta\circ(\ccc _{\beta_2} \ccc _{2\beta_1+\beta_2})^{-1}}&\text{ if }\beta\in\{\beta_2,2\beta_1+\beta_2\}\\X_{\beta\circ(\ccc _{\beta_2} \ccc _{2\beta_1+\beta_2})^{-1}}&\text{ if }\beta\in\{\beta_1+\beta_2\}\end{cases}.
%	 \ccc _{\beta_2} \ccc _{2\beta_1+\beta_2}v_{\beta}=\ccc _{2\beta_1+\beta_2} \ccc _{\beta_2} v_{\beta} &=  X_{\beta}\text{ if }\beta\in\{\alpha_1+\alpha_2\}
%	 \ccc _{\beta_2} \ccc _{2\beta_1+\beta_2}[v_{\beta},v_{-\beta}] =\ccc _{2\beta_1+\beta_2} \ccc _{\beta_2} [v_{\beta},v_{-\beta}] &=  H_{\beta}\\
%	\ccc _{2\beta_1+\beta_2} \ccc _{\beta_2} (U_0+U_3) &=\ii H_1\\
%	\ccc _{2\beta_1+\beta_2} \ccc _{\beta_2} (U_0-U_3) &= \ii H_2
\end{align*}
%\begin{align*}
%	P_1 &= v_{1,1,1}=\sqrt{2}i(H_1+i(X_{2\alpha_1+\alpha_2}+Y_{2\alpha_1+\alpha_2}))\\
%	P_2 &= v_{1,1,-1}=\sqrt{2}i(H_2 + i(X_{\alpha_2}+Y_{\alpha_2}))\\
%	P_3 &= v_{1,1,0}= X_{\alpha_1+\alpha_2}+Y_{\alpha_1+\alpha_2} - i(X_{\alpha_1}+Y_{\alpha_1})\\
%	P_4 &= v_{1,-1,-1}=\sqrt{2}i(H_1-i(X_{2\alpha_1+\alpha_2}+Y_{2\alpha_1+\alpha_2}))\\
%	P_5 &= v_{1,-1,1}=\sqrt{2}i(H_2 - i(X_{\alpha_2}+Y_{\alpha_2}))\\
%	P_6 &= v_{1,-1,0}=X_{\alpha_1+\alpha_2}+Y_{\alpha_1+\alpha_2} + i(X_{\alpha_1}+Y_{\alpha_1})
%\end{align*}
\subsubsection{Weyl Group}
From the embeddings $\Phi_{\alpha_1},\Phi_{\alpha_2}$ of $SL(2,\RR)$ into $Sp(4,\RR)$ given by the simple roots $\alpha_1,\alpha_2$, the simple Weyl reflections corresponding to these two roots are:
\begin{align*}
	w_{\alpha_i} = \Phi_{\alpha_i}\begin{spmatrix}0&1\\-1&0\end{spmatrix}=\exp\left(\frac{\pi}{2}(X_{\alpha_i}-X_{-\alpha_i})\right).
\end{align*}
The Lie algebra $\frakk$ has basis $\{X_{\alpha}-X_{-\alpha}\}_{\alpha\in \Delta^+}$, and the relationship between this basis and the basis $\{U_i\}$ is: 
\begin{align*}
\begin{matrix}
	U_0+U_3= X_{2\alpha_1+\alpha_2}-X_{-2\alpha_1-\alpha_2}&
	2U_1 = X_{\alpha_1+\alpha_2}-X_{-\alpha_1-\alpha_2}\\
	2U_2 = X_{\alpha_1}-X_{-\alpha_1}&
	U_0-U_3 = X_{\alpha_2}-X_{-\alpha_2}.
\end{matrix}
\end{align*}
Therefore, the simple reflections $w_{\alpha_i}$ can be expressed as
\begin{align}
		w_{\alpha_1} = \exp(\pi U_2), w_{\alpha_2} = \exp\left(\frac{\pi}{2}(U_0-U_3)\right).\label{Sp4Weyl}
\end{align}
Under the basis $\{\alpha_1,\alpha_2\}$, the actions of the simple reflections on an element $n_1\alpha_1+n_2\alpha_2\in \frakh^*_\CC$ are:
\begin{align*}
	w_{\alpha_1}(n_1\alpha_1+n_2\alpha_2) &= (2n_2-n_1)\alpha_1+n_2\alpha_2\\
	w_{\alpha_2}(n_1\alpha_1+n_2\alpha_2) &= n_1\alpha_1+(n_1-n_2)\alpha_2.
\end{align*}
If we represent $\lambda\in \frakh_\CC^*$ by $(\lambda_1,\lambda_2)$, the actions of $w_{\alpha_i}$ on $(\lambda_1,\lambda_2)$ where $\lambda_i = \lambda(H_i)$ is
\begin{align*}
	w_{\alpha_1}(\lambda_1,\lambda_2)&= (\lambda_2,\lambda_1)\\
	w_{\alpha_2}(\lambda_1,\lambda_2)&=(\lambda_1,-\lambda_2).
\end{align*}
The action of the simple reflections on the nilpotent radical satisfies:
\begin{align*}
	\overline{\frakn}\cap \Ad(w_{\alpha_i})^{-1}\frakn &= \RR X_{-\alpha_i}.
\end{align*}

\subsubsection{Harish-Chandra isomorphism}
Since $\mathfrak{sp}(4,\CC)$ is rank 2, there are two Casimir elements $\Omega_2$ and $\Omega_4$ with degree 2 and degree 4 respectively. They generate the center $Z(\frakg_\CC)$ of the complexified universal enveloping algebra. Consider the adjoint representation $\ad$ on $\frakg_\CC$, then the corresponding element in the center is \cite{Yang}:
\begin{align*}
	\Omega_2 =& \sum_{i,j}\Tr(\ad(X_i)\ad(X_j)) \tilde{X}_i\tilde{X}_j=\sum_{i,j}B(X_i,X_j) \tilde{X}_i\tilde{X}_j\\
	=&\frac{1}{12}(H_1^2+H_2^2+4H_1+2 H_2 + 2(X_{-\alpha_1}X_{\alpha_1}+2X_{-\alpha_2}X_{\alpha_2} \\
	&+ X_{-\alpha_1-\alpha_2}X_{\alpha_1+\alpha_2}+2 X_{-2\alpha_1-\alpha_2}X_{2\alpha_1+\alpha_2}))
\end{align*}
where $\{\tilde{X}_i\}$ is a dual basis of the basis $\{X_i\}$ of $\frakg_\CC$ under the Killing form on $\frakg$. The central element $\Omega_4$ corresponding to the standard representation is a quartic element in $U(\frakg_\CC)$ which we will only write down the $U(\frakh_\CC)$ part explicitly:
\begin{align*}
	\Omega_4 =& \frac{1}{5184} ((H_1+H_2+3)^4-2(2H_1H_2+2H_1+4H_2+7)((H_1+2)^2+&\\(H_2+1)^2)-11)
	&+\text{other terms}
\end{align*}
Since $\rho = \varpi_1+\varpi_2$, the Harish-Chandra homomorphism $\gamma'$ will map each $H_i$ to:
\[
	\gamma'(H_1) = H_1-2, \gamma'(H_2) = H_2-1
\]
The image of $\Omega_2$ and $\Omega_4$ under the Harish-Chandra isomorphism $\gamma'$ are:
\begin{align*}
	\gamma'(\Omega_2) &= \frac{1}{12}(H_1^2+H_2^2-5)\\
	\gamma'(\Omega_4) &= \frac{1}{5184} (H_1^4+H_2^4+6H_1^2H_2^2-6(H_1^2+H_2^2)-11).
\end{align*}

\section{Iwasawa decomposition on the group level}\label{IwasawaRefl}
We would like to compute the Iwasawa decomposition $kman$ of an element $\bar{n}\in \overline{N}$. For a root $\alpha$ with $\chi_{-\alpha}(t) = \exp(t X_{-\alpha})\in\overline{N}$, the corresponding Lie group $SL(2,\CC)$ embedded into $G_\CC$ has generators:
%\begin{lemma}
\[
	\chi_{\alpha}(t) = \exp(t X_{\alpha}), \chi_{-\alpha}(t) = \exp(t X_{-\alpha}), h_{\alpha}(t) = \exp(t H_\alpha)
\]
In terms of these matrix generators of $SL(2,\CC)$, the Iwasawa decomposition of $\chi_{-\alpha}(t)$ is
	\begin{align}
		\chi_{-\alpha}(t) &= \kappa_\alpha(t) h_{\alpha}\left(\sqrt{1+t^2}\right) \chi_{\alpha}\left(\frac{t}{1+t^2}\right)\label{IwasawaGroup}
	\end{align}
where
\[
	\kappa_\alpha(t) = \exp \left(\arctan(-t)(X_{\alpha}-X_{-\alpha})\right)
\]
is an element in the maximal compact subgroup $K_0\cong SO(2,\RR)$.
In general, for any simple real root $\alpha$ for $G = Sp(4,\RR)$, we would like to study at the embedding $\phi_{\alpha}$ of a Lie algebra $\fraksl(2,\RR)$ into $\frakg$, and the corresponding homomorphism of a Lie group $\Phi_{\alpha}:SL(2,\RR)\longrightarrow G$. The image of $\Phi_{\alpha}$ is fixed by $\theta$, and the embedding for the Lie algebra satisfies:
\begin{enumerate}
\item $\phi_{\alpha}\begin{spmatrix}0&1\\0&0\end{spmatrix} = X_{\alpha},\phi_{\alpha}\begin{spmatrix}0&0\\1&0\end{spmatrix}=X_{-\alpha}, \phi_{\alpha}\begin{spmatrix}1&0\\0&-1\end{spmatrix} = H_{\alpha}$,
\item $\Phi_{\alpha}\begin{spmatrix}1&t\\0&1\end{spmatrix} = \chi_{\alpha}(t),\Phi_{\alpha}\begin{spmatrix}1&0\\t&1\end{spmatrix}=\chi_{-\alpha}(t), \Phi_{\alpha}\begin{spmatrix}t&0\\0&1/t\end{spmatrix} = h_{\alpha}(t)$.
\end{enumerate}
We can see that the homomorphisms $\Phi_{\alpha}$ and $\phi_{\alpha}$ respect the Cartan involution: if we denote $\theta'$ as the Cartan involution on $SL(2,\RR)$ and $\fraksl(2,\RR)$, we can see that $\Phi_{\alpha}\circ \theta' = \theta\circ\Phi_{\alpha}$ and $\phi_{\alpha}\circ \theta' = \theta\circ\phi_{\alpha}$. Therefore, $\mathrm{Im}(\Phi_{\alpha})\cap K = \Phi_{\alpha} (K_0)$ where $K_0$ is the maximal compact subgroup $SO(2,\RR)$ in $SL(2,\RR)$.

Recall that under any simple reflection $w_{\alpha_i}$, there is a root vector $X_{-\alpha_i}$ such that $\bar{\frakn}\cap w_{\alpha_i}^{-1} \frakn w_{\alpha_i} = \RR X_{-\alpha_i}$. We can factor the nilpotent group $w_{\alpha_i}^{-1}N w_{\alpha_i}$ into
\[
	w_{\alpha_i}^{-1}N w_{\alpha_i} = N_{-\alpha} N'\text{ where }N'\subset N\text{ and } N_{-\alpha} = \{e^{t X_{-\alpha_i}}|t\in\RR\}.
\]
Using this factorization of the group $w_{\alpha_i}^{-1}N w_{\alpha_i}$, we can find a coordinate system on $w_{\alpha_i}^{-1}N w_{\alpha_i}$ which is consistent with the Iwasawa decomposition. Namely, for any $n\in N$, there is a $t\in\RR$ and $n'\in N'$ such that
\[
	w_{\alpha_i}^{-1}n w_{\alpha_i} = \kappa_{\alpha}(t)h_{\alpha}\left(\sqrt{1+t^2}\right)\chi_{\alpha}\left(\frac{t}{1+t^2}\right)n'.
\]

\section{Parabolic Subgroups of $Sp(4,\RR)$}
\subsection{Minimal Parabolic Subgroup and Jacobi Parabolic Subgroup}
Denote $P_0$ as the \emph{minimal parabolic} subgroup of $Sp(4,\RR)$, with the Lie algebra having a Levi decomposition:
\[
	\frakp_0 = \fraka \oplus \frakn_0, \text{ where }\frakn_0 = \bigoplus_{\alpha\in\Delta^+}\frakg_\alpha.
\]
For any real root $\alpha$ of $\Delta(\frakg_\CC,\mathfrak{h}_\CC)$, consider the homomorphism $\phi_\alpha: \fraksl(2,\RR)\longrightarrow \frakg$, which exponentiates to a homomorphism $\Phi_\alpha:\SL(2,\CC)\longrightarrow G_\CC$. Letting $\gamma_\alpha$ be the image of $\begin{spmatrix}-1&0\\0&-1\end{spmatrix}$, then
\begin{align*}
	\gamma_\alpha=\exp\frac{ 2\pi i H_\alpha}{|\alpha|^2} = \exp\pi(X_\alpha-\theta X_\alpha)
\end{align*}
Recalling the statements in \cite[Section~VII.4]{KnappLieGroup}, every nonidentity element of the centralizer $M$ of $\fraka$ in $K$ has order 2, and is generated by $\gamma_\alpha$ for all real roots $\alpha$. In the case of $Sp(4,\RR)$,
\begin{align*}
	M =\{\gamma_{\alpha_2}^{\epsilon_1}\gamma_{2\alpha_1+\alpha_2}^{\epsilon_2}|\epsilon_i = 0\text{ or }1\}\cong \ZZ_2\times\ZZ_2,
\end{align*}
The generators $\gamma_{\alpha_2},\gamma_{2\alpha_1+\alpha_2}$ of $M$ can also be expressed in terms of $U_i$:
\[
	\gamma_{\alpha_2} = \exp(\pi (U_0-U_3)), \gamma_{2\alpha_1+\alpha_2} = \exp(\pi (U_0+U_3)).
\]
The minimal parabolic subgroup $P_0$ has a Levi decomposition $P_0 = MA_0 N_0$  $M$ is the centralizer of $\fraka$ in $K$, and $A_0$ and $N_0$ are the analytic subgroups formed by exponentiating $\fraka$ and $\frakn_0$, respectively. \\
There are two proper standard parabolic subgroups that contain $P_0$. The parabolic subgroup $P_1$ with an abelian nilpotent radical is called a \emph{Siegel} subgroup. The Levi decomposition $P_1=M_1 N_1$ has the form:
\begin{align*}
	M_1 &= \{\begin{spmatrix}A&0\\0&(A^t)^{-1}\end{spmatrix}, A\in GL(2,\RR)\}\cong GL(2,\RR)\\
	N_1 &= \{\begin{spmatrix}1&0&x_4&x_3\\0&1&x_3&x_2\\0&0&1&0\\0&0&0&1\end{spmatrix}, x_i\in\RR\}\cong \RR^3.
\end{align*}
Their Lie algebras have restricted root space decompositions:
\begin{align*}
\frakm_1 &= \fraka\oplus\frakg_{\alpha_1}\oplus\frakg_{-\alpha_1}\\
\frakn_1 &=\bigoplus_{\alpha\in\{\alpha_2,\alpha_1+\alpha_2,2\alpha_1+\alpha_2\}}\frakg_\alpha.
\end{align*} There is another class of parabolic subgroup called the \emph{Jacobi} (or \emph{Heisenberg} and in some literature also called \emph{Klingen}) parabolic subgroup $P_2=M_2N_2$, having a Levi decomposition:
\begin{align*}
	M_2 &= \{\begin{spmatrix}h_1&0&0&0\\0&a&0&b\\0&0&1/h_1&0\\0&c&0&d\end{spmatrix},\begin{spmatrix}a&b\\c&d\end{spmatrix}\in SL(2,\RR), h_1\in\RR^{\times}\}\cong SL(2,\RR)\times \RR^{\times}\\
	N_2 &=  \{\begin{spmatrix}1&x_1&x_4&x_3\\0&1&x_3&0\\0&0&1&0\\0&0&-x_1&1\end{spmatrix}, x_i\in\RR\}.
\end{align*}
The Lie algebras of each subgroup have restricted root space decompositions:
\begin{align*}
\frakm_2 &= \fraka\oplus\frakg_{\alpha_2}\oplus\frakg_{-\alpha_2}\\
\frakn_2 &=\bigoplus_{\alpha\in\{\alpha_1,\alpha_1+\alpha_2,2\alpha_1+\alpha_2\}}\frakg_\alpha.
\end{align*}
Thus $\frakn_2 = \bigoplus_{\alpha\in\{\alpha_1,\alpha_1+\alpha_2,2\alpha_1+\alpha_2\}} \RR X_{\alpha}$. Since the root vectors $X_{\alpha}$ satisfy the commutation relations
\begin{align*}
	[X_{\alpha},X_{2\alpha_1+\alpha_2}] &=0\text{ for }\alpha \in\{\alpha_1,\alpha_1+\alpha_2\}\\
	[X_{\alpha_1}, X_{\alpha_1+\alpha_2}] &= 2X_{2\alpha_1+\alpha_2},
\end{align*}
the group $N_2$ is isomorphic to the Heisenberg group $\mathcal{H}_3$.
\subsection{Induction from the Minimal Parabolic Subgroup}
In this section consider the minimal principal series obtained by induction from the minimal parabolic subgroup $P_0$. Consider the Cartan involution $\theta(g)=(g^t)^{-1}$ on $Sp(2n,\RR)$, which acts on the split Cartan subalgebra
	\[
		\frakh_\RR = \bigoplus_{i=1}^{n}\CC(E_{i,i}-E_{n+i,n+i})
	\]
	 by sending all elements to their negatives. In this case, $\theta = -\mathrm{id}$ on the root lattice $X$. Therefore,
	 \[
	 	X/(1-\theta)X = \bigoplus_{i=1}^n\ZZ/2\ZZ \varpi_i	 \]
	 and -1 eigenspace of $\theta$ on $X$ is $X^{-\theta} = X$. The discrete character $\delta$ and the continuous character $\lambda$ can thus be determined by two vectors
	 \begin{align}
	 	\delta = \sum_{i=1}^n\delta_i\varpi_i &\longleftrightarrow (\delta_1,\ldots,\delta_n) & \delta_i\in\ZZ/2\ZZ\label{paramdelta}\\
		\lambda = \sum_{i=1}^n\lambda_i\varpi_i &\longleftrightarrow (\lambda_1,\ldots,\lambda_n) & \lambda\in\CC\label{paramlambda}.
	 \end{align}
	 If all of the $\lambda_i$'s are integers satisfying $\lambda_i\equiv\delta_i\text{ }\mathrm{mod}\text{ } 2$, then there exists an integral character $\chi_{\delta,\lambda}=\sum_{i=1}^n \lambda_i\varpi_i$ such that $\delta$ is the image of $\chi$ in the quotient $X/(1-\theta)X$, and $\lambda = \left(\frac{1-\theta}{2}\right) \chi$.\\
	 
In the case of $Sp(2n,\RR)$, a minimal principal series representation $I_{P_0}(\chi_{\delta,\lambda})$ of $Sp(4,\RR)$ is determined by the following data:
\begin{enumerate}
	\item A continuous character $\lambda: \fraka\longrightarrow \CC$ represented by the pair $(\lambda_1,\lambda_2)$ with $\lambda_i = \lambda(H_i)$. This character can be extended linearly to a character on $\fraka\otimes\CC = \frakh$, also denoted by $\lambda$;
	\item A character $\delta: M\longrightarrow \{\pm 1\}$ represented by a pair $(\delta_1,\delta_2)$ where $\delta_i\in\{0,1\}$ such that $\delta(\gamma_{\alpha_2}^{\epsilon_1}\gamma_{2\alpha_1+\alpha_2}^{\epsilon_2}) = (-1)^{\delta_1\epsilon_1+\delta_2\epsilon_2}$.
\end{enumerate}
If the numbers $\lambda_i$ are integers, with  $\delta_i\equiv \lambda_i \text{ }\mathrm{mod}\text{ } 2$, the parameters defined in (\ref{paramdelta})-(\ref{paramlambda}) combine to an algebraic character $\chi_{\delta,\lambda}$ on the complex Cartan subalgebra $\fraka_\CC = \CC H_1\oplus \CC H_2$. \\

The exponential map from $\fraka$ to the split torus $A_0$ is a bijection. Every element $a\in A_0$ can be written as $a = \exp H_a$ for some $H_a\in A_0$. For any $\nu\in X$, we introduce the notation $a^\nu = \exp \nu(H_a)$. The algebraic character $\chi_{\delta,\lambda}$ exponentiates to a character $\chi_{\delta, \lambda}:MA_0\longrightarrow \CC^*$ given by
\[
	(\chi_{\delta, \lambda+\rho})(ma) = \delta(m)a^{\lambda+\rho}.
\]
We define the \emph{principal series} representation $I_{P_0}(\delta,\lambda)$ induced from the minimal parabolic subgroup $P\subset G$ as the following vector space of functions on $G$:
\begin{align}
	I_P(\chi_{\delta, \lambda})= \{f:G\longrightarrow\CC|f(gman)=(\chi_{\delta, \lambda+\rho})^{-1}(ma)f(g)\}\label{princip}.
\end{align}
The action $\pi_{P_0}(\chi_{\delta, \lambda})$ of $G$ on $f\in I_{P_0}(\chi_{\delta, \lambda})$ is given by the \emph{left} regular representation
\[
	(\pi_{P_0}(\chi_{\delta, \lambda})(g) f)(h) = f(g^{-1}h).
\]

The $(\frakg,K)$-module of the principal series representation $I_{P_0}(\chi_{\delta,\lambda})$ embeds into the space
\[
	C_{\delta}(K) = \{f:K\longrightarrow\CC|f\text{ smooth and }f(km) = \delta(m)^{-1}f(k)\}.
\]
This space is isomorphic to the space of smooth global sections of the line bundle $K\otimes_{M}\CC_{\delta}$, where $\CC_{\delta}$ is the vector space on which $M$ acts by $\delta^{-1}$. To get a basis for the space $C_{\delta}(K)$ from the Wigner $D$-functions, we consider the right action by $m\in M$ on an arbitrary function $f\in C_\delta(K)$. $f$ can be written as the finite linear combination of Wigner $D$-functions with coefficients $a^{(j,n)}_{m_1,m_2}$:
\[
	f(k) = \sum_{\substack{j,n\\-j\leq m_i\leq j}}a^{(j,n)}_{m_1,m_2} W^{(j,n)}_{m_1,m_2}(k).
\]
The action by $m\in M$ on the right gives
\[
	f(km) = \sum_{\substack{j,n\\-j\leq m_i\leq j}}a^{(j,n)}_{m_1,m_2} W^{(j,n)}_{m_1,m_2}(k m).
\]
Recall that the Wigner $D$-functions are matrix coefficients of $U(2)$-representations. Their values on the product of two elements $km$ come from the multiplication of two matrices:
\[
	W^{(j,n)}_{m_1,m_2}(k m) = \sum_{m_3} W^{(j,n)}_{m_1,m_3}(k)W^{(j,n)}_{m_3,m_2}(m).
\]
The action of a general element \[
m = \gamma_{\alpha_2}^{\epsilon_1}\gamma_{2\alpha_1+\alpha_2}^{\epsilon_2} = e^{\pi(\epsilon_1+\epsilon_2)U_0}e^{\pi(-\epsilon_1+\epsilon_2)U_3}\in M\] on the Wigner $D$-functions is diagonal:
\[
	W^{(j,n)}_{m_3,m_2}(e^{\pi(\epsilon_1+\epsilon_2)U_0}e^{\pi(-\epsilon_1+\epsilon_2)U_3}) =(-1)^{ -(n-m_2) \epsilon_1-(n+m_2)\epsilon_2}\delta_{m_3,m_2}.
\]
Therefore, $f(km)$ can be written as
\[
	f(km) = \sum_{\substack{j,n\\-j\leq m_i\leq j}}a^{(j,n)}_{m_1,m_2}(-1)^{ -(n-m_2) \epsilon_1-(n+m_2)\epsilon_2}W^{(j,n)}_{m_1,m_2}(k).
\]
Because of the linear independence of different Wigner $D$-functions, the equality $f(km) = \delta(m)^{-1}f(k) = (-1)^{-\delta_1\epsilon_1-\delta_2\epsilon_2} f(k)$ holds for all $k\in K$ if and only if $j,n,m_1,m_2$ satisfy the compatibility condition
\[
	(-1)^{ (n-m_2) \epsilon_1+ (n+m_2)\epsilon_2}= (-1)^{\delta_1\epsilon_1+\delta_2\epsilon_2}.
\]
Therefore, the space $C_{\delta}(K)$ can be written as the direct sum
\begin{equation}
	C_{\delta}(K)= \bigoplus_{(j,n)\in\mathtt{KTypes}(\delta_1,\delta_2)}\bigoplus_{\substack{m_1\in\{-j,-j+1,\ldots,j\}\\m_2\in \mathtt{M}(j,n;\delta_1,\delta_2)}}\CC W^{(j,n)}_{m_1,m_2}
\end{equation}
in which the two sets of admissible $j,n,m_1,m_2$ are defined as:
\begin{align}
	\mathtt{KTypes}(\delta_1,\delta_2) &= \{(j,n)\in \frac{1}{2}\ZZ_{\geq 0}\times \frac{1}{2}\ZZ|2j \equiv 2n \equiv \delta_1+\delta_2 \text{ }\mathrm{mod}\text{ } 2\}\label{KCond}
\end{align}
\begin{align}
	\mathtt{M}(j,n;\delta_1,\delta_2)&= \{m_2\in\{-j,-j+1,\ldots,j\}| n-m_2\equiv \delta_1\text{ and }n+m_2\equiv \delta_2 \text{ }\mathrm{mod}\text{ } 2\}.\label{MCond}
\end{align}
Similarly to the case of $SU(2,1)$, for each $(j,n)\in\mathtt{KTypes}(\delta_1,\delta_2)$, we can denote 
\begin{align}
	\tau^{(j,n)} = \bigoplus_{\substack{m_1\in\{-j,-j+1,\ldots,j\}\\m_2\in \mathtt{M}(j,n;\delta_1,\delta_2)}}\CC W^{(j,n)}_{m_1,m_2}\label{PrincipalDecomposition2}
\end{align}
as the $K$-isotypic subspace of $I_{P_0}(\delta,\lambda)$ which decomposes into copies of irreducible $K$-representations of highest weight $(j,n)$. The restriction of the $(\frakg,K)$-module of the principal series $I_{P_0}(\delta,\lambda)$ to $K$ can be decomposed as a direct sum of the $K$-isotypic spaces $\tau^{(j,n)}$:
\begin{align}
	I_{P_0}(\delta,\lambda) = \bigoplus_{(j,n)\in\mathtt{KTypes}(\delta_1,\delta_2)}
	\tau^{(j,n)}.\label{PrincipalDecomposition1}
\end{align}
The different copies of irreducible $K$-representations are distinguished by the index $m_2$, and the action of $\fraku(2)$ raising and lowering operators $U_1\pm\ii U_2$ moves each $m_1$ to $m_1\pm 1$.  For each $K$-isotypic space $\tau^{(j,n)}$, the cardinality of the set $\mathtt{M}(j,n;\delta_1,\delta_2)$ is equal to the multiplicity of $K$-types with highest weight $(j,n)$. In fact, if we assume
\[
	(\delta_1,\delta_2) \in\{ (0,0),(1,1)\},
\]
the set $\mathtt{M}(j,n;\delta_1,\delta_2)$ is
%\[
%	\mathtt{M}(j,n;\delta_1,\delta_2) = \begin{cases}\underbrace{\{\ldots j-4,j-2,j\}}_{\lfloor j\rfloor+1\text{ elements}}&j-n+\delta_1\equiv 0\text{ }\mathrm{mod}\text{ } 2\\\underbrace{\{\ldots,j-5,j-3,j-1\}}_{\lfloor j-\frac{1}{2}\rfloor+1\text{ elements}}&j-n+\delta_1\equiv 1\text{ }\mathrm{mod}\text{ } 2\end{cases}.
%\]
\begin{align}
	\mathtt{M}(j,n;\delta_1,\delta_2) = \left\{\begin{smallmatrix}\{\ldots j-4,j-2,j\}&j-n+\delta_1\equiv 0\text{ }\mathrm{mod}\text{ } 2\\\{\ldots j-5,j-3,j-1\}&j-n+\delta_1\equiv 1\text{ }\mathrm{mod}\text{ } 2\end{smallmatrix}\right.
\end{align}
and its cardinality is
%\begin{align}
%	|\mathtt{M}(j,n;\delta_1,\delta_2)| = \left\{\begin{smallmatrix}\lfloor j\rfloor+1&j-n+\delta_1\equiv 0\text{ }\mathrm{mod}\text{ } 2\\\lfloor j-\frac{1}{2}\rfloor+1&j-n+\delta_1\equiv 1\text{ }\mathrm{mod}\text{ } 2\end{smallmatrix}\right..
%\end{align}
\begin{align}
	|\mathtt{M}(j,n;\delta_1,\delta_2)| = \left\{\begin{smallmatrix} j+1&j-n+\delta_1\equiv 0\text{ }\mathrm{mod}\text{ } 2\\\lfloor j-\frac{1}{2}\rfloor+1&j-n+\delta_1\equiv 1\text{ }\mathrm{mod}\text{ } 2\end{smallmatrix}\right..
\end{align}
We will make use of these facts to calculate the long intertwining operator for the principal series of $Sp(4,\RR)$.

%\subsection{Induction from Siegel and Heisenberg Parabolics}
%Consider a Siegel or Heisenberg parabolic subgroup $P_i = M_i A_i N_i$, with an irreducible representation $(\sigma,V_\sigma)$ on $M_i$ and a character $\chi$ on $A_i$, we look at the $(\frakg,K)$-module of the induced representation $\Ind_{P_i}^G(\sigma\otimes\chi \otimes 1)$:
%\[
%	\Ind_{P_i}^G(\sigma\otimes\chi\otimes 1):=\{f:G\longrightarrow V_{\sigma}: f(gman):=\sigma(m)^{-1}a^{-\chi-\rho}f(g)\}
%\]
%The Levi subgroups for Siegel and Heisenberg parabolic subgroups are either $SL(2,\RR)$ or $GL(2,\RR)$. According to the classification theorem of irreducible $(\frakg,K)$-modules of $SL(2,\RR)$ and $GL(2,\RR)$, we can take $\sigma$ to be the following possibilities:
%\begin{enumerate}
%	\item The case of Heisenberg parabolic, where $M_2\cong SL(2,\RR)\times \RR^{\times}$, with a Lie algebra generated by:
%	\[
%		\frakm_1 = \underbrace{\CC(U_0-U_3)\oplus \CC(-\sqrt{2}\ii u_{\alpha_2})\oplus \CC(-\sqrt{2}\ii u_{-\alpha_2})}_{\fraksl(2,\RR)}\oplus \CC H_{\alpha_1}.
%	\]
%	Using the following classification of $SL(2,\RR)$-representations:
%		\begin{enumerate}
%			\item Principal series
%			\item Discrete series
%		\end{enumerate}
%	\item	The case of Siegel parabolic, where $M_1\cong GL(2,\RR)$:\\
%		The Lie algebra of $\frakm_1$ have the following decomposition $\CC U_2\oplus \CC$
%		\begin{enumerate}
%			\item Principal series
%			\item Discrete series
%		\end{enumerate}
%\end{enumerate}

\section{The $(\frakg,K)$-Module Structure}
\subsection{Normalization of basis and the Iwasawa decomposition}
The set of positive and negative noncompact roots are \[\begin{matrix}\Delta^+_{nc} = \{\beta_2,\beta_1+\beta_2,2\beta_1+\beta_2\}& \Delta^-_{nc} = \{-\beta_2,-\beta_1-\beta_2,-2\beta_1-\beta_2\}\end{matrix}.\]
We define the vectors $u_{\beta}\in \frakp_\CC$ by multiplying the $v_\beta$ defined in (\ref{V1})-(\ref{V3}) by a factor:
\[
	u_\beta = \begin{cases}\sqrt{2}\ii v_\beta& \beta\in\{\beta_2,2\beta_1+\beta_2,-\beta_2,-2\beta_1-\beta_2\}\\
	v_\beta &\text{otherwise}
	\end{cases}.
\]
% is a root of the root system $\Delta(\frakg_\CC,\frakt_\CC)$, and equal to 0 otherwise. In fact, we have:
%\[
%	u_\beta = \begin{cases}\sqrt{2}\ii v_\beta& \beta\in\{\alpha_2,2\alpha_1+\alpha_2,-\alpha_2,-2\alpha_1-\alpha_2\}\\
%	v_\beta &\text{otherwise}
%	\end{cases}.
%\]
Under this normalization, the Lie algebra Iwasawa decomposition of the basis vectors $u_\beta$ is:
\begin{align}
	u_{\pm(2\beta_1+\beta_2)}&=\frac{1}{\sqrt{2}}(\pm(U_0+U_3)+\ii H_1\mp 2X_{2\alpha_1+\alpha_2})\label{IwasawaSp4:1}\\
	u_{\pm\beta_2}&=\frac{1}{\sqrt{2}}( \pm(U_0-U_3)+\ii H_2 \mp 2X_{\alpha_2})\label{IwasawaSp4:2}\\
	u_{\pm(\beta_1+\beta_2)}&= -( U_1\mp \ii U_2)+ X_{\alpha_1+\alpha_2} \mp \ii X_{\alpha_1}\label{IwasawaSp4:3}
	%u_{-2\alpha_1-\alpha_2}&=\frac{\sqrt{2}}{2}(-(U_0+U_3) + \ii H_1 + 2X_{2\alpha_1+\alpha_2})\\
	%u_{-\alpha_2}&=\frac{\sqrt{2}}{2}(-(U_0-U_3) + \ii H_2 +2 X_{\alpha_2})\\
	%u_{-\alpha_1-\alpha_2}&=- (U_1+i U_2)+ X_{\alpha_1+\alpha_2} + iX_{\alpha_1}
\end{align}
The representation of $K=U(2)$ on $\frakp_\CC$ decomposes into two irreducible subrepresentations
\[
	\begin{matrix}\frakp_\CC^+ = \oplus_{\beta\in\Delta^+_{nc}}\CC u_\beta & \frakp_\CC^- = \oplus_{\beta\in\Delta^-_{nc}}\CC u_\beta\end{matrix}.
\]
The $u_\beta$'s are the normalized weight vectors on which the adjoint action by $U_1\pm \ii U_2$ acts as
\[
	\ad(U_1\pm\ii U_2) u_{\beta} = q_{\beta,\beta\pm\beta_1} u_{\beta\pm\beta_1}, 
\]
where the coefficient $q_{\beta,\beta\pm\beta_1}$ turns out to be
\[
	q_{\beta,\beta\pm\beta_1}=\begin{cases}-\ii\sqrt{2} &\text{ if }\beta\pm\beta_1\in \Delta(\frakg_\CC,\frakt_\CC)\\
	0&\text{ otherwise}\end{cases}.
\]

\subsection{Left Action of $\frakp_\CC$}
For any noncompact root $\beta$, define two integers:
\begin{align}
	m_\beta  &= -\ii \beta(U_3) \in\{-1,0,1\}\label{Sp4M}\\
	n_\beta &=  -\ii \beta(U_0) \in\{-1,1\}\label{Sp4N}
\end{align}
Each weight $\beta$ is uniquely determined by the pair of integers $(m_\beta,n_\beta)$ specified in the following chart:\\
\begin{table}[h]
\begin{center}
\begin{tabular}{cccc}
	 & $m_{\beta} = -1$ & $m_{\beta}=0$ & $m_{\beta} = 1$\\\hline
	$n_{\beta} = -1$ & $-2\beta_1-\beta_2$ & $-\beta_1-\beta_2$ & $-\beta_2$\\
	$n_{\beta} = 1$ & $\beta_2$ & $\beta_1+\beta_2$ & $2\beta_1+\beta_2$
\end{tabular}
\end{center}
\caption{Correspondence between $(m_\beta,n_\beta)$ and noncompact roots}
\label{CorrespondenceRoots}
\end{table}\\
%Therefore, for any weight vector $v^{(j,n)}_{m}$ in the irreducible representation $V^{(j,n)}$ of $U(2)$, the weight $(U_1\pm \ii U_2)v^{(j,n)}_{m}$ is 
% For any irreducible representation $V^\mu$, consider a weight vector $v^\mu_\nu$, the weight of  $(U_1\pm \ii U_2)v^\mu_\nu$ is $\nu\pm\beta_1$.\\
By embedding the principal series representation $I(\chi_{\delta,\lambda})$ into $C_{\delta}(K)$ (compare to the identical method in \cite{MillerButtcane} and the $SU(2,1)$ case in Section 5), we would like to understand the action of
\[
	\dd l(u_{\beta}) = \dd r(-\Ad(k^{-1}) u_\beta)
\]
on any basis vector $W^{(j,n)}_{m_1,m_2}(k)$ in $C_{\delta}(K)$. By the definition of the principal series in (\ref{princip}), the right action on any vector $f\in I(\chi_{\delta,\lambda})$ by $H_i$ is always a scalar multiplication by $\lambda_i+ \rho(H_i)$
\begin{align}
	\dd r(H_i) f=-(\lambda_i+ \rho(H_i))f,\label{Sp4CartanAction}
\end{align}
and the right action by any element from $N$ annihilates $f$. We can therefore use the decomposition (\ref{IwasawaSp4:1})-(\ref{IwasawaSp4:3}) and the differential operators (\ref{CompactAction:1})-(\ref{CompactAction:3}) in $\frakg_\CC$ of $u_{\beta}$ to calculate the right action of $u_{\beta}$ on the basis vectors $W^{(j,n)}_{m_1,m_2}$ of $I(\chi_{\delta,\lambda})$. Also, recall from 6.1.1 that $\rho(H_1)=2, \rho(H_2)=1$, we have
\begin{align}
	\dd r (u_{\pm(2\beta_1+\beta_2)})W^{(j,n)}_{m_1,m_2}&=\frac{\ii}{\sqrt{2}}(\mp n\mp m_2- (\lambda_1+\rho(H_1)))W^{(j,n)}_{m_1,m_2}\nonumber\\
	&=\frac{\ii(\mp n\mp m_2- (\lambda_1+2))}{\sqrt{2}}W^{(j,n)}_{m_1,m_2}\label{Sp4RightAction:1}\\
	\dd r(u_{\pm\beta_2})W^{(j,n)}_{m_1,m_2}&=\frac{\ii}{\sqrt{2}}(\mp n\pm m_2- (\lambda_2+\rho(H_2)))W^{(j,n)}_{m_1,m_2}\nonumber\\
	&=\frac{\ii(\mp n\pm m_2- (\lambda_2+1))}{\sqrt{2}}W^{(j,n)}_{m_1,m_2}\label{Sp4RightAction:2}\\
	\dd r(u_{\pm(\beta_1+\beta_2}))W^{(j,n)}_{m_1,m_2}
	&= -\ii \sqrt{(j\pm m_2)(j\mp m_2+1)}W^{(j,n)}_{m_1,m_2\mp 1}.\label{Sp4RightAction:3}
%	\dd r(u_{-2\beta_1-\beta_2})W^{(j,n)}_{m_1,m_2}
%	&=\frac{\ii}{\sqrt{2}}(-n-m_1- \langle \check{\alpha}_1, \lambda+\rho_0 \rangle)W^{(j,n)}_{m_1,m_2}\\
%	\dd r(u_{-\beta_2})W^{(j,n)}_{m_1,m_2}
%	&=\frac{\ii}{\sqrt{2}}(-n+m_1- \langle \check{\alpha}_2, \lambda+\rho_0 \rangle)W^{(j,n)}_{m_1,m_2}\\
%	\dd r(u_{-\beta_1-\beta_2})W^{(j,n)}_{m_1,m_2}
%	&= \ii \sqrt{(j-m_1)(j+m_1+1)}W^{(j,n)}_{m_1+1,m_2}
\end{align}
Recalling the correspondence between the weight $\beta$ and $(m_\beta,n_\beta)$ discussed in 
(\ref{Sp4M})-(\ref{Sp4N}) and Table $\ref{CorrespondenceRoots}$, since the irreducible constituents of $\frakp^\pm_\CC$ have highest weights $(j,n) = (1,\pm 1)$ respectively, the right action of $u_\beta$ with $\beta\in\Delta_{nc}^\pm$ can be transferred to the left by observing
\begin{align*}
	\Ad(k^{-1}) u_{\beta} &= \sum_{\nu\in\Delta_{nc}^\pm} \overline{W^{(1,\pm 1)}_{m_\nu,m_{\beta}}(k^{-1})}u_{\nu}= \sum_{\nu\in\Delta_{nc}^\pm} W^{(1,\pm 1)}_{m_\beta,m_{\nu}}(k)u_{\nu}.
\end{align*}
Based on the correspondence in Table \ref{CorrespondenceRoots} between the weights of $u_{\beta}$ and the pair of integers $(m_\beta,n_\beta)$, it is clear that the left action of $u_\beta$  for $\beta\in \Delta^\pm _{nc}$ on the Wigner $D$-functions $W^{(j,n)}_{m_1,m_2}\in C_{\delta}(K)$ can be written explicitly as follows,
\begin{align}
	&\dd l(u_{\beta})W^{(j,n)}_{m_1,m_2} = \nonumber\\
%	\begin{cases}
	&\left(W^{(1,\pm 1)}_{m_\beta,\mp 1} \dd r(u_{\pm \beta_2})+ W^{(1,\pm 1)}_{m_\beta,0} \dd r(u_{\pm(\beta_1+\beta_2)})
	 + W^{(1,\pm 1)}_{m_\beta,\pm 1} \dd r(u_{\pm (2\beta_1+\beta_2)})\right)W^{(j,n)}_{m_1,m_2}.
%	&= (-1)^{m_1}\ii\sum_{|j-1|\leq J \leq j+1} \frac{1}{\sqrt{2}}(n-m_1- \langle \check{\alpha}_2, \lambda+\rho_0 \rangle)\left(\begin{smallmatrix}J,m_1+1\\j,m_1,1,1\end{smallmatrix}\right)\left(\begin{smallmatrix}J,m_2-m_\beta\\j,m_2,1,-m_\beta\end{smallmatrix}\right)W^{(J,n-1)}_{m_1+1,m_2-m{(\nu)}}\\
%	&+(- \sqrt{(j+m_1)(j-m_1+1)} \left(\begin{smallmatrix}J,m_1-1\\j,m_1-1,1,0\end{smallmatrix}\right)+ \frac{1}{\sqrt{2}}(n+m_1- \langle \check{\alpha}_1, \lambda+\rho_0 \rangle)\left(\begin{smallmatrix}J,m_1-1\\j,m_1,1,-1\end{smallmatrix}\right)) \left(\begin{smallmatrix}J,m_1-m_\beta\\j,m_2,1,-m_\beta\end{smallmatrix}\right)W^{(J,n-1)}_{m_1-1,m_2-m{(\nu)}}\\
%	(-1)^{m_1} (W^{(1,-1)}_{1, -m_\beta} \dd r(u_{-2\beta_1-\beta_2}) - W^{(1,-1)}_{0,-m_\beta} \dd r(u_{-\beta_1-\beta_2}) \\
%	+ W^{(1,-1)}_{-1, -m_\beta} \dd r(u_{-\beta_2}))W^{(j,n)}_{m_1,m_2} & \nu\in \Delta^-_{nc}\\
%	\end{cases}
%	&= (-1)^{m_1}\ii\sum_{|j-1|\leq J \leq j+1}\\
%	& (- \sqrt{(j-m_1)(j+m_1+1)}\left(\begin{smallmatrix}J,m_1+1\\j,m_1+1,1,0\end{smallmatrix}\right)+ \frac{1}{\sqrt{2}}(-n-m_1- \langle \check{\alpha}_1, \lambda+\rho_0 \rangle)\left(\begin{smallmatrix}J,m_1+1\\j,m_1,1,1\end{smallmatrix}\right)) \left(\begin{smallmatrix}J,m_1-m_\beta\\j,m_2,1,-m_\beta\end{smallmatrix}\right)W^{(J,n+1)}_{m_1+1,m_2-m{(\nu)}} \\
%	&+\frac{1}{\sqrt{2}}(-n+m_1- \langle \check{\alpha}_2, \lambda+\rho_0 \rangle)\left(\begin{smallmatrix}J,m_1-1\\j,m_1,1,-1\end{smallmatrix}\right)\left(\begin{smallmatrix}J,m_2-m_\beta\\j,m_2,1,-m_\beta\end{smallmatrix}\right)W^{(J,n+1)}_{m_1-1,m_2-m{(\nu)}}\\
\end{align}
According to (\ref{Sp4RightAction:1})-(\ref{Sp4RightAction:3}), we apply the right action $\dd r(u_\beta)$ to 
$W^{(j,n)}_{m_1,m_2}$ and get
\begin{align}
	&\dd l(u_{\beta})W^{(j,n)}_{m_1,m_2} = \nonumber\\
%	\begin{cases}
	&\ii\left(\frac{\mp n\pm m_2- (\lambda_2+\rho(H_2))}{\sqrt{2}}W^{(1,\pm 1)}_{m_\beta,\mp 1}
	 + \frac{\mp n\mp m_2- (\lambda_1+\rho(H_1))}{\sqrt{2}}W^{(1,\pm 1)}_{m_\beta,\pm 1}\right)W^{(j,n)}_{m_1,m_2}\nonumber\\
	 &-\ii  \sqrt{(j\pm m_2)(j\mp m_2+1)}W^{(1,\pm 1)}_{m_\beta,0} W^{(j,n)}_{m_1,m_2\mp 1}.
\end{align}

%\begin{align*}
%	&\dd l(u_{\nu})W^{(j,n)}_{m_1,m_2} = (-1)^{m_1} (W^{(1,-1)}_{1, -m_\beta} \dd r(u_{\alpha_2}) - W^{(1,-1)}_{0,-m_\beta} \dd r(u_{\alpha_1+\alpha_2}) + W^{(1,-1)}_{-1, -m_\beta} \dd r(u_{2\alpha_1+\alpha_2}))W^{(j,n)}_{m_1,m_2}\\
%	&= (-1)^{m_1} (\frac{\ii}{\sqrt{2}}(n-m_1- \langle \check{\alpha}_2, \lambda+\rho_0 \rangle)W^{(1,-1)}_{1, -m_\beta}W^{(j,n)}_{m_1,m_2}\\
%	&-  \ii \sqrt{(j+m_1)(j-m_1+1)} W^{(1,-1)}_{0,-m_\beta}W^{(j,n)}_{m_1-1,m_2}+ \frac{\ii}{\sqrt{2}}(n+m_1- \langle \check{\alpha}_1, \lambda+\rho_0 \rangle)W^{(1,-1)}_{-1, -m_\beta} W^{(j,n)}_{m_1,m_2})
%\end{align*}
We can replace the products of Wigner $D$-functions by a linear combination of Wigner $D$-functions with Clebsch-Gordan coefficients as described in formula (\ref{ProductClebschGordan}) of Section \ref{cptgrp}. The left action of $u_\beta$ for $\beta\in \Delta^{\pm}_{nc}$ on Wigner $D$-functions can thus be expressed as
\begin{align}
	&\dd l(u_{\beta})W^{(j,n)}_{m_1,m_2} \nonumber\\
	=&
%	&\begin{cases} 
	\ii\sum_{j_0\in\{j-1,j,j+1\}}\left(\begin{smallmatrix}j+j_0,m_1+m_\beta\\j,m_1,1,m_\beta\end{smallmatrix}\right)\left(\frac{\mp n\mp m_2- (\lambda_1+\rho(H_1))}{\sqrt{2}}\left(\begin{smallmatrix}j+j_0,m_2\pm 1\\j,m_2,1,\pm 1\end{smallmatrix}\right)W^{(j+j_0,n\pm 1)}_{m_1+m_\beta,m_2\pm 1}+\right.\nonumber\\
	&\left.\left(- \sqrt{(j\pm m_2)(j\mp m_2+1)} \left(\begin{smallmatrix}j+j_0,m_2\mp 1\\j,m_2\mp 1,1,0\end{smallmatrix}\right)+ \frac{\mp n\pm m_2- (\lambda_2+\rho(H_2))}{\sqrt{2}}\left(\begin{smallmatrix}j+j_0,m_2\mp 1\\j,m_2,1,\mp 1\end{smallmatrix}\right)\right)\right.\nonumber\\
	&\left.W^{(j+j_0,n\pm 1)}_{m_1+m_\beta,m_2\mp 1}\right)
\end{align}
After computing all the Clebsch-Gordan coefficients using the formulas listed in Table \ref{TableCG2}, the action of weight vectors $u_\beta$ of $\frakp_\CC$ on the left when $\beta\in \Delta^\pm_{nc}$ can be expressed as the following linear combination:
\begin{align}
	\dd l(u_{\beta})W^{(j,n)}_{m_1,m_2} =& \frac{\ii}{2}\sum_{\substack{j_0\in\{-1,0,1\}\\\epsilon=\pm 1}}\begin{spmatrix}j+j_0,m_1+m_\beta\\j,m_1,1,m_\beta\end{spmatrix}C_{j+j_0}q_{j_0,\varepsilon}\kappa_{\pm,j_0,\varepsilon}(j,n,m_1;\lambda) W^{(j+j_0,n\pm 1)}_{m_1+m_\beta,m_2+\varepsilon}\label{Sp4GKAction}
\end{align}
with the coefficients given in the tables:\\
\begin{center}
	\begin{tabular}{cc}
	$C_{j+j_0}$ & \\\hline
		$j_0 = -1$ & $j^{-\frac{1}{2}}(2j+1)^{-\frac{1}{2}}$\\
		$j_0=0$ & $j^{-\frac{1}{2}}(j+1)^{-\frac{1}{2}} $\\
		$j_0=1$ & $(j+1)^{-\frac{1}{2}}(2j+1)^{-\frac{1}{2}}$
\end{tabular}\\

  \begin{tabular}{ccc}
		$q_{ j_0,\varepsilon}$& $\varepsilon = -1$ & $\varepsilon = 1$\\\hline
		$j_0 = -1$ & $\sqrt{(j+m_2-1) (j+m_2)}$ & $\sqrt{(j-m_2-1) (j-m_2)}$ \\
		$j_0=0$ & $\sqrt{(j+m_2)(j-m_2+1)}$ & $\sqrt{(j-m_2)(j+m_2+1)}$\\
		$j_0=1$ &$\sqrt{(j-m_2+1)(j-m_2+2)}$ & $\sqrt{(j+m_2+1)(j+m_2+2)}$
\end{tabular}\\

  \begin{tabular}{ccc}
		 $\kappa_{+, j_0,\varepsilon}$ & $\varepsilon = -1$ & $\varepsilon = 1$\\\hline
		$j_0 = -1$ & $2+2j -m_2-n-\left( \lambda_2+\rho(H_2)\right) $ & $-n-m_2-\left( \lambda_1+\rho(H_1)\right)$ \\
		$j_0=0$ & $2-n-m_2-\left( \lambda_2+\rho(H_2)\right)$ & $n+m_2+\left( \lambda_1+\rho(H_1)\right)$\\
		$j_0=1$ &$-2j-m_2-n-\left( \lambda_2+\rho(H_2)\right)$ & $-n-m_2- \left( \lambda_1+\rho(H_1)\right)$\\
\end{tabular}\\

 \begin{tabular}{ccc}
		$\kappa_{-, j_0,\varepsilon}$ & $\varepsilon = -1$ & $\varepsilon = 1$\\\hline
		$j_0 = -1$ & $n+m_2-\left( \lambda_1+\rho(H_1)\right)$ & $2+2j+m_2+n- \left( \lambda_2+\rho(H_2)\right)$ \\
		$j_0=0$ & $n+m_2-\left( \lambda_1+\rho(H_1)\right)$ & $-2-n-m_2+ \left( \lambda_2+\rho(H_2)\right)$\\
		$j_0=1$ &$n+m_2-\left( \lambda_1+\rho(H_1)\right)$ & $-2j+m_2+n-\left( \lambda_2+\rho(H_2)\right)$\\
\end{tabular}.

\end{center}
The $(\frakg,K)$-action on $Sp(4,\RR)$ principal series $I(\chi_{\delta,\lambda})$ is completely determined by formula (\ref{Sp4GKAction}) and the four tables above.
%\section{Reducibility of $I(P,\chi_{\delta_1,\delta_2,\lambda})$}

\section{Intertwining Operators}
The longest element $w_0=w_{\alpha_2}w_{\alpha_1}w_{\alpha_2}w_{\alpha_1}$ in the Weyl group of $Sp(4,\RR)$  corresponds to the long intertwining operator:
\[
	A(w_0,\chi_{\delta,\lambda})f(k) = \int_{\overline{N}\cap w^{-1} N w}f(k w \overline{n})\dd \overline{n}.
\]
Applying the Langlands' Lemma  on factorization of intertwining operators in \cite{ShahidiEisenstein}, $A(w_0,\chi_{\delta,\lambda})$ can be factored into 4 intertwining operators corresponding to simple reflections:
\begin{align}
	A(w_0,\chi_{\delta,\lambda}) &= A(w_{\alpha_2},w_{\alpha_1}w_{\alpha_2}w_{\alpha_1}\chi_{\delta,\lambda})A(w_{\alpha_1},w_{\alpha_2}w_{\alpha_1}\chi_{\delta,\lambda})A(w_{\alpha_2},w_{\alpha_1}\chi_{\delta,\lambda})\times\nonumber \\
	&A(w_{\alpha_1},\chi_{\delta,\lambda}).\label{langlandslemma}
\end{align}
Since there is an embedding of the Harish-Chandra module of $I(\chi_{\delta,\lambda})$ into the space $C_\delta(K)$, we can express any function $f(k)$ in the principal series under such embedding as a linear combination of Wigner $D$-functions:
\[
	f(k) = \sum_{(j,n)}\sum_{(m_1,m_2)}a^{(j,n)}_{m_1,m_2}W^{(j,n)}_{m_1,m_2}(k).
\]
Therefore, it is sufficient to compute the matrix coefficients of the intertwining operator on the basis $W^{(j,n)}_{m_1,m_2}$ of the space $C_\delta(K)$:
\begin{align*}
	A(w,\chi_{\delta,\lambda})W^{(j,n)}_{m_1,m_2}(k) &= \int_{\bar{N}\cap N^{w}}W^{(j,n)}_{m_1,m_2}(kw\bar n)\dd\bar{n}\\
    &= \int_{\bar{N}\cap N^{w}}a(\bar n)^{-(\lambda+\rho)}W^{(j,n)}_{m_1,m_2}(k w k(\bar n))\dd\bar{n}\\
    &=\sum_{m_3} \left(\int_{\bar{N}\cap N^{w}}a(\bar n)^{-(\lambda+\rho)}W^{(j,n)}_{m_3,m_2}( wk(\bar n))\dd\bar{n}\right)W^{(j,n)}_{m_1,m_3}(k),
\end{align*}
and the matrix coefficients of the intertwining operator $A(w,\delta,\lambda)$ are given by the formula
\begin{equation}\label{ComputeIntertwining}
	\langle W^{(j,n)}_{m_1,m_3},A(w,\chi_{\delta,\lambda})W^{(j,n)}_{m_1,m_2}\rangle = \int_{\bar{N}\cap N^{w}}a(\bar n)^{-(\lambda+\rho)}W^{(j,n)}_{m_3,m_2}( wk(\bar n))\dd\bar{n}.
\end{equation}
We will use this formula to compute the intertwining operators explicitly for the group $Sp(4,\RR)$ in this section. Combining the Langlands' Lemma and (\ref{ComputeIntertwining}), we have the following proposition:

\begin{proposition}\label{GindikinKarpelevich}
Let $w_0=w_{\alpha_2}w_{\alpha_1}w_{\alpha_2}w_{\alpha_1}$ be the longest element in the Weyl group $W$ of $Sp(4,\RR)$. The matrix for the long intertwining operator
$A(w_0,\chi_{\delta,\lambda})$ under the basis $W^{(j,n)}_{m_1,m_2}$ of $C_{\delta}(K)$ has a factorization:
\begin{align*}
	A(w_0,\chi_{\delta,\lambda})|_{j,n} = A_4(\nu)\cdot A_3(\lambda)\cdot A_2(\lambda)\cdot A_1(\lambda)
\end{align*}
If we define 
\begin{equation}\label{QDef}
	Q(z,n) = \frac{\pi 2^{2-2z}\Gamma(2z-1)}{\Gamma(z+n)\Gamma(z-n)},
\end{equation}
and let
\begin{align*}
	S^{j,n}_{m_3,m_2}(z)&=\sum_{-j\leq m_4\leq j}\ii^{-2m_4}M^{j,n}_{m_3,m_4}N^{j,n}_{m_4,m_2}Q(z,m_4)\\
   T^n_{m_1}(z)&=\ii^{-n+m_1}Q\left(z,\frac{m_1-n}{2}\right)
\end{align*}
%\substack{m_3\in \mathtt{M}(j,n;\delta_1,\delta_2)}
where $N^{j,n}_{m_1,m_3}$ is the inverse matrix of  $M^{j,n}_{m_1,m_3}$, with
\begin{align*}
	M^{j,n}_{m_3,m_4}
    =&c^j_{m_3}c^j_{m_4}\begin{cases}
	\frac{\ii^{m_3-m_4}(-1)^{2j}2^{-j}}{(j-m_3)!(m_3-m_4)!(j+m_4)!}{}_2F_1(\begin{smallmatrix}-j+m_3,-j-m_4\\1+m_3-m_4\end{smallmatrix};-1)&m_3>m_4\\
    \frac{\ii^{m_4-m_3}(-1)^{2j}2^{-j}}{(j+m_3)!(m_4-m_3)!(j-m_4)!}{}_2F_1(\begin{smallmatrix}-j-m_3,-j+m_4\\1-m_3+m_4\end{smallmatrix},-1)&m_3\leq m_4
\end{cases}
\end{align*}
then the operators $A_i(\lambda)$ act as
\begin{align}
	A_1(\lambda)W^{(j,n)}_{m_1,m_2}&=\sum_{\substack{m_3\in \mathtt{M}(j,n;\delta_2,\delta_1)}}S^{j,n}_{m_3,m_2}\left(\frac{\lambda_1-\lambda_2+1}{2}\right)W^{(j,n)}_{m_1,m_3}\label{A1}\\
   A_2(\lambda)W^{(j,n)}_{m_1,m_2}&=T^n_{m_2}\left(\frac{\lambda_1+1}{2}\right)W^{(j,n)}_{m_1,m_2}\label{A2}\\
   A_3(\lambda)W^{(j,n)}_{m_1,m_2}&=\sum_{\substack{m_3\in \mathtt{M}(j,n;\delta_2,\delta_1)}}S^{j,n}_{m_3,m_2}\left(\frac{\lambda_1+\lambda_2+1}{2}\right)W^{(j,n)}_{m_1,m_3}\label{A3}\\
 A_4(\lambda)W^{(j,n)}_{m_1,m_2}&=T^n_{m_2}\left(\frac{\lambda_2+1}{2}\right)W^{(j,n)}_{m_1,m_2}.\label{A4}
\end{align}
\end{proposition}

We will prove the Proposition \ref{GindikinKarpelevich} in the following two sections \ref{ProofOfProp1} and \ref{ProofOfProp2}. \begin{remark}The simple reflection $w_{\alpha_1}$ sends the character $\delta = (\delta_1,\delta_2)$ on $M$ to $(\delta_2,\delta_1)$. For the $w_{\alpha_1}$ intertwining operator, it is important to recall that the parity condition $\mathtt{M}(j,n;\delta_2,\delta_1)$ of $m_3$ restricts the allowed Wigner $D$-functions $W^{(j,n)}_{m_1,m_3}$ in $I_{P_0}(w_{\alpha_1}\chi_{\delta,\lambda})$ as in (\ref{PrincipalDecomposition2}) and (\ref{PrincipalDecomposition1}).\end{remark}

\subsection{Rank 1 Intertwining Operators}\label{ProofOfProp1}
Starting from any character $\lambda$ on $\fraka_\CC$, the rank 1 intertwining operator $A(w_{\alpha}, \mu)$ associated to a simple reflection $w_{\alpha}$ can be written as
\begin{align*}
	A(w_{\alpha},\lambda)f(k) &= \int_{\overline{N}\cap w_{\alpha}^{-1} N w_{\alpha}} f(kw_{\alpha} \bar{n})\dd \overline{n}= \int_{-\infty}^{\infty} f(kw_{\alpha} \exp(t X_{-\alpha}))\dd t.
\end{align*}
By the Iwasawa decomposition of $\exp(t X_{-\alpha})$ given in (\ref{IwasawaGroup}), \[\exp(t X_{-\alpha}) = \kappa_{\alpha}(t)h_{\alpha}(\sqrt{1+t^2}) \chi_{\alpha}\left(\frac{t}{1+t^2}\right),\]
where $\kappa_{\alpha}(t) = \exp \left(\arctan(-t)(X_{\alpha}-X_{-\alpha})\right)$.
If $f$ is any vector in the principal series representation $I(\chi_{\delta,\lambda})$, the action of $\exp(t X_{-\alpha})$ on the right is
\begin{align*}
	f(kw_{\alpha}\exp(t X_{-\alpha})) &= f\left(kw_{\alpha}\kappa_{\alpha}(t)h_{\alpha}(\sqrt{1+t^2}) \chi_{\alpha}\left(\frac{t}{1+t^2}\right)\right)\\
	&=(1+t^2)^{-\frac{\langle \check{\alpha},\lambda+\rho\rangle}{2}}f\left(kw_{\alpha}e^{\arctan(-t)(X_{\alpha}-X_{-\alpha})}\right).
\end{align*}
Denoting by $\theta(t) = \arctan(-t)$, since $X_{\alpha_1}-X_{-\alpha_1}=2 U_2, X_{\alpha_2}-X_{-\alpha_2}=U_0-U_3$, and recalling from (\ref{Sp4Weyl}) the expressions of simple reflections $w_{\alpha_i}$ in terms of Euler angles, the action of $\exp(t X_{-\alpha_1})$ and $\exp(t X_{-\alpha_2})$ on a Wigner $D$-function $f = W^{(j,n)}_{m_1,m_2}$ is thus
\begin{align}
	&W^{(j,n)}_{m_1,m_2}(kw_{\alpha_1}e^{\theta(t)(X_{\alpha_1}-X_{-\alpha_1})})\nonumber\\
	=&
	\sum_{-j \leq m_3\leq j}W^{(j,n)}_{m_1,m_3}(k)W^{(j,n)}_{m_3,m_2}(0,0,-\pi-2\theta(t),0)\label{TargetWigner1}\\
	&W^{(j,n)}_{m_1,m_2}(kw_{\alpha_2}e^{\theta(t)(X_{\alpha_2}-X_{-\alpha_2})})\nonumber\\
	=&
	\sum_{-j \leq m_3\leq j}W^{(j,n)}_{m_1,m_3}(k)W^{(j,n)}_{m_3,m_2}\left(-\frac{\pi}{2}-\theta(t),\frac{\pi}{2}+\theta(t),0,0\right)\label{TargetWigner2}.
\end{align}
Therefore, the simple intertwining operators $A(w_{\alpha_1},\lambda)W^{(j,n)}_{m_1,m_2}$ and $A(w_{\alpha_2},\lambda)W^{(j,n)}_{m_1,m_2}$ can be expressed in terms of integrals involving Wigner $D$-functions:
\begin{align}
	&\left(A(w_{\alpha_1},\lambda)W^{(j,n)}_{m_1,m_2} \right)(k)=\nonumber\\
	&\sum_{-j \leq m_3\leq j} W^{(j,n)}_{m_1,m_3}(k)\int_{-\infty}^\infty (1+t^2)^{-\frac{\langle \check{\alpha_1},\lambda+\rho\rangle}{2}}W^{(j,n)}_{m_3,m_2}(0,0,-\pi-2\theta(t),0)\dd t\label{NewA1}\\
	&\left(A(w_{\alpha_2},\lambda)W^{(j,n)}_{m_1,m_2} \right)(k)=\nonumber\\
	&\sum_{-j \leq m_3\leq j} W^{(j,n)}_{m_1,m_3}(k)\int_{-\infty}^\infty (1+t^2)^{-\frac{\langle \check{\alpha_2},\lambda+\rho\rangle}{2}}W^{(j,n)}_{m_3,m_2}\left(-\frac{\pi}{2}-\theta(t),\frac{\pi}{2}+\theta(t),0,0\right)\dd t\label{NewA2}
\end{align}
Observe that $\theta(t)$ is an odd function. It is important to mention that, from the formula of the Wigner $D$-function (\ref{DPart}), we have
\begin{align*}
d^{(j,n)}_{m_3,m_2}(-\pi-2\theta(t))&=\sum_{p=\max(0,m_3-m_2)}^{\min(j-m_2,j+m_3)}\frac{(-1)^{2j+m_2-m_3+p}}{(j+m_3-p)!p!(m_2-m_3+p)!(j-m_2-p)!}\nonumber\\&\cos^{m_2-m_3+2p}\left(\theta(t)\right)\sin^{2 j + m_3 - m_2 - 2 p}\left(\theta(t)\right).
\end{align*}
If $2j+m_3-m_2\equiv 1\text{ mod } 2$, the integrand of (\ref{NewA1}) is an odd function, which makes the integral (\ref{NewA1}) zero. If $W^{(j,n)}_{m_1,m_2}\in I_{P_0}(\chi_{\delta,\lambda})$, we must have
\[
	2j\equiv\delta_1+\delta_2\text{ mod }2
\]
and
\[
	m_2\in \mathtt{M}(j,n;\delta_1,\delta_2).
\]
The set $\mathtt{M}(j,n;\delta_1,\delta_2)$ has been defined in (\ref{MCond}). In order to make the integral (\ref{NewA1}) nonzero, the function $d^{(j,n)}_{m_3,m_2}(-\pi-2\theta(t))$ must be an even function. In this case, the exponent $2j+m_3-m_2 \equiv 0\text{ mod } 2$. Therefore,
\begin{align}
	2j+m_3-m_2\equiv \delta_1+\delta_2+m_3-m_2 \equiv 0\text{ mod } 2.\label{m3m2rel}
\end{align}
Thus, if $m_2$ satisfies $n-m_2\equiv \delta_1\text{ mod } 2$ and $n+m_2\equiv \delta_2\text{ mod } 2$, by (\ref{m3m2rel}), we must have
\begin{align*}
	n-m_3\equiv n-m_2+\delta_1+\delta_2&\equiv \delta_2\text{ mod } 2\\
	n+m_3\equiv n+m_2 +\delta_1+\delta_2&\equiv \delta_1\text{ mod } 2.
\end{align*}
Thus the parity condition $m_3$ is given by the set  $\mathtt{M}(j,n;\delta_2,\delta_1)$, with $\delta_1$ and $\delta_2$ flipped from the parity condition of $m_2$. In particular, $A(w_{\alpha_1},\lambda)W^{(j,n)}_{m_1,m_2}$ lies in the space $I_{P_0}(w_{\alpha_1}\chi_{\delta,\lambda})$. Hence the sum in (\ref{NewA1}) is in fact a sum over $m_3\in \mathtt{M}(j,n;\delta_2,\delta_1)$:
\begin{align}
	&\left(A(w_{\alpha_1},\lambda)W^{(j,n)}_{m_1,m_2} \right)(k)=\nonumber\\
	&\sum_{m_3\in \mathtt{M}(j,n;\delta_2,\delta_1)} W^{(j,n)}_{m_1,m_3}(k)\int_{-\infty}^\infty (1+t^2)^{-\frac{\langle \check{\alpha_1},\lambda+\rho\rangle}{2}}W^{(j,n)}_{m_3,m_2}(0,0,-\pi-2\theta(t),0)\dd t.\label{NewNewA1}
\end{align}
Note that the Wigner $D$-function $W^{(j,n)}_{m_3,m_2}(0,0,-\pi-2\theta(t),0)$ in (\ref{TargetWigner1}) and (\ref{NewNewA1}) has a nonzero $U_2-$Euler angle. In the following section, we will diagonalize the matrix of Wigner $D$-functions and transform the $U_2-$Euler angle to a $U_3-$ Euler angle in order to compute the intertwining operators more easily.\\

\subsection{Diagonalization of the intertwining operators}\label{ProofOfProp2}
From the commutation relation of Pauli matrices, $U_2$ and $U_3$ can be related in the following way:
\begin{equation}
	U_2 = \Ad(e^{-\frac{3\pi}{2} U_1})U_3.\label{Comm1}
\end{equation}
We can use this relation to diagonalize the action of 
\[
e^{(\pi+2\theta(t))U_2} = \Ad(e^{-\frac{3\pi}{2} U_1})e^{(\pi+2\theta(t))U_3}
\]
appearing in the Iwasawa decomposition of $\exp(t X_{-\alpha_1})$ of the Wigner $D$-functions. %If we substitute the exponential of $U_2$ by the exponential of the right hand side of (\ref{Comm1}), we have
By the multiplicativity of the Wigner $D$-function, we have
\begin{align}
	&W^{(j,n)}_{m_1,m_2}(kw_{\alpha_1}e^{2\theta(t)U_2})=W^{(j,n)}_{m_1,m_2}(ke^{(\pi+2\theta(t))U_2})\nonumber\\
	=&\sum_{m_3\in \mathtt{M}(j,n;\delta_2,\delta_1)}W^{(j,n)}_{m_1,m_3}(k)W^{(j,n)}_{m_3,m_2}(e^{-\frac{3\pi}{2} U_1}e^{(\pi+2\theta(t))U_3}\exp^{\frac{3\pi}{2} U_1})\nonumber\\
	=&\sum_{m_3\in \mathtt{M}(j,n;\delta_2,\delta_1)}W^{(j,n)}_{m_1,m_3}(k)\sum_{m_4,m_5}W^{(j,n)}_{m_3,m_4}(e^{-\frac{3\pi}{2} U_1})W^{(j,n)}_{m_5,m_2}(e^{\frac{3\pi}{2} U_1})W^{(j,n)}_{m_4,m_5}(e^{(\pi+2\theta(t))U_3})\nonumber\\
	=&\sum_{m_3\in \mathtt{M}(j,n;\delta_2,\delta_1)}W^{(j,n)}_{m_1,m_3}(k)\sum_{m_4,m_5}W^{(j,n)}_{m_3,m_4}(e^{-\frac{3\pi}{2} U_1})W^{(j,n)}_{m_5,m_2}(e^{\frac{3\pi}{2} U_1})\times\nonumber\\
	&W^{(j,n)}_{m_4,m_5}(0,-\pi-2\theta(t),0,0).\label{DiagonalizedEuler}
\end{align}
We define the function $S'^{(j,n)}_{m_3,m_2}(z)$ by
\begin{align}
	S'^{(j,n)}_{m_3,m_2}(z)=&\sum_{m_4,m_5}W^{(j,n)}_{m_3,m_4}(e^{-\frac{3\pi}{2} U_1})W^{(j,n)}_{m_5,m_2}(e^{\frac{3\pi}{2} U_1})\nonumber\\
	&\int_{-\infty}^\infty (1+t^2)^{-z}W^{(j,n)}_{m_4,m_5}(0,-\pi-2\theta(t),0,0).\label{SDef}
\end{align}
\subsubsection*{Change-of-basis matrix}
Similarly to the relation (\ref{Comm1}) between $U_2$ and $U_3$, we can relate $U_1$ and $U_2$ by conjugating a multiple of $U_3$:
\begin{align}
	U_1 = \Ad(e^{\frac{\pi}{2}U_3})U_2.\label{Comm2}
\end{align}
Recalling the notation from from (\ref{WignerDefinition}) and (\ref{DPart}) that $c^j_m = \sqrt{(j+m)!(j-m)!}$, the change-of-basis matrices $W^{(j,n)}_{m_3,m_4}(e^{-\frac{3\pi}{2} U_1})$ and $W^{(j,n)}_{m_5,m_2}(e^{\frac{3\pi}{2} U_1})$ can be expressed in terms of the value of a Wigner $D$-function:
\begin{align}
	W^{(j,n)}_{m_3,m_4}&(e^{-\frac{3\pi}{2}U_1})=W^{(j,n)}_{m_3,m_4}(e^{\frac{\pi}{2}U_3}e^{-\frac{3\pi}{2}U_2}e^{-\frac{\pi}{2}U_3})=W^{(j,n)}_{m_3,m_4}\left(-\frac{\pi}{2},\frac{3\pi}{2},\frac{\pi}{2}\right)\nonumber\\
    &=c^j_{m_3}c^j_{m_4}\begin{cases}
	\frac{\ii^{m_3-m_4}(-1)^{2j}2^{-j}}{(j-m_3)!(m_3-m_4)!(j+m_4)!}{}_2F_1(\begin{smallmatrix}-j+m_3,-j-m_4\\1+m_3-m_4\end{smallmatrix};-1)&m_3>m_4\\
    \frac{\ii^{m_4-m_3}(-1)^{2j}2^{-j}}{(j+m_3)!(m_4-m_3)!(j-m_4)!}{}_2F_1(\begin{smallmatrix}-j-m_3,-j+m_4\\1-m_3+m_4\end{smallmatrix},-1)&m_3\leq m_4
\end{cases}.\label{ChangeOfBasis}
\end{align}
We define \begin{align}M^{j}_{m_3,m_4} = W^{(j,n)}_{m_3,m_4}(e^{-\frac{3\pi}{2}U_1})\label{MMDef}\\N^{j}_{m_5,m_2} = W^{(j,n)}_{m_5,m_2}(e^{\frac{3\pi}{2}U_1}).\label{NNDef}\end{align} As we have seen in (\ref{WignerJacobi}), the entries of $M^{j}_{m_3,m_4}$ and $N^{j}_{m_3,m_4}$ are related to values $P^{m_1-m_4,m_1+m_4}_{j-m_1}(0)$ of Jacobi polynomials $P^{\alpha,\beta}_n(x)$:
\begin{align}
	M^{j}_{m_3,m_4} = \frac{c^j_{m_4}}{c^j_{m_3}}(-1)^{2j}\ii^{-m_3+m_4}2^{-m_4}P^{-m_3+m_4,m_3+m_4}_{j-m_4}(0)\\
	N^{j}_{m_5,m_2} =\frac{c^j_{m_2}}{c^j_{m_5}}(-1)^{2j}\ii^{m_5-m_2}2^{-m_2}P^{-m_5+m_2,m_5+m_2}_{j-m_2}(0).
\end{align}
After introducing these notations, the function $S'^{(j,n)}_{m_3,m_2}(z)$ defined in (\ref{SDef}) becomes
\begin{equation}
	S'^{(j,n)}_{m_3,m_2}(z)=\sum_{m_4,m_5}M^{(j,n)}_{m_3,m_4}N^{(j,n)}_{m_5,m_2}\int_{-\infty}^\infty (1+t^2)^{-z}W^{(j,n)}_{m_4,m_5}(0,-\pi-2\theta(t),0,0).\label{SWrap}
\end{equation}
Among the two equivalent definitions (\ref{Jacobi1}) from \cite{AbStegun} and (\ref{Jacobi2}) in \cite{WolframAlpha}, for the sake of simplicity in expressions we choose the definition (\ref{Jacobi2}) using hypergeometric functions
\begin{align*}
P^{\alpha,\beta}_n(x)&=\binom{n+\alpha}{n}\left(\frac{x+1}{2}\right)^n {}_2F_1\left(\begin{smallmatrix}-n,-n-\beta\\\alpha+1\end{smallmatrix};\frac{x-1}{x+1}\right).
\end{align*}
These Jacobi polynomials have a generating function \cite{GeneratingFunctions}
\begin{align}
\sum_{n=0}^{\infty}P^{(\alpha-n,\beta-n)}_n(x) t^n &= \left(1+\frac{1}{2}(x+1)t\right)^{\alpha} \left(1+\frac{1}{2}(x-1)t\right)^{\beta}\label{GeneratingJacobi}
\end{align}
for $|x|<1$. The series on the left hand side converges absolutely for $|t|<1$. For any meromorphic function $f(t)$, we denote $[f(t)]_0$ as the zeroth Laurent series coefficient of $f(t)$. The change-of-basis matrices $M^{j}_{m_3,m_4}, N^{j}_{m_5,m_2}$ can thus be expressed as the constant term in  Laurent series,
\begin{align}
M^j_{m_3,m_4}&=\left[(-1)^{2j}\frac{c^j_{m_4}}{c^j_{m_3}}\ii^{-m_3+m_4}2^{-m_4}t^{m_4-j}\left(1+\frac{t}{2}\right)^{j-m_3}\left(1-\frac{t}{2}\right)^{j+m_3}\right]_0\label{MLabel}\\
N^j_{m_5,m_2}&=\left[(-1)^{2j}\frac{c^j_{m_2}}{c^j_{m_5}}\ii^{m_5-m_2}2^{-m_2}t^{m_2-j}\left(1+\frac{t}{2}\right)^{j-m_5}\left(1-\frac{t}{2}\right)^{j+m_5}\right]_0\label{NLabel}.
\end{align}
\subsubsection*{The Singular Integrals}
In (\ref{ComputeIntertwining}), we have described the procedure of calculating the matrix entries \begin{equation}[A(w,\lambda)]^{(j,n)}_{m_3,m_2} = \langle W^{(j,n)}_{m_1,m_3}, A(w,\lambda)W^{(j,n)}_{m_1,m_2}\rangle\label{AMatDef}\end{equation} of the intertwining operator $A(w,\lambda)$. The right hand side of (\ref{AMatDef}) is independent of $m_1$. The calculation of simple intertwining operators reduces to the calculation of integrals (\ref{NewNewA1}) and (\ref{NewA2}).  Combining the diagonalized operator (\ref{DiagonalizedEuler}) with (\ref{NewA1}) and (\ref{NewA2}),  the problem of the calculation of the intertwining operators $A(w,\lambda)$ reduces to the following two integrals:
\begin{align}\label{SingularInt1}
	&\int_{-\infty}^{\infty}(1+t^2)^{-\frac{\langle \check{\alpha}_1,\lambda+\rho\rangle}{2}}W^{(j,n)}_{m_4,m_5}\left(0,-\pi-2\theta(t),0,0\right)\dd t\\
	\label{SingularInt2}&\int_{-\infty}^{\infty}(1+t^2)^{-\frac{\langle \check{\alpha}_2,\lambda+\rho\rangle}{2}}W^{(j,n)}_{m_3,m_2}\left(-\frac{\pi}{2}-\theta(t),\frac{\pi}{2}+\theta(t),0,0\right)\dd t
\end{align}
Recall from the definition of Wigner $D$-functions and their values on Euler angles that
\[
	W^{(j,n)}_{m_1,m_2}(\zeta,\psi,0,0)=e^{\ii n\zeta+\ii m_1\psi}\delta_{m_1,m_2}.
\]
The Wigner $D$-function part of the integrands of the above two integrals (\ref{SingularInt1}) and (\ref{SingularInt2}) are
\begin{align*}
	W^{(j,n)}_{m_4,m_5}&\left(0,-\pi-2\theta(t),0,0\right) =\ii^{-2m_4}(1+\ii t)^{m_4}(1-\ii t)^{-m_4}\delta_{m_4,m_5}\\
	W^{(j,n)}_{m_3,m_2}&\left(-\frac{\pi}{2}-\theta(t),\frac{\pi}{2}+\theta(t),0,0\right) = \ii^{-n+m_2}(1+\ii t)^{\frac{n-m_2}{2}}(1-\ii t)^{-\frac{n-m_2}{2}}\delta_{m_3,m_2}
\end{align*}
The integral $\int_{-\infty}^{\infty}(1+\ii t)^{s_1}(1-\ii t)^{s_2}\dd t$ is convergent for $\mathrm{Re}(s_1+s_2)<-1$, and can be meromorphically continued to the whole complex plane as
\begin{align}\label{BetaFun}
	\int_{-\infty}^{\infty}(1+\ii t)^{s_1}(1-\ii t)^{s_2}\dd t = \pi 2^{s_1+s_2+2}\frac{\Gamma(-s_1-s_2-1)}{\Gamma(-s_1)\Gamma(-s_2)}.
\end{align}
We can use (\ref{BetaFun}) to express the two integrals (\ref{SingularInt1}) and (\ref{SingularInt2}) in terms of $\Gamma$-functions. Thus the integral (\ref{SingularInt1}) becomes
\begin{align*}
	&\int_{-\infty}^{\infty}(1+t^2)^{-\frac{\langle \check{\alpha}_1,\lambda+\rho\rangle}{2}}W^{(j,n)}_{m_4,m_5}(0,-\pi-2\theta(t),0,0)\dd t \\
	=& \ii^{-2m_4}\delta_{m_4,m_5}\int_{-\infty}^{\infty}(1-\ii t)^{-\frac{\langle \check{\alpha}_1,\lambda+\rho\rangle}{2}-m_4}(1+\ii t)^{-\frac{\langle \check{\alpha}_1,\lambda+\rho\rangle}{2}+m_4}\dd t\\
	=& \ii^{-2m_4}\pi 2^{-\langle \check{\alpha}_1,\lambda+\rho\rangle+2}\frac{\Gamma(\langle \check{\alpha}_1,\lambda+\rho\rangle-1)}{\Gamma\left(\frac{\langle \check{\alpha}_1,\lambda+\rho\rangle}{2}-m_4\right)\Gamma\left(\frac{\langle \check{\alpha}_1,\lambda+\rho\rangle}{2}+m_4\right)}\delta_{m_4,m_5}
\end{align*}
and the integral (\ref{SingularInt2}) is
\begin{align*}
	&\int_{-\infty}^{\infty}(1+t^2)^{-\frac{\langle \check{\alpha}_2,\lambda+\rho\rangle}{2}}W^{(j,n)}_{m_3,m_2}\left(-\frac{\pi}{2}-\theta(t),\frac{\pi}{2}+\theta(t),0,0\right)\dd t\\
	=& \ii^{-n+m_2}\delta_{m_3,m_2}\int_{-\infty}^{\infty}(1-\ii t)^{-\frac{\langle \check{\alpha}_2,\lambda+\rho\rangle}{2}-\frac{n-m_2}{2}}(1+\ii t)^{-\frac{\langle \check{\alpha}_2,\lambda+\rho\rangle}{2}+\frac{n-m_2}{2}}\dd t\\
	=&\ii^{-n+m_2}\pi 2^{-\langle \check{\alpha}_2,\lambda+\rho\rangle+2}\frac{\Gamma(\langle \check{\alpha}_2,\lambda+\rho\rangle-1)}{\Gamma\left(\frac{\langle \check{\alpha}_2,\lambda+\rho\rangle}{2}-\frac{n-m_2}{2}\right)\Gamma\left(\frac{\langle \check{\alpha}_2,\lambda+\rho\rangle}{2}+\frac{n-m_2}{2}\right)}\delta_{m_3,m_2}.
\end{align*}
We have defined the function $Q(z,n) = \frac{\pi 2^{2-2z}\Gamma(2z-1)}{\Gamma(z+n)\Gamma(z-n)}$ in (\ref{QDef}), so the integrals (\ref{SingularInt1}) and (\ref{SingularInt2}) can be expressed as:
\begin{align}
	\int_{-\infty}^{\infty}(1+t^2)^{-\frac{\langle \check{\alpha}_1,\lambda+\rho\rangle}{2}}&W^{(j,n)}_{m_4,m_5}\left(0,-\pi-2\theta(t),0,0\right)\dd t\nonumber\\
	&=\ii^{-2m_4}Q\left(\frac{\langle\check{\alpha}_1,\lambda+\rho\rangle}{2},m_4\right)\delta_{m_4,m_5}\label{SInt}\\
	\int_{-\infty}^{\infty}(1+t^2)^{-\frac{\langle \check{\alpha}_2,\lambda+\rho\rangle}{2}}&W^{(j,n)}_{m_3,m_2}\left(-\frac{\pi}{2}-\theta(t),\frac{\pi}{2}+\theta(t),0,0\right)\dd t\nonumber\\
	&=\ii^{-n+m_2}Q\left(\frac{\langle\check{\alpha}_2,\lambda+\rho\rangle}{2},\frac{m_2-n}{2}\right)\delta_{m_3,m_2}.\label{TMat}
\end{align}
If we denote by $T^n_{m_2}(z) = \ii^{-n+m_2}Q\left(z,\frac{m_2-n}{2}\right)$, the result of the integral in (\ref{TMat}) is equal to $T^n_{m_2}(\frac{\langle\check{\alpha}_2,\lambda+\rho\rangle}{2}) = \ii^{-n+m_2}Q\left(\frac{\langle\check{\alpha}_2,\lambda+\rho\rangle}{2},\frac{m_2-n}{2}\right)$. By collecting the results from (\ref{ChangeOfBasis}) and (\ref{SInt}), it turns out that the $S'^{(j,n)}_{m_3,m_2}(z)$ defined in (\ref{SDef}) and (\ref{SWrap}) is exactly the same as the $S^{(j,n)}_{m_3,m_2}(z)$ in Proposition 6.1:
\begin{align*}
	S^{(j,n)}_{m_3,m_2}\left(\frac{\langle\check{\alpha}_1,\lambda+\rho\rangle}{2}\right)&=S'^{(j,n)}_{m_3,m_2}\left(\frac{\langle\check{\alpha}_1,\lambda+\rho\rangle}{2}\right)\\
	&=\sum_{m_4}M^{(j,n)}_{m_3,m_4}N^{(j,n)}_{m_4,m_2}\ii^{-2m_4}Q\left(\frac{\langle\check{\alpha}_1,\lambda+\rho\rangle}{2},m_4\right).
\end{align*}
Therefore, by (\ref{NewNewA1}), the simple intertwining operator $A(w_{\alpha_1},\lambda)$ acts on Wigner $D$-functions by
\[
	A(w_{\alpha_1},\lambda)W^{(j,n)}_{m_1,m_2} = \sum_{-j \leq m_3\leq j}S^{(j,n)}_{m_3,m_2}\left(\frac{\langle\check{\alpha}_1,\lambda+\rho\rangle}{2}\right)W^{(j,n)}_{m_1,m_3}
\]
From (\ref{NewA2}) and (\ref{TMat}), we can see that the simple intertwining operator $A(w_{\alpha_2},\lambda)$ acts by
\[
	A(w_{\alpha_2},\lambda)W^{(j,n)}_{m_1,m_2} = T^{n}_{m_2}\left(\frac{\langle\check{\alpha}_2,\lambda+\rho\rangle}{2}\right)W^{(j,n)}_{m_1,m_2}.
\]
By the Langlands' lemma (\ref{langlandslemma})
\[
	A(w_0,\lambda) = A(w_{\alpha_2},w_{\alpha_1}w_{\alpha_2}w_{\alpha_1}\lambda)A(w_{\alpha_1},w_{\alpha_2}w_{\alpha_1}\lambda)A(w_{\alpha_2},w_{\alpha_1}\lambda)A(w_{\alpha_1},\lambda),
\] we can replace $\lambda$ in the formula above by $w_{\alpha_1}\lambda$, $w_{\alpha_2}w_{\alpha_1}\lambda$ and $w_{\alpha_1}w_{\alpha_2}w_{\alpha_1}\lambda$ to calculate the other stages in the composition of simple intertwining operators. Recall that if we express $\lambda$ by the pair $(\lambda_1,\lambda_2)$ as in Section 6.1.1, where $\lambda_i=\lambda(H_i)$, the composition of actions by the Weyl group action will send the pair $(\lambda_1,\lambda_2)$ to:
\[
	(\lambda_1,\lambda_2)\xrightarrow{w_{\alpha_1}}(\lambda_2,\lambda_1)\xrightarrow{w_{\alpha_2}}(\lambda_2,-\lambda_1)\xrightarrow{w_{\alpha_1}}(-\lambda_1,\lambda_2)\xrightarrow{w_{\alpha_2}}(-\lambda_1,-\lambda_2).
\]
The pairing of $\lambda+\rho$ with the coroots $\check{\alpha}_1$ and $\check{\alpha}_2$ are given by
\begin{align*}
	\langle\check{\alpha}_1,\lambda+\rho\rangle &= \lambda_1-\lambda_2+1\\
	\langle\check{\alpha}_2,\lambda+\rho\rangle &= \lambda_2+1
\end{align*}
Therefore, we can conclude that the matrix coefficients for the simple intertwining operators are
\begin{align*}
	A(w_{\alpha_1},\lambda)W^{(j,n)}_{m_1,m_2} &= \sum_{m_3\in \mathtt{M}(j,n;\delta_2,\delta_1)}S^{(j,n)}_{m_3,m_2}\left(\frac{\langle\check{\alpha}_1,\lambda+\rho\rangle}{2}\right)W^{(j,n)}_{m_1,m_3}\\
 A(w_{\alpha_2},w_{\alpha_1}\lambda)W^{(j,n)}_{m_1,m_2} &= T^{n}_{m_2}\left(\frac{\langle\check{\alpha}_2,w_{\alpha_1}\lambda+\rho\rangle}{2}\right)W^{(j,n)}_{m_1,m_2}\\
   A(w_{\alpha_1},w_{\alpha_2}w_{\alpha_1}\lambda)W^{(j,n)}_{m_1,m_2} &= \sum_{m_3\in \mathtt{M}(j,n;\delta_2,\delta_1)}S^{(j,n)}_{m_3,m_2}\left(\frac{\langle\check{\alpha}_1,w_{\alpha_2}w_{\alpha_1}\lambda+\rho\rangle}{2}\right)W^{(j,n)}_{m_1,m_3}\\
 A(w_{\alpha_2},w_{\alpha_1}w_{\alpha_2}w_{\alpha_1}\lambda)W^{(j,n)}_{m_1,m_2} &= T^{n}_{m_2}\left(\frac{\langle\check{\alpha}_2,w_{\alpha_1}w_{\alpha_2}w_{\alpha_1}\lambda+\rho\rangle}{2}\right)W^{(j,n)}_{m_1,m_2}
\end{align*}
Collecting the above result with the change-of-basis matrix defined in (\ref{ChangeOfBasis}), we have finished the proof of the Proposition \ref{GindikinKarpelevich}.
\subsection{Expression of $S^{j,n}_{m_1,m_4}(z)$ as Hypergeometric Functions}
In this section, we assume
\[
	(\delta_1,\delta_2) = (0,0) \text{ or }(1,1),
\]
in which case $j,n$ are integers. We also introduce the rising factorial or the \emph{Pochhammer symbol}
\begin{equation}
	(a)^{(n)}=\frac{\Gamma (a+n)}{\Gamma (a)}
\end{equation}
\begin{definition}\label{DefMellinTransform}Part I of \cite{TableMellinTransform}.
\begin{enumerate}
\item  The \emph{Mellin transform} of a function $f(x)$ is formally defined as
\begin{equation}
	\mathcal{M}(f(x))(z) = \int_0^\infty f(x)x^{z-1}\dd x.
\end{equation}
\item Consider a function $F(z)$ of one complex variable $z = \sigma+\ii\tau$, such that
\begin{enumerate}
	\item $F(z)$ is holomorphic on the strip $S = \{z\in\CC|a<\sigma<b\}$, such that $F(z)\to 0$ uniformly in the strip $S_\epsilon = \{z\in\CC|a+\epsilon<\sigma<b-\epsilon\}$ for arbitrarily small $\epsilon >0$, and
	\item $\int_{-\infty}^{\infty}|F(\sigma+\ii\tau)|\dd\tau<\infty$ for all $\sigma\in(a,b)$.
\end{enumerate}
We define the \emph{inverse Mellin transform} of the function $F(z)$ as
\begin{equation}\label{InverseMellinTransform}
	\mathcal{M}^{-1}(F(z))(x) = \frac{1}{2\pi \ii}\int_{\gamma-\ii\infty}^{\gamma+\ii\infty}F(z)x^{-z}\dd z,
\end{equation}
for $x>0$ and some fixed $\gamma\in (a,b)$. It satisfies the property that
\begin{equation}
	\mathcal{M}(\mathcal{M}^{-1}(F))(z) = F(z).
\end{equation}
\end{enumerate}
\end{definition}
The function $F(z)$ to our interest is the ratio of $\Gamma$-functions $\frac{\Gamma(z)}{\Gamma\left(\frac{1+z-\nu}{2}\right)\Gamma\left(\frac{1+z+\nu}{2}\right)}$. We would like to calculate the integral
\begin{equation}
	\frac{1}{2\pi \ii}\int_{\gamma-\ii\infty}^{\gamma+\ii\infty}\frac{\Gamma(z)}{\Gamma\left(\frac{1+z-\nu}{2}\right)\Gamma\left(\frac{1+z+\nu}{2}\right)}x^{-z}\dd z\label{MyMellin}
\end{equation}
along some well chosen contour on which the integral converges. Though this integral doesn't satisfy the assumption on $F(z)$ above, we still can take the inverse Mellin transform of this function, resulting in a $\mathcal{M}^{-1}(F(z))(x)$ with discontinuities. We refer to \cite[Page 49, Section 1.19]{Bateman} for a general \emph{Mellin-Barnes integral}
\begin{equation}
	\frac{1}{2\pi \ii}\int_{\gamma-\ii\infty}^{\gamma+\ii\infty}\frac{\prod_{i=1}^m\Gamma(a_i+A_i z)\prod_{j=1}^n\Gamma(b_j-B_j z)}{\prod_{s=1}^p\Gamma(c_s+C_s z)\prod_{t=1}^q\Gamma(d_t-D_t z)}x^{-z}\dd z,
\end{equation}
where $\gamma\in\RR$ and $A_i,B_j,C_s, D_t >0$. If $x\in\RR$, the integrand is asymptotically $e^{-\frac{1}{2}\alpha\pi |y|}|y|^{\beta\gamma+\lambda}|\frac{x}{\rho}|^{-\gamma}$, where the numbers $\alpha, \beta, \lambda,\rho$ are defined by
\begin{align*}
	\alpha &= \sum_{i=1}^m A_i + \sum_{j=1}^n B_j - \sum_{s=1}^p C_s - \sum_{t=1}^q D_t\\
	\beta &= \sum_{i=1}^m A_i - \sum_{j=1}^n B_j - \sum_{s=1}^p C_s + \sum_{t=1}^q D_t\\
	\lambda &= \mathrm{Re}(\sum_{i=1}^m a_i + \sum_{j=1}^n b_j - \sum_{s=1}^p c_s - \sum_{t=1}^q d_j) - \frac{m+n-p-q}{2}\\
	\rho &= \prod_{i=1}^m A_i^{A_i} \prod_{j=1}^n {B_j}^{-B_j} \prod_{s=1}^p {C_s}^{-C_s}  \prod_{t=1}^q D_t^{D_t},
\end{align*}
and $y  = \mathrm{Im}(z)$.
Our integral (\ref{MyMellin}) satisfies $\alpha=\beta=0$ and $\lambda = -\frac{1}{2}$, $\rho = 2$, falling into the fourth type in the description on the convergence of the Mellin-Barnes integral in \cite{Bateman}, which states that the integral (\ref{MyMellin}) conditionally converges to an analytic function in $x$ on the intervals $|x|<2$ and $|x|>2$, with points of discontinuity at $x=\pm 2$. In \cite[Section II, (5.21)]{TableMellinTransform}, for $\mathrm{Re}(z)>0$, we have an explicit formula for (\ref{MyMellin}):
\begin{equation}
	\mathcal{M}^{-1}\left(\frac{\Gamma(z)}{\Gamma\left(\frac{1+z-\nu}{2}\right)\Gamma\left(\frac{1+z+\nu}{2}\right)}\right)(x) =\left\{\begin{smallmatrix} 2\pi^{-1}(4-x^2)^{-1/2}\cos\left(\nu\arccos(x/2)\right)&|x|<2\\0&|x|>2\end{smallmatrix}\right.
\end{equation}
In the rest of this chapter, we will only interchange finite sums with this conditionally convergent integral, and thus there are no analytic issues justifying the calculation in the remaining portion of this thesis. 
In (\ref{QDef}) we have defined the function $Q(z,n) =  \frac{\pi 2^{2-2z}\Gamma(2z-1)}{\Gamma(z+n)\Gamma(z-n)}$. We would like to compute the inverse Mellin transform $\mathcal{M}^{-1}(Q(z,m_3))(x)$ of $Q(z,m_3)$:
\begin{align}
	\mathcal{M}^{-1}\left(\frac{\pi 2^{2-2z}\Gamma(2z-1)}{\Gamma(z+m_3)\Gamma(z-m_3)}\right)&(x) = \frac{4\pi}{2\pi\ii}\int_{\gamma-\ii\infty}^{\gamma+\ii\infty} \frac{\Gamma(2z-1)}{\Gamma(z+m_3)\Gamma(z-m_3)}(4x)^{-z}\dd z\nonumber\\
	&= \frac{2\pi}{2\pi\ii}\int_{2\gamma-1-\ii\infty}^{2\gamma-1+\ii\infty} \frac{\Gamma(z)}{\Gamma(\frac{z+1+2m_3}{2})\Gamma(\frac{z+1-2m_3}{2})}(4x)^{-\frac{z+1}{2}}\dd z\nonumber\\
	&=\frac{\pi}{2\pi\ii\sqrt{x}}\int_{2\gamma-1-\ii\infty}^{2\gamma-1+\ii\infty} \frac{\Gamma(z)}{\Gamma(\frac{z+1+2m_3}{2})\Gamma(\frac{z+1-2m_3}{2})}(2\sqrt{x})^{-z}\dd z\nonumber\\
	&=\frac{\pi}{\sqrt{x}}\mathcal{M}^{-1}\left(\frac{\Gamma(z)}{\Gamma\left(\frac{1+z-2m_3}{2}\right)\Gamma\left(\frac{1+z+2m_3}{2}\right)}\right)(2\sqrt{x})\nonumber\\
	&=\frac{\pi}{\sqrt{x}}\left\{\begin{smallmatrix} \pi^{-1}(1-x)^{-1/2}\cos\left(2m_3\arccos(\sqrt{x})\right)&0<x<1\\0&x>1\end{smallmatrix}\right.\nonumber\\
	&=\frac{1}{\sqrt{x(1-x)}}\left\{\begin{smallmatrix}\cos\left(2m_3\arcsin(\sqrt{1-x})\right)&|x|<1\\0&|x|>1\end{smallmatrix}\right..
\end{align}
If we define the \emph{Heaviside step function} $\theta(x) = \left\{\begin{smallmatrix}0&x<0\\1&x>0\end{smallmatrix}\right.$, the result above can be written as
\begin{align}
	\mathcal{M}^{-1}(Q(z,m_3))(x) &= \frac{1-\theta (\left| x\right| -1)}{  \sqrt{(1-x) x}}\cos \left(2 m_3 \arcsin
  \left(\sqrt{1-x}\right)\right)\nonumber\\
  &=\frac{1-\theta (\left| x\right| -1)}{  \sqrt{(1-x) x}}\left(\frac{1}{2} \left(\sqrt{x}+\ii \sqrt{1-x}\right)^{-2 m_3}+\frac{1}{2} \left(\sqrt{x}+\ii
   \sqrt{1-x}\right)^{2 m_3}\right)
\end{align}
Therefore, the inverse Mellin transform $\mathcal{M}^{-1}(S^{j,n}_{m_1,m_4}(z))(x)$ of the matrix entries of intertwining operators $S^{j,n}_{m_3,m_2}(z)$ is
\begin{align}
	\mathcal{M}^{-1}(S^{j,n}_{m_1,m_4}(z))(x)&=\frac{1-\theta (\left| x\right| -1)}{2  \sqrt{(1-x) x}}\sum_{m_3=-j}^j M^j_{m_1,m_3}N^j_{m_3,m_4}\nonumber\\
	&\ii^{-2m_3}\left(\left(\sqrt{x}+\ii \sqrt{1-x}\right)^{-2 m_3}+\left(\sqrt{x}+\ii \sqrt{1-x}\right)^{2 m_3}\right).\label{SumOfCoef}
\end{align}
By (\ref{MLabel}) and (\ref{NLabel}), the change-of-basis matrix $M^j_{m_1,m_3}$ and $N^j_{m_1,m_3}$ are the zeroth Laurent coefficient of some analytic functions. By substituting the matrix entries $M^j_{m_1,m_3}$ and $N^j_{m_1,m_3}$ in (\ref{SumOfCoef}) with the corresponding Laurent series (\ref{MLabel}) and (\ref{NLabel}), the sum over $m_3$ from $-j$ to $j$ as in (\ref{SumOfCoef}) actually gives rise to a geometric series. Recalling formulas (\ref{MLabel}) and (\ref{NLabel}) and defining
\begin{align*}
M^j_{m_1,m_3}(t)&=(-1)^{2j}\frac{c^j_{m_3}}{c^j_{m_1}}\ii^{-m_1+m_3}2^{-m_3}t^{m_3-j}\left(1+\frac{t}{2}\right)^{j-m_1}\left(1-\frac{t}{2}\right)^{j+m_1}\\
N^j_{m_3,m_4}(s)&=(-1)^{2j}\frac{c^j_{m_4}}{c^j_{m_3}}\ii^{m_3-m_4}2^{-m_4}s^{m_4-j}\left(1+\frac{s}{2}\right)^{j-m_3}\left(1-\frac{s}{2}\right)^{j+m_3},
\end{align*}
where $s,t$ are formal variables, we define
\begin{align}
	\sigma_1(s,t) = &\sum_{m_3=-j}^j M^j_{m_1,m_3}(t)N^j_{m_3,m_4}(s)\ii^{-2m_3}\left(\sqrt{x}+\ii \sqrt{1-x}\right)^{-2 m_3} \nonumber\\
	=&\frac{c^j_{m_4}}{c^j_{m_1}}\ii^{-m_1-m_4} 2^{-m_4} \left(1-\frac{t}{2}\right)^{j+m_1}
   \left(1+\frac{t}{2}\right)^{j-m_1} s^{m_4-j}t^{-j}\left(1+\frac{s}{2}\right)^{j}\left(1-\frac{s}{2}\right)^{j}\times\nonumber\\
   &\sum_{m_3=-j}^j \left(\frac{t \left(1-\frac{s}{2}\right)}{2\left(1+\frac{s}{2}\right)}\right)^{m_3}\left(\sqrt{x}+\ii \sqrt{1-x}\right)^{-2 m_3}\nonumber\\
   =&\frac{c^j_{m_4}}{c^j_{m_1}}\ii^{-m_1-m_4} 2^{-m_4} \left(1-\frac{t}{2}\right)^{j+m_1}
   \left(1+\frac{t}{2}\right)^{j-m_1} s^{m_4-j}t^{-j}\left(1+\frac{s}{2}\right)^{j}\left(1-\frac{s}{2}\right)^{j}\times\nonumber\\
   &\frac{\left(\frac{(1-s/2)t}{2(1+s/2)(\ii\sqrt{1-x}+\sqrt{x})^2}\right)^{-j} - \left(\frac{(1-s/2)t}{2(1+s/2)(\ii\sqrt{1-x}+\sqrt{x})^2}\right)^{j+1}}
   {1-\frac{(1-s/2)t}{2(1+s/2)(\ii\sqrt{1-x}+\sqrt{x})^2}}
  \end{align}
and
  \begin{align}
\sigma_2(s,t) = &\sum_{m_3=-j}^j M^j_{m_1,m_3}(t)N^j_{m_3,m_4}(s)\ii^{-2m_3}\left(\sqrt{x}+\ii \sqrt{1-x}\right)^{2 m_3} \nonumber\\
	=&\frac{c^j_{m_4}}{c^j_{m_1}}\ii^{-m_1-m_4} 2^{-m_4} \left(1-\frac{t}{2}\right)^{j+m_1}
   \left(1+\frac{t}{2}\right)^{j-m_1} s^{m_4-j}t^{-j}\left(1+\frac{s}{2}\right)^{j}\left(1-\frac{s}{2}\right)^{j}\times\nonumber\\
   &\sum_{m_3=-j}^j \left(\frac{t \left(1-\frac{s}{2}\right)}{2\left(1+\frac{s}{2}\right)}\right)^{m_3}\left(\sqrt{x}+\ii \sqrt{1-x}\right)^{2 m_3}\nonumber\\
   =&\frac{c^j_{m_4}}{c^j_{m_1}}\ii^{-m_1-m_4} 2^{-m_4} \left(1-\frac{t}{2}\right)^{j+m_1}
   \left(1+\frac{t}{2}\right)^{j-m_1} s^{m_4-j}t^{-j}\left(1+\frac{s}{2}\right)^{j}\left(1-\frac{s}{2}\right)^{j}\times\nonumber\\
   &\frac{\left(\frac{(1-s/2)t(\ii\sqrt{1-x}+\sqrt{x})^2}{2(1+s/2)}\right)^{-j} - \left(\frac{(1-s/2)t(\ii\sqrt{1-x}+\sqrt{x})^2}{2(1+s/2)}\right)^{j+1}}
   {1-\frac{(1-s/2)t(\ii\sqrt{1-x}+\sqrt{x})^2}{2(1+s/2)}}.
 \end{align}
If we denote by $[\sigma_i(s,t)]_{0,0}$ for the constant term in the Laurent series expansion of $\sigma_i(s,t)$, then (\ref{SumOfCoef}) can be expressed as
\begin{align}
	\mathcal{M}^{-1}(S^{j,n}_{m_1,m_4}(z))(x)&=\frac{1-\theta (\left| x\right| -1)}{2  \sqrt{(1-x) x}}([\sigma_1(s,t)+\sigma_2(s,t)]_{0,0}).
\end{align}
To make the calculation simpler, and since we only care about the \emph{constant} term of the Laurent series $\sigma_i(s,t)$, we can feel free to make change of variables so that $\sigma_1(s,t)$ and $\sigma_2(s,t)$ have the same denominators. Two such choices are
\begin{align}
	&\sigma_1(-s,2/t) \nonumber\\
	=&\frac{c^j_{m_4}}{c^j_{m_1}}\ii^{-m_1-m_4} 2^{-j-m_4} (-1)^{m_4-j}\left(1-\frac{1}{t}\right)^{j+m_1}
   \left(1+\frac{1}{t}\right)^{j-m_1} s^{m_4-j}t^{j}\left(1+\frac{s}{2}\right)^{j}\left(1-\frac{s}{2}\right)^{j}\times\nonumber\\
    &\frac{\left(\frac{(1+s/2)}{t(1-s/2)(\ii\sqrt{1-x}+\sqrt{x})^2}\right)^{-j} - \left(\frac{(1+s/2)}{t(1-s/2)(\ii\sqrt{1-x}+\sqrt{x})^2}\right)^{j+1}}
   {1-\frac{(1+s/2)}{t(1-s/2)(\ii\sqrt{1-x}+\sqrt{x})^2}}\nonumber\\
   =&\frac{c^j_{m_4}}{c^j_{m_1}}\ii^{-m_1-m_4} 2^{-j-m_4} (-1)^{m_4-j}\left(t-1\right)^{j+m_1}
   \left(t+1\right)^{j-m_1} s^{m_4-j}t^{-j}\left(1+\frac{s}{2}\right)^{j}\left(1-\frac{s}{2}\right)^{j}\times\nonumber\\
    &\frac{\left(\frac{(1-s/2)t(\ii\sqrt{1-x}+\sqrt{x})^2}{(1+s/2)}\right)^{-j} - \left(\frac{(1-s/2)t(\ii\sqrt{1-x}+\sqrt{x})^2}{(1+s/2)}\right)^{j+1}}
   {1-\frac{(1-s/2)t(\ii\sqrt{1-x}+\sqrt{x})^2}{(1+s/2)}}
 \end{align}
 and
 \begin{align}
	&\sigma_2(s,2t) \nonumber\\
	=&\frac{c^j_{m_4}}{c^j_{m_1}}\ii^{-m_1-m_4} 2^{-j-m_4} (-1)^{j+m_1}\left(t-1\right)^{j+m_1}
   \left(1+t\right)^{j-m_1} s^{m_4-j}t^{-j}\left(1+\frac{s}{2}\right)^{j}\left(1-\frac{s}{2}\right)^{j}\times\nonumber\\
   &\frac{\left(\frac{(1-s/2)t(\ii\sqrt{1-x}+\sqrt{x})^2}{(1+s/2)}\right)^{-j} - \left(\frac{(1-s/2)t(\ii\sqrt{1-x}+\sqrt{x})^2}{(1+s/2)}\right)^{j+1}}
   {1-\frac{(1-s/2)t(\ii\sqrt{1-x}+\sqrt{x})^2}{(1+s/2)}}
\end{align}
After taking the constant Laurent series coefficient of the function $\sigma_1(-s,2/t)+\sigma_2(s,2t)$, we will get the following expression for $\mathcal{M}^{-1}(S^{j,n}_{m_1,m_4}(z))$: 
\begin{align}
\mathcal{M}^{-1}&(S^{j,n}_{m_1,m_4}(z))(x)=\frac{c^j_{m_4}}{c^j_{m_1}}\ii^{-m_1-m_4}2^{-j-m_4}\frac{(-1)^{j+m_1}+(-1)^{m_4-j}}{2}\nonumber\\
&\frac{1-\theta (\left| x\right| -1)}{\sqrt{(1-x)x}}\left[\frac{1}{1-(\ii \sqrt{1-x}+\sqrt{x})^2\frac{(1-\frac{s}{2})t}{1+\frac{s}{2}}}\right.\nonumber\\
    &\left. \left(\left(1+\frac{s}{2}\right)^{2 j} \left(1-t\right)^{j+m_1}
   \left(1+t\right)^{j-m_1} s^{m_4-j}t^{-2 j}\left(\ii \sqrt{1-x}+\sqrt{x}\right)^{-2 j}
  \right.  \right.\nonumber\\
    &\left.\left.- \left(1-\frac{s}{2}\right)^{2 j+1}
    \left(1-t\right)^{j+m_1}
   \left(1+t\right)^{j-m_1}\left(1+\frac{s}{2}\right)^{-1} s^{m_4-j} t\right.\right.\nonumber\\
   &\left.\left.\left(\ii \sqrt{1-x}+\sqrt{x}\right)^{2 j+2} \right)\right]_{0,0}\label{GeometricSum}
\end{align}
Notice that in the Laurent series expansion of the function in $s,t$ in (\ref{GeometricSum}), the exponent of $t$ in the Laurent expansion of the second term \[\left(1-\frac{s}{2}\right)^{2 j+1}
    \left(1-\frac{t}{2}\right)^{j+m_1}
   \left(1+\frac{t}{2}\right)^{j-m_1}\left(1+\frac{s}{2}\right)^{-1} s^{m_4-j} t\] in the parentheses is always strictly greater than 0. Therefore it does not contribute to the zeroth Laurent series coefficient in the expansion of (\ref{GeometricSum}):
 \begin{align}
\mathcal{M}^{-1}&(S^{j,n}_{m_1,m_4}(z))(x)=\frac{c^j_{m_4}}{c^j_{m_1}}\ii^{-m_1-m_4}2^{-j-m_4}\frac{(-1)^{j+m_1}+(-1)^{m_4-j}}{2}\nonumber\\
&\frac{1-\theta (\left| x\right| -1)}{\sqrt{(1-x)x}}\left[\frac{1}{1-(\ii \sqrt{1-x}+\sqrt{x})^2\frac{(1-\frac{s}{2})t}{1+\frac{s}{2}}}\right.\nonumber\\
    &\left. \left(\left(1+\frac{s}{2}\right)^{2 j} \left(1-t\right)^{j+m_1}
   \left(1+t\right)^{j-m_1} s^{m_4-j}t^{-2 j}\left(\ii \sqrt{1-x}+\sqrt{x}\right)^{-2 j}
  \right.  \right]_{0,0}\nonumber\\
  &=\frac{c^j_{m_4}}{c^j_{m_1}}\ii^{-m_1-m_4}2^{-j-m_4}\frac{(-1)^{j+m_1}+(-1)^{m_4-j}}{2}\frac{1-\theta (\left| x\right| -1)}{\sqrt{(1-x)x}}\nonumber\\
&\left[\sum_{q\geq 0,q\in\ZZ} \left(\left(1+\frac{s}{2}\right)^{2 j-q}\left(1-\frac{s}{2}\right)^{q} \left(1-t\right)^{j+m_1}
   \left(1+t\right)^{j-m_1}\times\right.\right.\nonumber\\
   &\left.\left. s^{m_4-j}t^{q-2 j}\left(\ii \sqrt{1-x}+\sqrt{x}\right)^{2(q- j)}
  \right)  \right]_{0,0},
\end{align}
which expands to
\begin{align}
  &\mathcal{M}^{-1}(S^{j,n}_{m_1,m_4}(z))(x)=\frac{c^j_{m_4}}{c^j_{m_1}}\ii^{-m_1-m_4}\frac{(-1)^{2j+m_1-m_4}+1}{2}\frac{1-\theta (\left| x\right| -1)}{\sqrt{(1-x)x}}\times\nonumber\\
  &\sum_{q,\kappa_1,\kappa_2\in\ZZ}\begin{spmatrix}2j-q\\j-m_4-\kappa_1\end{spmatrix}\begin{spmatrix}q\\\kappa_1\end{spmatrix}\begin{spmatrix}j+m_1\\\kappa_2\end{spmatrix}\begin{spmatrix}j-m_1\\2j-q-\kappa_2\end{spmatrix}(-1)^{\kappa_1+\kappa_2}\left(\ii \sqrt{1-x}+\sqrt{x}\right)^{2(q- j)}.\label{HypSum1}
\end{align}
The summations over $\kappa_1$ and $\kappa_2$ can be expressed as hypergeometric values
\begin{align}
\sum_{\kappa_1\in\ZZ}\begin{spmatrix}2j-q\\\kappa_1\end{spmatrix}\begin{spmatrix}q\\j+m_4-\kappa_1\end{spmatrix}(-1)^{\kappa_1}&=\begin{spmatrix}2j-q\\j-m_4\end{spmatrix}{}_2F_1(\begin{smallmatrix}-q, -j+m_4\\1+j-q+m_4\end{smallmatrix};-1)\\
\sum_{\kappa_2\in\ZZ}\begin{spmatrix}j+m_1\\\kappa_2\end{spmatrix}\begin{spmatrix}j-m_1\\2j-q-\kappa_2\end{spmatrix}(-1)^{\kappa_2} &=\begin{spmatrix}j-m_1\\2j-q\end{spmatrix}{}_2F_1(\begin{smallmatrix}-2j+q, -j-m_1\\1-j+q-m_1\end{smallmatrix};-1).
\end{align}
 The sum in (\ref{HypSum1}) becomes
 \begin{align}
 	\sum_{0\leq q\leq 2j, q\in\ZZ}\begin{spmatrix}2j-q\\j-m_4\end{spmatrix}\begin{spmatrix}j-m_1\\2j-q\end{spmatrix}{}_2F_1(\begin{smallmatrix}-q, -j+m_4\\1+j-q+m_4\end{smallmatrix};-1){}_2F_1(\begin{smallmatrix}-2j+q, -j-m_1\\1-j+q-m_1\end{smallmatrix};-1)\left(\ii \sqrt{1-x}+\sqrt{x}\right)^{2(q- j)}.
 \end{align}
Again, by writing the values of hypergeometric function $_{}2F_1$ at -1 in terms of Jacobi polynomials, we have
\begin{align}
	{}_2F_1(\begin{smallmatrix}-q, -j+m_4\\1+j-q+m_4\end{smallmatrix};-1)&=2^{j-m_4}\frac{(j-m_4)!(j+m_4-q)!}{(2j-q)!}P^{j+m_4-q,-j+m_4+q}_{j-m_4}(0)\\
	{}_2F_1(\begin{smallmatrix}-2j+q, -j-m_1\\1-j+q-m_1\end{smallmatrix};-1)&=2^{j+m_1}\frac{(j+m_1)!(-j-m_1+q)!}{q!}P^{j-m_1-q,-j-m_1+q}_{j+m_1}(0),
\end{align}
the sum (\ref{HypSum1}) is equal to
\begin{align}
	\sum_{0\leq q\leq 2j, q\in\ZZ}&2^{2j+m_1-m_4}\frac{(j-m_1)!(j+m_1)!}{q!(2j-q)!}P^{j+m_4-q,-j+m_4+q}_{j-m_4}(0)P^{j-m_1-q,-j-m_1+q}_{j+m_1}(0)\times\nonumber\\
	&\left(\ii \sqrt{1-x}+\sqrt{x}\right)^{2(q- j)}.
\end{align}
Again by (\ref{GeneratingJacobi}), (\ref{HypSum1}) is the constant term of the Laurent seris expansion of the following function in $s$ and $t$
\begin{align}
	&s^{m_4-j}t^{-j-m_1}\sum_{0\leq q\leq 2j, q\in\ZZ}\frac{(j-m_1)!(j+m_1)!}{q!(2j-q)!}(1+s)^{2j-q}
	(1-s)^{q}(1-t)^{2j-q}\times\nonumber\\&(1+t)^{q}\left(\ii \sqrt{1-x}+\sqrt{x}\right)^{2(q- j)}\nonumber\\
	=&s^{m_4-j}t^{-j-m_1}\frac{(j-m_1)!(j+m_1)!}{(2j)!}\times\nonumber\\
	&\left(\left(\ii \sqrt{1-x}+\sqrt{x}\right)^{-1}(1+s)(1-t)+\left(\ii \sqrt{1-x}+\sqrt{x}\right)(1-s)(1+t)\right)^{2j}.\label{Simplify01}
\end{align}
For simplicity of notations, we set $\xi = \ii \sqrt{1-x}+\sqrt{x}$. The parenthesis of (\ref{Simplify01}) can be reorganized as
\begin{align}
	&\left(\xi^{-1}(1+s)(1-t)+\xi(1-s)(1+t)\right)^{2j} \nonumber\\
	=&\left(\xi+\xi^{-1}) + (\xi^{-1}-\xi)s - (\xi^{-1}-\xi) t- st (\xi+\xi^{-1}\right)^{2j}\nonumber\\
	=&\sum_{\substack{\nu_1,\nu_2,\nu_3,\nu_4\\
	\sum \nu_i = 2j}}\frac{(2j)!}{\nu_1!\nu_3!\nu_2!\nu_4!}(\xi+\xi^{-1})^{\nu_1+\nu_4}(\xi^{-1}-\xi)^{\nu_3+\nu_2}(-1)^{\nu_2+\nu_4}s^{\nu_3+\nu_4}t^{\nu_2+\nu_4}.
\end{align}
Therefore, the constant term of (\ref{Simplify01}) can be rewritten into a sum involving multinomial coefficients, noting that $2j+m_1-m_4+\nu_1-\nu_4=2j$:
\begin{align}
	&(-1)^{j+m_1}\sum_{\nu_1,\nu_4}\frac{(j-m_1)!(j+m_1)!}{\nu_1!(j-m_4-\nu_4)!\nu_4!(j+m_1-\nu_4)!}(\xi+\xi^{-1})^{\nu_1+\nu_4}(\xi^{-1}-\xi)^{2j+m_1-m_4-2\nu_4}\nonumber\\
	=&(-1)^{j+m_1}\sum_{\nu_1,\nu_4}\frac{(j-m_1)!(j+m_1)!}{\nu_1!(j-m_4-\nu_4)!\nu_4!(j+m_1-\nu_4)!}(\xi+\xi^{-1})^{\nu_1+\nu_4}(\xi^{-1}-\xi)^{2j-\nu_1-\nu_4}.
\end{align}
Reorganizing,
\begin{align}
	&\left(\frac{c^j_{m_1}}{c^j_{m_4}}\right)^2(-1)^{j+m_1}\sum_{\nu_1,\nu_4}\frac{(j-m_4)!(j+m_4)!}{\nu_1!(j-m_4-\nu_4)!\nu_4!(j+m_4-\nu_1)!}(\xi+\xi^{-1})^{\nu_1+\nu_4}(\xi^{-1}-\xi)^{2j-\nu_1-\nu_4}\nonumber\\
	=&\left(\frac{c^j_{m_1}}{c^j_{m_4}}\right)^2(-1)^{j+m_1}\sum_{\nu_1,\nu_4}\frac{(j-m_4)!(j+m_4)!}{\nu_1!(j-m_4-\nu_4)!\nu_4!(j+m_4-\nu_1)!}(\sqrt{x})^{\nu_1+\nu_4}(-\ii\sqrt{1-x})^{2j-\nu_1-\nu_4}
\end{align}
We can observe that the function in the summation above can be simplified to the zeroth Laurent coefficient of a function in one variable $u$:
\begin{align}
	&\sum_{\nu_1,\nu_4}\frac{(j-m_4)!(j+m_4)!}{\nu_1!(j-m_4-\nu_4)!\nu_4!(j+m_4-\nu_1)!}(\sqrt{x})^{\nu_1+\nu_4}(-\ii\sqrt{1-x})^{2j-\nu_1-\nu_4}\nonumber\\
	=&\left[u^{-j-m_1}\sum_{\nu_1,\nu_4} \begin{spmatrix}j+m_4\\\nu_1\end{spmatrix}\begin{spmatrix}j-m_4\\\nu_4\end{spmatrix} (-\ii\sqrt{1-x} u)^{j+m_4-\nu_1}(-\ii\sqrt{1-x})^{j-m_4-\nu_2}u^{\nu_4}(\sqrt{x})^{\nu_1+\nu_4}\right]_0\nonumber\\
	=&\left[u^{-j-m_1}(\sqrt{x}-\ii u \sqrt{1-x})^{j+m_4}(u\sqrt{x}-\ii \sqrt{1-x})^{j-m_4}\right]_0
\end{align}
Plugging the formula above back into (\ref{HypSum1}),
\begin{align}
  &\mathcal{M}^{-1}(S^{j,n}_{m_1,m_4}(z))(x)=\frac{c^j_{m_1}}{c^j_{m_4}}\ii^{-m_1-m_4}\frac{(-1)^{2j+m_1-m_4}+1}{2}\frac{1-\theta (\left| x\right| -1)}{\sqrt{(1-x)x}}\times\nonumber\\
  &(-\ii)^{2j}
  (-1)^{j+m_1}\left[u^{-j-m_1}(\sqrt{x}+ u \sqrt{1-x})^{j+m_4}(\ii u\sqrt{x}+\sqrt{1-x})^{j-m_4}\right]_0,
\end{align}
and recall that $j,m_1,m_4\in\ZZ$ and $m_1-m_4$ is even, we finally have
\begin{align}
\mathcal{M}^{-1}(S^{j,n}_{m_1,m_4}(z))(x)& = \frac{c^j_{m_1}}{c^j_{m_4}}\frac{((-1)^{2j+m_1-m_4}+1)(-1)^{\frac{m_1-m_4}{2}}(1-\theta (\left| x\right| -1))}{2\sqrt{(1-x)x}}\nonumber\\
    &\left[\frac{(u\sqrt{1-x}+\ii\sqrt{x})^{j+m_4}(\sqrt{1-x}+\ii u\sqrt{x})^{j-m_4}}{u^{j+m_1}}\right]_{0}.
\end{align}
From the generating function of Jacobi polynomials which we introduced in (\ref{GeneratingJacobi}), this inverse Mellin transform is in fact
\begin{align}
&\mathcal{M}^{-1}(S^{j,n}_{m_1,m_4}(z))(x)\nonumber\\
=&\frac{c^j_{m_1}}{c^j_{m_4}}\frac{((-1)^{2j+m_1}\ii^{2j}+(-1)^{j+m_4})(1-\theta (\left| x\right| -1))}{2(1-x)^{\frac{m_1+m_4+1}{2}}x^{-\frac{2j+m_1+m_4-1}{2}}}P^{(-m_1-m_4,-2j-1)}_{j+m_1}\left(\frac{2}{x}-1\right).\label{JacobiForm1}
\end{align}
We refer to formulas 05.06.26.0002.01 and 07.23.17.0055.01 of \cite{WolframAlpha},
\begin{align}
	P_n^{(a,b)}(z)&=\frac{\Gamma (a+n+1)}{\Gamma (a+1) \Gamma (n+1)} \,
   _2F_1\left(\begin{smallmatrix}-n,a+b+n+1\\a+1\end{smallmatrix};\frac{1-z}{2}\right)\\
   	\, _2F_1(a,b;c;z)&=(1-z)^{-a} \, _2F_1\left(\begin{smallmatrix}a,c-b\\c\end{smallmatrix};\frac{z}{z-1}\right).
\end{align}
Applying these two formulas to (\ref{JacobiForm1}), we get
\begin{align}
&\mathcal{M}^{-1}(S^{j,n}_{m_1,m_4}(z))(x)\nonumber\\
=&\frac{c^j_{m_1}}{c^j_{m_4}}\frac{((-1)^{2j}+(-1)^{m_4-m_1})(1-\theta (\left| x\right| -1))}{2(1-x)^{\frac{m_1+m_4+1}{2}}x^{\frac{m_1-m_4+1}{2}}}P^{(m_4-m_1,-m_1-m_4)}_{j+m_1}\left(1-2x\right)
\end{align}
According to formula 9.43 in \cite{TableMellinTransform}, we can apply the Mellin transform $\mathcal{M}$ on the Jacobi polynomials,
\begin{align}
	&\mathcal{M}\left(\left\{\begin{smallmatrix}(b-x)^{\mu-1}P^{(\alpha,\beta)}_n(1-\gamma x) &x<b\\0&x>b\end{smallmatrix}\right.\right)(z)\nonumber\\
	=& \frac{\Gamma(\alpha+n+1)\Gamma(\mu)\Gamma(z)b^{z+\mu-1}}{n!\Gamma(1+\alpha)\Gamma(\mu+z)}{}_3F_2(\begin{smallmatrix}-n,1+n+\alpha+\beta,z\\1+\alpha,\mu+z\end{smallmatrix};\frac{1}{2}\gamma b)
\end{align}
After applying the Mellin transform, the matrix entries $S^{j,n}_{m_1,m_4}(z)$ can be transformed to
\begin{align}
&S^{j,n}_{m_1,m_4}(z)=\nonumber\\
&\frac{(-1)^{2 j}+(-1)^{m_4-m_1}}{2}\frac{c^j_{m_1}}{c^j_{m_4}}\frac{\Gamma \left(\frac{-m_1-m_4+1}{2}
  \right)  \Gamma \left(z+\frac{-m_1+m_4-1}{2}
 \right) }{\Gamma \left(-m_1+m_4+1\right) \Gamma \left(z-m_1\right)}\nonumber\\&\,
   _3F_2\left(\begin{smallmatrix}-j-m_1,j-m_1+1,z-\frac{m_1}{2}+\frac{m_4}{2}-\frac{1}{2}\\z-m_1,-m_1+m_4+1\end{smallmatrix};1\right
   ).
\end{align}
According to 07.27.17.0042.01 of \cite{WolframAlpha}, we have the following transformation property of $_3F_2$ value at 1:
\begin{align}
    {}_3F_2(-n,b,c;d,e;1)=\frac{(b)^{(n)} (-b-c+d+e)^{(n)} }{(d)^{(n)} (e)^{(n)}}{}_3F_2(\begin{smallmatrix}d-b,e-b,-n\\-b-c+d+e,-b-n+1\end{smallmatrix};1).
\end{align}
Apply this transform to the formula for $S^{j,n}_{m_1,m_4}(z)$,
\begin{align*}
&S^{j,n}_{m_1,m_4}(z)\\
  =&\frac{((-1)^{2 j}+(-1)^{m_4-m_1})\pi (2j)!}{2c^j_{m_1}c^j_{m_4}}\frac{\sec \left(\frac{m_1+m_4}{2} \pi  \right)\Gamma
   \left(\frac{2 z-m_1+m_4-1}{2} \right)}{ \Gamma \left(\frac{-2 j-m_1+m_4+1}{2}\right) \Gamma
   (j+z)}\\&\, _3F_2\left(\begin{smallmatrix}-j+z-1,-j-m_1,m_4-j\\-2
   j,-j-\frac{m_1}{2}+\frac{m_4}{2}+\frac{1}{2}\end{smallmatrix};1\right).
\end{align*}
We also need to consider that $m_1$ and $m_4$ satisfy the same parity condition, that $m_1,m_4\in \mathtt{M}(j,n;\delta_2,\delta_1)$, and therefore $\frac{m_1+m_4}{2}$ is an integer. Hence $S^{j,n}_{m_1,m_4}$ can be expressed in a simpler form in terms of hypergeometric function $_3F_2$:
\begin{align}
&S^{j,n}_{m_1,m_4}(z)=
  \frac{(-1)^{\frac{m_1+m_4}{2}} (2j)! \pi}{c^j_{m_1}c^j_{m_4}}\frac{ \Gamma
   \left(\frac{2 z-m_1+m_4-1}{2}\right)}{ \Gamma \left(\frac{-2 j-m_1+m_4+1}{2}\right) \Gamma
   (j+z)}\, _3F_2\left(\begin{smallmatrix}-j+z-1,-j-m_1,m_4-j\\-2
   j,-j-\frac{m_1}{2}+\frac{m_4}{2}+\frac{1}{2}\end{smallmatrix};1\right)\label{SExpression}
\end{align}
\subsection{Normalizations}
The paper \cite{KnappSteinIntertwining} discussed the normalization of intertwining operators in detail. There is a scalar valued meromorphic function $\eta(w,\chi_{\delta,\lambda})$ in $\lambda$ such that:
\begin{align}
	A(w_0,\chi_{\delta,\lambda})A(w_0^{-1},w_0\chi_{\delta,\lambda})=\eta(w_0,\chi_{\delta,\lambda}) I .\label{IntertwiningNorm1}
\end{align}
The meromorphic function $\eta(w_0,\chi_{\delta,\lambda})$ factor into rank-one factors:
\[
	\eta(w_0,\chi_{\delta,\lambda}) = \prod_{\substack{\alpha\text{ reduced}\\\langle\alpha_i,\alpha\rangle>0\text{ for positive }\alpha_i}}\eta_{\delta,\alpha}\left(\frac{\langle\lambda,\rho^{(\alpha)}\rangle}{\langle \rho^{(\alpha)},\rho^{(\alpha)}\rangle}\right).
\]
Moreover, these $\eta_{\delta,\alpha}(z)$'s factor into meromorphic functions $\gamma_{\delta,\alpha}(z)$:
\begin{align}
	\eta_{\delta,\alpha}(z) = \gamma_{\delta,\alpha}(z)\overline{\gamma_{\delta,\alpha}(-\bar{z})}.\label{IntertwiningNorm2}
\end{align}
In terms of the long intertwining operator $A(w_0,w_0\chi_{\delta,\lambda})$ corresponding to $w_0$, it is related to the adjoint of the intertwining operator corresponding to $w_0$:
\[
	A(w_0,w_0\chi_{\delta,\lambda})^* = A(w_0,\chi_{\delta,-\bar{\lambda}}).
\]
We can thus normalize the intertwining operator $A(w_0,\chi_{\delta,\lambda})$ by dividing it by the meromorphic factor:
\[
	\gamma(w_0,\chi_{\delta,\lambda}) =\prod_{\substack{\alpha\text{ reduced}\\\langle\alpha_i,\alpha\rangle>0\text{ for positive }\alpha_i}} \gamma_{\delta,\alpha}\left(\frac{\langle\lambda,\rho^{(\alpha)}\rangle}{\langle \rho^{(\alpha)},\rho^{(\alpha)}\rangle}\right).
\]
We denote the new normalized intertwining operator by \[A'(w_0,\chi_{\delta,\lambda})=\gamma(w_0,\chi_{\delta,\lambda})^{-1}A(w_0,\chi_{\delta,\lambda}).\]

In the case of $Sp(4,\RR)$ principal series, according to (\ref{IntertwiningNorm1}) and (\ref{IntertwiningNorm2}), there is a normalization factor $\gamma(w_{\alpha_1},\lambda)$ of the intertwining operator entries $S^{j,n}_{m_1,m_4}(z)$, such that the operator $A(w_{\alpha_1},\lambda)$ has the property 
\[
	A(w_{\alpha_1},w_{\alpha_1}\lambda)A(w_{\alpha_1},\lambda) = \gamma(w_{\alpha_1},\lambda)\gamma(w_{\alpha_1},w_{\alpha_1}\lambda) I.
\]
We define the normalized intertwining operator as $A'(w_{\alpha_1},\lambda) = \frac{1}{\gamma(w_{\alpha_1},\lambda)}A(w_{\alpha_1},\lambda)$. Since we have assumed that $(\delta_1,\delta_2) \in\{(0,0),(1,1)\}$, the $j$ and $n$ are integers. We can check that for $m_1,m_3\in \mathtt{M}(j,n;\delta_2,\delta_1)$,
\begin{align}
	&\sum_{m_2\in \mathtt{M}(j,n;\delta_2,\delta_1)}S^{j,n}_{m_1,m_2}(z)S^{j,n}_{m_2,m_3}(1-z)\nonumber\\ 
	=& \sum_{\substack{m_2\in \mathtt{M}(j,n;\delta_2,\delta_1)\\-j\leq m_4,m_5\leq j}}(-1)^{m_4+m_5}M^{j,n}_{m_1,m_4}N^{j,n}_{m_4,m_2}M^{j,n}_{m_2,m_5}N^{j,n}_{m_5,m_3}Q(z,m_4)Q(1-z,m_5).\label{QSum}
\end{align}
Since $\frac{1+(-1)^k}{2}$ is $0$ for $k\equiv 1\text{ mod }2$, and is $1$ for $k\equiv 0\text{ mod }2$, the sum on $m_2$ can be reduced to
\begin{align}
\sum_{m_2\in \mathtt{M}(j,n;\delta_2,\delta_1)}N^{j,n}_{m_4,m_2}M^{j,n}_{m_2,m_5} &= \sum_{-j\leq m_2\leq j}N^{j,n}_{m_4,m_2}M^{j,n}_{m_2,m_5} \frac{1+(-1)^{n-m_2-\delta_1}}
{2}\frac{1+(-1)^{n+m_2-\delta_2}}{2}\label{MNSum}
\end{align}
Taking into account the condition $\delta_1=\delta_2\in\{0,1\}$ and and its consequence that $m_3,n\in\ZZ$, we can expand the terms 
\begin{align*}
\frac{1+(-1)^{n-m_2-\delta_1}}
{2}\frac{1+(-1)^{n+m_2-\delta_2}}{2} =& \frac{1}{2} + \frac{(-1)^{n-\delta_1}}{2}\frac{(-1)^{m_2}+(-1)^{-m_2}}{2}\\
=&\frac{1}{2} + \frac{(-1)^{n-\delta_1}}{2}(-1)^{m_2}\\
=&\frac{1}{2} + \frac{(-1)^{n-\delta_1}}{2}e^{\ii m_2 \pi}.
\end{align*}
Since $W^{(j,n)}_{m_2,m_2} (0,-\pi,0,0)= e^{\ii m_2 \pi}$, by (\ref{Comm1}) and the definition (\ref{MMDef}) and (\ref{NNDef}) of $M^{j,n}_{m_2,m_5}$ and $N^{j,n}_{m_4,m_2}$, we have \[\sum_{m_2\in \mathtt{M}(j,n;\delta_2,\delta_1)}N^{j,n}_{m_4,m_2}M^{j,n}_{m_2,m_5}W^{(j,n)}_{m_2,m_2}(0,-\pi,0,0) = W^{(j,n)}_{m_4,m_5}(0,0,-\pi,0).\] 
The sum (\ref{MNSum}) can thus be written as
\begin{align*}
\sum_{m_2\in \mathtt{M}(j,n;\delta_2,\delta_1)}N^{j,n}_{m_4,m_2}M^{j,n}_{m_2,m_5} &= \frac{1}{2}\delta_{m_4,m_5} + \frac{(-1)^{n-\delta_1}}{2} W^{(j,n)}_{m_4,m_5}(e^{\pi U_2})\\
&= \frac{1}{2}\delta_{m_4,m_5} + \frac{(-1)^{j-m_4-\delta_1+n}}{2} \delta_{-m_4,m_5},
\end{align*}
and the sum (\ref{QSum}) becomes
\begin{align}
	&\sum_{m_2\in \mathtt{M}(j,n;\delta_2,\delta_1)}S^{j,n}_{m_1,m_2}(z)S^{j,n}_{m_2,m_3}(1-z)\nonumber\\ 
	=& \sum_{\substack{-j\leq m_4,m_5\leq j}}(-1)^{m_4+m_5}M^{j,n}_{m_1,m_4}N^{j,n}_{m_5,m_3}Q(z,m_4)Q(1-z,m_5)\frac{\delta_{m_4,m_5}}{2}+\nonumber \\
	&\sum_{\substack{-j\leq m_4,m_5\leq j}}(-1)^{m_4+m_5}M^{j,n}_{m_1,m_4}N^{j,n}_{m_5,m_3}Q(z,m_4)Q(1-z,m_5) \frac{(-1)^{j-m_4-\delta_1+n}}{2} \delta_{-m_4,m_5}\nonumber\\
	=& \frac{1}{2}\sum_{\substack{-j\leq m_4\leq j}}M^{j,n}_{m_1,m_4}N^{j,n}_{m_4,m_3}Q(z,m_4)Q(1-z,m_4)+ \frac{1}{2}\sum_{\substack{-j\leq m_4\leq j}}M^{j,n}_{m_1,m_4}N^{j,n}_{-m_4,m_3}\nonumber\\
	&(-1)^{j-m_4-\delta_1+n}Q(z,m_4)Q(1-z,-m_4)\label{ProductOfS}.
\end{align}
We observe that the product of the $Q$-functions are simply products of $\Gamma$-functions:
\begin{equation}
Q(z,n)Q(1-z,n)=Q(z,n)Q(1-z,-n)=\pi \frac{\Gamma(z-1/2)}{\Gamma(z)}\frac{\Gamma(-z+1/2)}{\Gamma(1-z)},
\end{equation}
hence the formula (\ref{ProductOfS}) is in fact
\begin{align}
	&\sum_{m_2\in \mathtt{M}(j,n;\delta_2,\delta_1)}S^{j,n}_{m_1,m_2}(z)S^{j,n}_{m_2,m_3}(1-z)\nonumber\\ 
	=& \frac{\pi}{2}\frac{\Gamma(z-1/2)}{\Gamma(z)}\frac{\Gamma(-z+1/2)}{\Gamma(1-z)}\sum_{\substack{-j\leq m_4\leq j}}\left(M^{j,n}_{m_1,m_4}N^{j,n}_{m_4,m_3}+(-1)^{j-m_4-\delta_1+n}M^{j,n}_{m_1,m_4}N^{j,n}_{-m_4,m_3}\right)\nonumber\\
	=&\frac{\pi}{2}\frac{\Gamma(z-1/2)}{\Gamma(z)}\frac{\Gamma(-z+1/2)}{\Gamma(1-z)}(\delta_{m_1,m_3}+(-1)^{-m_1-\delta_1+n}\delta_{m_1,m_3})\nonumber\\
	=&\pi\frac{\Gamma(z-1/2)}{\Gamma(z)}\frac{\Gamma(-z+1/2)}{\Gamma(1-z)}\frac{1+(-1)^{-m_1-\delta_1+n}}{2}\delta_{m_1,m_3}\label{NormalizeS}.
\end{align}
Therefore, it is natural to consider be the action of the intertwining operators on the $K$-types with $j=0$:
\begin{align}
	S^{(0,n)}_{0,0}(z) = \frac{\sqrt{\pi}\Gamma\left(z-1/2\right)}{\Gamma(z)}
\end{align}
and set the normalization factor $\gamma(w_{\alpha},\lambda)$ to be \[\gamma(w_\alpha,\lambda) = S^{(0,n)}_{0,0}(\frac{\langle\check{\alpha},\lambda+\rho\rangle}{2}).\]
The matrix entries of the normalized intertwining operator $A'(w_{\alpha_1},\lambda)$ are given by the function 
\begin{align}
\mathcal{S}^{j,n}_{m_1,m_2} (z):= \frac{S^{j,n}_{m_1,m_2}(z)}{S^{(0,n)}_{0,0}(z)} = \left(\frac{\sqrt{\pi}\Gamma\left(z-1/2\right)}{\Gamma(z)}\right)^{-1}S^{j,n}_{m_1,m_2}(z),
\end{align}
and by (\ref{SExpression}) we can further simplify the expression of the matrix entries in terms of Pochhammer symbols:
\begin{align}
\mathcal{S}^{j,n}_{m_1,m_2}(z)
  =&\frac{(-1)^{\frac{m_1+m_2}{2}} (2j)! \sqrt{\pi}}{c^j_{m_1}c^j_{m_2}\Gamma\left(\frac{1-2j-m_1+m_2}{2}\right)}\frac{(z-1/2)^{(-\frac{m_1-m_2}{2})}}{(z)^{(j)}}\, _3F_2\left(\begin{smallmatrix}-j+z-1,-j-m_1,-j+m_2\\-2
   j,-j-\frac{m_1-m_2}{2}+\frac{1}{2}\end{smallmatrix};1\right).
\end{align}
According to (\ref{NormalizeS}), if $m_1,m_3\in \mathtt{M}(j,n;\delta_1,\delta_2)$, the normalized intertwining operator entries $\mathcal{S}^{j,n}_{m_1,m_2}(z)$ satisfy the property
\begin{align}
	\sum_{m_2\in \mathtt{M}(j,n;\delta_2,\delta_1)}\mathcal{S}^{j,n}_{m_1,m_2}(z)\mathcal{S}^{j,n}_{m_2,m_3}(1-z) = \delta_{m_1,m_3}.
\end{align}
%From the generating function of Hahn polynomial, we have:
%\begin{align*}
%	\sum_{m_4=-j}^{j}&\frac{2\sqrt{(j-m_1)!(j+m_1)!(j-m_2)!(j+m_2)!}}{(-1)^{2j}+(-1)^{m_4-m_1}}\\
%    &\left(-\frac{\sin (2 \pi  j) \cos \left(\frac{1}{2} \pi  \left(m_1+m_4\right)\right) \Gamma \left(-j-m_4\right) \Gamma
%   \left(\frac{1}{2} \left(-2 j-m_1+m_4+1\right)\right){}^2 \Gamma (j+z)}{\pi ^2 \Gamma \left(\frac{1}{2}
%   \left(-m_1-m_4+1\right)\right) \Gamma \left(j-m_4+1\right) \Gamma \left(z+\frac{1}{2} \left(-m_1+m_4-1\right)\right)}\right)S_{m_1,m}(z)q^{j-m}\\
%   &\, _1F_1\left(-j-z+1;\frac{1}{2} \left(-2 j-m_1+m_4+1\right);q\right) \, _1F_1\left(-j+z-1;\frac{1}{2} \left(-2
%   j-m_1+m_4+1\right);-q\right)
%\end{align*}
%\subsection{Generating Functions of $\mathcal{S}^{j,n}_{m_1,m_2}(z)$}
If $n\geq 0$, the $n$-th Laurent series coefficient of the function in terms of $t$ is the $_3F_2$ hypergeometric function (see \cite[Page 94]{GeneratingFunctions}):
\begin{align}
	\left[(1-t)^{-\lambda}{}_2F_1\left(\begin{smallmatrix}
		a,b\\c
	\end{smallmatrix};zt\right)\right]_n = \frac{\Gamma(\lambda+n)}{\Gamma(n+1)\Gamma(\lambda)}{}_3F_2\left(\begin{smallmatrix}
		-n,a,b\\c,1-n-\lambda
	\end{smallmatrix};z\right).
\end{align}
Therefore, $\mathcal{S}^{j,n}_{m_1,m_2}(z)$ is the constant Laurent series coefficient of either of the following two functions:
\begin{align}
	H^j_{m_1,m_2}(z;t)&=\frac{(2j)!(j+m_1)!}{c^j_{m_1}c^j_{m_2}\sqrt{\pi}}\Gamma\left(\frac{1-m_1-m_2}{2}\right)\frac{(z-1/2)^{(-\frac{m_1-m_2}{2})}}{(z)^{(j)}}\times\nonumber\\
    &\frac{(1-t)^{\frac{-1+m_1+m_2}{2}}}{(-t)^{j+m_1}}{}_2F_1\left(\begin{smallmatrix}
	-j+m_2,-1-j+z\\-2j
	\end{smallmatrix};t\right)\label{HDef}\\
    G^j_{m_1,m_2}(z;t)&=\frac{(2j)!(j-m_2)!}{c^j_{m_1}c^j_{m_2}\sqrt{\pi}}\Gamma\left(\frac{1+m_1+m_2}{2}\right)\frac{(z-1/2)^{(-\frac{m_1-m_2}{2})}}{(z)^{(j)}}\times\nonumber\\
    &\frac{(1-t)^{\frac{-1-m_1-m_2}{2}}}{(-t)^{j-m_2}}{}_2F_1\left(\begin{smallmatrix}
	-j-m_1,-1-j+z\\-2j
	\end{smallmatrix};t\right).\label{GDef}
\end{align}

\subsection{Computing the Product Matrix}
In this section we compute the product of the four matrices:
\[
	A'(w_0,\chi_{\delta,\lambda})=A_4(\lambda)\cdot A_3(\lambda)\cdot A_2(\lambda)\cdot A_1(\lambda)
\]
when $\delta = (0,0)$ or $(1,1)$. We would also like to normalize each individual matrices $A_i(\lambda)$ such that they satisfy the condition:
\[
	A_i(-\lambda)A_i(\lambda) = I.
\]
Recalling the parity conditions (\ref{KCond}) and (\ref{MCond}) satisfied by $j,n,m_1,m_2$, since $(\delta_1,\delta_2) \in\{ (0,0),(1,1)\}$, the pair $(j,n)$ are integers, and $m_1,m_2$ satisfy
\[
	\begin{smallmatrix}n-m_i+\delta_1\equiv 0\\n+m_i+\delta_2\equiv 0\end{smallmatrix}\text{ }\mathrm{mod}\text{ } 2.
\]
Despite the difficulty of calculating the product of the 4 matrices, a lot of terms in the sum can be reduced using the parity condition. We replace the function $T^n_{m_1}(z)$ by the normalized function defined by
\[
  \mathcal{T}^n_{m_1}(z)=\frac{\ii^{-n+m_1}}{(z)^{(\frac{m_1-n}{2})}(z)^{(\frac{-m_1+n}{2})}}.
\] The matrix entry $[A'(\lambda)]^{j,n}_{m_1,m_2}$ of the normalized intertwining operator can be expressed as the constant term of the Laurent series in $t_1,t_2$ of the sum:
\begin{align}
	[A(\lambda)]^{j,n}_{m_1,m_2}(t_1,t_2) &= \mathcal{T}^{n}_{m_2}(\frac{\lambda_2+1}{2})\times\nonumber\\
    &\sum_{m_3\in\mathtt{M}(j,n;\delta_1,\delta_2)}G^{j}_{m_2,m_3}(\frac{\lambda_1+\lambda_2+1}{2},t_2)\mathcal{T}^{n}_{m_3}(\frac{\lambda_1+1}{2})H^{j}_{m_3,m_1}(\frac{\lambda_1-\lambda_2+1}{2};t_1).\label{ASum}
\end{align}
Plugging (\ref{HDef}) and (\ref{GDef}) into (\ref{ASum}), we have
\begin{align}
	&[A(\lambda)]^{j,n}_{m_1,m_2}(t_1,t_2)\nonumber\\
 =&\frac{((2j)!)^2}{c^j_{m_1}c^j_{m_2}\pi} \frac{\ii^{-n+m_2}}{\left(\frac{\lambda_2+1}{2}\right)^{(\frac{m_2-n}{2})}\left(\frac{\lambda_2+1}{2}\right)^{(\frac{-m_2+n}{2})}}\frac{1}{(\frac{\lambda_1-\lambda_2+1}{2})^{(j)}(\frac{\lambda_1+\lambda_2+1}{2})^{(j)}}\times\nonumber\\
  &{}_2F_1\left(\begin{smallmatrix}
	-j+m_1,-1-j+\frac{\lambda_1-\lambda_2+1}{2}\\-2j
	\end{smallmatrix};t_1\right)
	{}_2F_1\left(\begin{smallmatrix}
	-j-m_2,-1-j+\frac{\lambda_1+\lambda_2+1}{2}\\-2j
	\end{smallmatrix};t_2\right)\times\nonumber\\
 &\sum_{m_3\in\mathtt{M}(j,n;\delta_1,\delta_2)}\frac{\ii^{-n+m_3}}{\left(\frac{\lambda_1+1}{2}\right)^{(\frac{m_3-n}{2})}\left(\frac{\lambda_1+1}{2}\right)^{(\frac{-m_3+n}{2})}}\Gamma\left(\frac{1+m_3+m_2}{2}\right)\Gamma\left(\frac{1-m_1-m_3}{2}\right)\times\nonumber\\
 &\left(\frac{\lambda_1-\lambda_2}{2}\right)^{(\frac{m_1-m_3}{2})}\left(\frac{\lambda_1+\lambda_2}{2}\right)^{(\frac{m_3-m_2}{2})}\nonumber\frac{(1-t_2)^{\frac{-1-m_3-m_2}{2}}}{(-t_2)^{j-m_3}}\frac{(1-t_1)^{\frac{-1+m_1+m_3}{2}}}{(-t_1)^{j+m_3}}.\nonumber
\end{align}
Based on the parity of $j-m_3$, we can separate the calculation into two cases:
\subsubsection{Summation when $j-n\equiv \delta_i\text{ mod }2$}

In the situation that $j-n\equiv \delta_i\text{ }\mathrm{mod}\text{ } 2$, the set of $m_3\in \mathtt{M}(j,n;\delta_2,\delta_1)$'s are:
\[
	m_3=j-2p\text{ where } p\in\{0,1,2,\ldots,j\}.
\]
Thus we can replace $m_3$ by $j-2p$, and the sum reduces to a sum of finitely many Pochhammer symbols:
\begin{align}
	&[A(\lambda)]^{j,n}_{m_1,m_2}(t_1,t_2) =\frac{((2j)!)^2}{c_{m_1}^j c_{m_2}^j\pi}\frac{\ii^{-n+m_2}}{\left(\frac{\lambda_2+1}{2}\right)^{(\frac{m_2-n}{2})}\left(\frac{\lambda_2+1}{2}\right)^{(-\frac{m_2-n}{2})}}\times\nonumber\\
	&{}_2F_1\left(\begin{smallmatrix}-j+m_1,-j+\frac{\lambda_1-\lambda-1}{2}\\-2j\end{smallmatrix};t_1\right)
	{}_2F_1\left(\begin{smallmatrix}-j-m_2,-j+\frac{\lambda_1+\lambda-1}{2}\\-2j\end{smallmatrix};t_2\right)\times\nonumber\\
	&\sum_{p=0}^j\frac{\ii^{-n+j-2p}}{\left(\frac{\lambda_1+1}{2}\right)^{(\frac{j-n-2p}{2})}\left(\frac{\lambda_1+1}{2}\right)^{(-\frac{j-n-2p}{2})}}
	\Gamma\left(\frac{1+j+m_2}{2}-p\right)\Gamma\left(\frac{1-j-m_1}{2}+p\right)\times\nonumber\\
	&\left(\frac{\lambda_1-\lambda_2}{2}\right)^{(\frac{-j+m_1}{2}+p)}\left(\frac{\lambda_1+\lambda_2}{2}\right)^{(\frac{j-m_2}{2}-p)} (1-t_1)^{\frac{-1+j+m_1}{2}-p} (1-t_2)^{\frac{-1-j-m_2}{2}+p}t_1^{-2j+2p}t_2^{-2p}.\label{Sp4IntertwiningEven}
\end{align}
By definition of Pochhammer symbols, we have the formula
\begin{equation}
	(x)^{(l\pm p)} = (x)^{(l)} (x+l)^{(\pm p)} \label{Rule1}
\end{equation}
which we can apply to the Pochhammer symbols in the summation (\ref{Sp4IntertwiningEven}) to single out the $\pm p$'s in the Pochhammer exponents. Also, by applying the formula
\begin{equation}
	(x)^{(-p)}(1-x)^{(p)} = (-1)^p\text{ for }p\in\ZZ\label{Rule2}
\end{equation}
to the Pochhammer symbols with $-p$ as exponent, we can manage to change all the exponents of the Pochhammer symbols in (\ref{Sp4IntertwiningEven}) to $(p)$. The explicit set of rules is the following
\begin{align}
	\left(\frac{\lambda_1+1}{2}\right)^{(\frac{j-n-2p}{2})}&\rightarrow \left(\frac{\lambda_1+1}{2}\right)^{(\frac{j-n}{2})} \left(\frac{\lambda_1+1+j-n}{2}\right)^{(-p)}\nonumber\\
	&\rightarrow (-1)^p \left(\frac{\lambda_1+1}{2}\right)^{(\frac{j-n}{2})} \left(1-\frac{\lambda_1+1+j-n}{2}\right)^{(p)}\label{Rule001}\\
	\left(\frac{\lambda_1+1}{2}\right)^{(-\frac{j-n-2p}{2})}&\rightarrow \left(\frac{\lambda_1+1}{2}\right)^{(-\frac{j-n}{2})} \left(\frac{\lambda_1+1-j+n}{2}\right)^{(p)}
\end{align}
\begin{align}
	\Gamma\left(\frac{1+j+m_2}{2}-p\right)&\rightarrow 	\Gamma\left(\frac{1+j+m_2}{2}\right) 	\left(\frac{1+j+m_2}{2}\right)^{(-p)}\nonumber\\
	&\rightarrow (-1)^p \Gamma\left(\frac{1+j+m_2}{2}\right) 	\left(1-\frac{1+j+m_2}{2}\right)^{(p)}\\
	\Gamma\left(\frac{1-j-m_1}{2}+p\right)&\rightarrow 	\Gamma\left(\frac{1-j-m_1}{2}\right) 	\left(\frac{1-j-m_1}{2}\right)^{(p)}
\end{align}
\begin{align}
	\left(\frac{\lambda_1-\lambda_2}{2}\right)^{(-\frac{j-m_1-2p}{2})}&\rightarrow \left(\frac{\lambda_1-\lambda_2}{2}\right)^{(-\frac{j-m_1}{2})} \left(\frac{\lambda_1-\lambda_2}{2}-\frac{j-m_1}{2}\right)^{(p)}\\
	\left(\frac{\lambda_1+\lambda_2}{2}\right)^{(\frac{j-m_2-2p}{2})}&\rightarrow \left(\frac{\lambda_1+\lambda_2}{2}\right)^{(\frac{j-m_2}{2})} \left(\frac{\lambda_1+\lambda_2}{2}+\frac{j-m_2}{2}\right)^{(-p)}\nonumber\\
	&\rightarrow (-1)^p\left(\frac{\lambda_1+\lambda_2}{2}\right)^{(\frac{j-m_2}{2})} \left(1-\frac{\lambda_1+\lambda_2}{2}+\frac{j-m_2}{2}\right)^{(p)}.\label{Rule002}
\end{align}
After applying (\ref{Rule001})-(\ref{Rule002}) and bringing out the factors not involved in the summation, we have
\begin{align}
	&[A(\lambda)]^{j,n}_{m_1,m_2}(t_1,t_2) =\nonumber\\
	&\frac{\ii^{j+m_2-2n}((2j)!)^2}{\pi c^j_{m_1}c^j_{m_2}}\frac{\left(\frac{\lambda_1-\lambda_2}{2}\right)^{(\frac{-j+m_1}{2})}\left(\frac{\lambda_1+\lambda_2}{2}\right)^{(\frac{j-m_2}{2})}}{\left(\frac{\lambda_1+1}{2}\right)^{(\frac{-j+n}{2})}\left(\frac{\lambda_1+1}{2}\right)^{(\frac{j-n}{2})}}\times\nonumber\\&\frac{\Gamma\left(\frac{1+j+m_2}{2}\right)\Gamma\left(\frac{1-j-m_1}{2}\right)}{\left(\frac{\lambda_1+\lambda_2+1}{2}\right)^{(j)}\left(\frac{\lambda_1-\lambda_2+1}{2}\right)^{(j)}\left(\frac{\lambda_2+1}{2}\right)^{(\frac{n-m_2}{2})}\left(\frac{\lambda_2+1}{2}\right)^{(\frac{-n+m_2}{2})}}\times\nonumber\\
	&{}_2F_1\left(\begin{smallmatrix}-j+m_1,\frac{\lambda_1-\lambda_2-2j-1}{2}\\-2j\end{smallmatrix};t_1\right)
	{}_2F_1\left(\begin{smallmatrix}-j-m_2,\frac{\lambda_1+\lambda_2-2j-1}{2}\\-2j\end{smallmatrix};t_2\right)\times\nonumber\\&\gamma'^{j,n}_{m_1,m_2}(\lambda; t_1,t_2)\label{Sp4IntertwiningBeforeTransformation},
\end{align}
where the function $\gamma'^{j,n}_{m_1,m_2}(\lambda; t_1,t_2)$ in two variables $t_1,t_2$ is defined as the hypergeometric sum
\begin{align}
	\gamma'^{j,n}_{m_1,m_2}(\lambda; t_1,t_2)&=(1-t_1)^{\frac{-1+j+m_1}{2}}(-t_1)^{-2j}(1-t_2)^{\frac{-1-j-m_2}{2}}\times\nonumber\\
	&\sum_{p=0}^{j}\frac{\left(\frac{1-j-m_1}{2}\right)^{(p)}\left(\frac{-j+m_1+\lambda_1-\lambda_2}{2}\right)^{(p)}\left(\frac{1-j+n-\lambda_1}{2}\right)^{(p)}(1)^{(p)}}{\left(\frac{1-j-m_2}{2}\right)^{(p)}\left(1-\frac{j-m_2+\lambda_1+\lambda_2}{2}\right)^{(p)}\left(\frac{1-j+n+\lambda_1}{2}\right)^{(p)}} \frac{1}{p!}\left(\frac{t_1^2(1-t_2)}{t_2^2(1-t_1)}\right)^p\label{thesum}.
\end{align}
We can use a partial sum formula for hypergeometric series from 16.2.4 of \cite{NIST}:
\begin{align}
	\sum_{k=0}^m\frac{(\mathbf{a})^{(k)}}{(\mathbf{b})^{(k)}}\frac{z^k}{k!} = \frac{(\mathbf{a})^{(m)} z^m}{(\mathbf{b})^{(m)} m!}{}_{q+2}F_p\left(\begin{smallmatrix}
	-m,1,1-m-\mathbf{b}\\1-m-\mathbf{a}
	\end{smallmatrix};\frac{(-1)^{p+q+1}}{z}\right)\label{NISTPARTIAL},
\end{align}
where $\mathbf{a} = (a_1,\ldots,a_p)$ and $\mathbf{b} = (b_1,\ldots,b_q)$ are $p$ and $q$ dimensional vectors of complex numbers, and the Pochhammer symbol $(\mathbf{a})^{(k)}$ and $(\mathbf{b})^{(k)}$ for these two vectors are defined as the product $(\mathbf{a})^{(k)} = (a_1)^{(k)}\ldots (a_p)^{(k)}$ and $(\mathbf{b})^{(k)} = (b_1)^{(k)}\ldots (b_q)^{(k)}$.
Then the summation term in $\gamma'^{j,n}_{m_1,m_2}(\lambda; t_1,t_2)$ can be simplified into a single compact form hypergeometric function ${}_5F_4$:
\begin{align*}
	&\sum_{p=0}^{j}\frac{\left(\frac{1-j-m_1}{2}\right)^{(p)}\left(\frac{-j+m_1+\lambda_1-\lambda_2}{2}\right)^{(p)}\left(\frac{1-j+n-\lambda_1}{2}\right)^{(p)}(1)^{(p)}}{\left(\frac{1-j-m_2}{2}\right)^{(p)}\left(1-\frac{j-m_2+\lambda_1+\lambda_2}{2}\right)^{(p)}\left(\frac{1-j+n+\lambda_1}{2}\right)^{(p)}} \frac{1}{p!}\left(\frac{t_1^2(1-t_2)}{t_2^2(1-t_1)}\right)^p\\
   =&\frac{ \left(\frac{1-j-m_1}{2}\right)^{(j)} \left(\frac{-j+m_1+\lambda _1-\lambda _2}{2}\right)
  ^{(j)} \left(\frac{-j+n-\lambda _1+1}{2}\right)^{(j)}}{\left(\frac{1-j-m_2}{2}\right)^{(j)} \left(\frac{-j+m_2-\lambda _1-\lambda _2}{2}
  +1\right)^{(j)} \left(\frac{-j+n+\lambda _1+1}{2}
  \right)^{(j)}}\left(\frac{ t_1^2\left(1-t_2\right)}{t_2^2 \left(1-t_1\right)}\right)^j\times\\
  &
   \,
   _5F_4\left(\begin{matrix}-j,1,\frac{1-j+m_2}{2},\frac{-j-n-\lambda_1+1}{2},\frac{-j-m_2+\lambda_1+\lambda_2}{2}\\\frac{1-j+m_1}{2},\frac{1-j-n+\lambda_1}{2},\frac{-j-m_1-\lambda_1+\lambda_2}{2}+1,-j\end{matrix};\frac{t_2^2 \left(1-t_1\right)}{t_1^2\left(1-t_2\right)
   }\right)%\\
%   =&\frac{ \left(\frac{1-j-m_1}{2}\right)^{(j)} \left(\frac{-j+m_1+\lambda _1-\lambda _2}{2}\right)
%  ^{(j)} \left(\frac{-j+n-\lambda _1+1}{2}\right)^{(j)}}{\left(\frac{1-j-m_2}{2}\right)^{(j)} \left(\frac{-j+m_2-\lambda _1-\lambda _2}{2}
%  +1\right)^{(j)} \left(\frac{-j+n+\lambda _1+1}{2}
%  \right)^{(j)}}\left(\frac{ t_1^2\left(1-t_2\right)}{t_2^2 \left(1-t_1\right)}\right)^j\times\\
%  &
%   \,
%   _4F_3\left(\begin{matrix}1,\frac{1-j+m_2}{2},\frac{-j-n-\lambda_1+1}{2},\frac{-j-m_2+\lambda_1+\lambda_2}{2}\\\frac{1-j+m_1}{2},\frac{1-j-n+\lambda_1}{2},\frac{-j-m_1-\lambda_1+\lambda_2}{2}+1\end{matrix};\frac{t_2^2 \left(1-t_1\right)}{t_1^2\left(1-t_2\right)
%   }\right),
\end{align*}
and if we replace the sum in (\ref{thesum}) by the expression above involving $_5F_4$, also we recover the Pochhammer symbols from the transforms (\ref{Rule001})-(\ref{Rule002}) based on (\ref{Rule1}) and (\ref{Rule2}), the original function (\ref{Sp4IntertwiningEven}) can be expressed as
\begin{align}
	&[A(\lambda)]^{j,n}_{m_1,m_2}(t_1,t_2)=\nonumber\\
	&\frac{(-1)^{n}((2j)!)^2}{\left(\frac{\lambda_1+1}{2}\right)^{(\frac{j+n}{2})}\left(\frac{\lambda_1+1}{2}\right)^{(-\frac{j+n}{2})}\left(\frac{\lambda_1-\lambda_2+1}{2}\right)^{(j)}\left(\frac{\lambda_1+\lambda_2+1}{2}\right)^{(j)}}\times\nonumber\\
	&\frac{1}{c^j_{m_1}c^j_{m_2}}\frac{\left(\frac{\lambda_1-\lambda_2}{2}\right)^{(\frac{j+m_1}{2})}\left(\frac{\lambda_1+\lambda_2}{2}\right)^{(\frac{-j-m_2}{2})}}{\left(\frac{\lambda_2+1}{2}\right)^{(\frac{m_2-n}{2})}\left(\frac{\lambda_2+1}{2}\right)^{(-\frac{m_2-n}{2})}}\frac{\Gamma\left(\frac{1+j-m_1}{2}\right)}{\Gamma\left(\frac{1+j-m_2}{2}\right)}\times\nonumber\\
	&(1-t_1)^{\frac{-1-j+m_1}{2}}(1-t_2)^{\frac{-1+j-m_2}{2}}t_2^{-2j}\times\nonumber\\
	&{}_2F_1\left(\begin{smallmatrix}-j+m_1,\frac{\lambda_1-\lambda_2-2j-1}{2}\\-2j\end{smallmatrix};t_1\right){}_2F_1\left(\begin{smallmatrix}-j-m_2,\frac{\lambda_1+\lambda_2-2j-1}{2}\\-2j\end{smallmatrix};t_2\right)\times\nonumber\\
	 &\,
   _5F_4\left(\begin{smallmatrix}-j,1,\frac{-j+m_2+1}{2},\frac{-j-n-\lambda_1+1}{2},\frac{-j-m_2+\lambda_1+\lambda_2}{2}\\-j,\frac{-j+m_1+1}{2},\frac{-j-n+\lambda_1+1}{2},\frac{-j-m_1-\lambda_1+\lambda_2}{2}+1\end{smallmatrix};\frac{t_2^2 \left(1-t_1\right)}{ t_1^2\left(1-t_2\right)}\right)\label{evencase}.
\end{align}
The $-j$'s as parameters to the function $_5F_4$ don't necessarily cancel since they are non-positive. They play an important role in making the function $_5F_4$ rational.
\subsubsection{Summation when $j-n\not\equiv \delta_i\text{ mod }2$}
On the other hand, if $j-n\not\equiv \delta_i\text{ }\mathrm{mod}\text{ } 2$, 
\[
	m = j-2p-1\text{ where }p\in\{0,1,\ldots,j-1\},
\]
we have
\begin{align}
	&[A(\lambda)]^{j,n}_{m_1,m_2}(t_1,t_2) =\frac{((2j)!)^2}{c_{m_1}^j c_{m_2}^j\pi}\frac{\ii^{-n+m_2}}{\left(\frac{\lambda_2+1}{2}\right)^{(\frac{m_2-n}{2})}\left(\frac{\lambda_2+1}{2}\right)^{(-\frac{m_2-n}{2})}}\times\nonumber\\
	&{}_2F_1\left(\begin{smallmatrix}-j+m_1,-j+\frac{\lambda_1-\lambda-1}{2}\\-2j\end{smallmatrix};t_1\right)
	{}_2F_1\left(\begin{smallmatrix}-j-m_2,-j+\frac{\lambda_1+\lambda-1}{2}\\-2j\end{smallmatrix};t_2\right)\times\nonumber\\
	&\sum_{p=0}^{j-1}\frac{\ii^{-n+j-2p-1}}{\left(\frac{\lambda_1+1}{2}\right)^{(\frac{j-n-2p-1}{2})}\left(\frac{\lambda_1+1}{2}\right)^{(-\frac{j-n-2p-1}{2})}}
	\Gamma\left(\frac{j+m_2}{2}-p\right)\Gamma\left(\frac{2-j-m_1}{2}+p\right)\times\nonumber\\
	&\left(\frac{\lambda_1-\lambda_2}{2}\right)^{(\frac{1-j+m_1}{2}+p)}\left(\frac{\lambda_1+\lambda_2}{2}\right)^{(\frac{j-m_2-1}{2}-p)} (1-t_1)^{\frac{-2+j+m_1}{2}-p}\times\nonumber\\& (1-t_2)^{\frac{-j-m_2}{2}+p}t_1^{-2j+2p+1}t_2^{-2p-1}.\label{Sp4IntertwiningOdd}
\end{align}
Comparing with the even case, after applying (\ref{Rule1}) and $(\ref{Rule2})$, the summation will change to
\begin{align}
	&[A(\lambda)]^{j,n}_{m_1,m_2}(t_1,t_2) =\nonumber\\
	&\frac{\ii^{j+m_2-2n}((2j)!)^2}{\pi c^j_{m_1}c^j_{m_2}}\frac{\left(\frac{\lambda_1-\lambda_2}{2}\right)^{(\frac{-j+m_1+1}{2})}\left(\frac{\lambda_1+\lambda_2}{2}\right)^{(\frac{j-m_2-1}{2})}}{\left(\frac{\lambda_1+1}{2}\right)^{(\frac{-j+n+1}{2})}\left(\frac{\lambda_1+1}{2}\right)^{(\frac{j-n-1}{2})}}\times\nonumber\\&\frac{\Gamma\left(\frac{j+m_2}{2}\right)\Gamma\left(\frac{2-j-m_1}{2}\right)}{\left(\frac{\lambda_1+\lambda_2+1}{2}\right)^{(j)}\left(\frac{\lambda_1-\lambda_2+1}{2}\right)^{(j)}\left(\frac{\lambda_2+1}{2}\right)^{(\frac{n-m_2}{2})}\left(\frac{\lambda_2+1}{2}\right)^{(\frac{-n+m_2}{2})}}\times\nonumber\\
	&{}_2F_1\left(\begin{smallmatrix}-j+m_1,\frac{\lambda_1-\lambda_2-2j-1}{2}\\-2j\end{smallmatrix};t_1\right)
	{}_2F_1\left(\begin{smallmatrix}-j-m_2,\frac{\lambda_1+\lambda_2-2j-1}{2}\\-2j\end{smallmatrix};t_2\right)\times\nonumber\\&\gamma'^{j,n}_{m_1,m_2}(\lambda; t_1,t_2)\label{Sp4IntertwiningBeforeTransformation2},
\end{align}
where
\begin{align}
	&\gamma''^{j,n}_{m_1,m_2}(\lambda; t_1,t_2)=(1-t_1)^{\frac{j+m_1-2}{2}}t_2^{-2j+1}t_1^{-1}(1-t_2)^{\frac{-j-m_2}{2}}\times\nonumber\\
	&{}_2F_1\left(\begin{smallmatrix}-j+m_1,\frac{\lambda_1-\lambda_2-2j-1}{2}\\-2j\end{smallmatrix};t_1\right)
	{}_2F_1\left(\begin{smallmatrix}-j-m_2,\frac{\lambda_1+\lambda_2-2j-1}{2}\\-2j\end{smallmatrix};t_2\right)\times\nonumber\\
	&\sum_{p=0}^{j-1}\frac{\left(\frac{2-j-m_1}{2}\right)^{(p)}\left(\frac{-j+m_2+\lambda_1-\lambda_2}{2}\right)^{(p)}\left(1-\frac{j-n+\lambda_1}{2}\right)^{(p)}(1)^{(p)}}{\left(\frac{2-j-m_2}{2}\right)^{(p)}\left(\frac{3}{2}-\frac{j-m_2-\lambda_1-\lambda_2}{2}\right)^{(p)}\left(\frac{2-j+n+\lambda_1}{2}\right)^{(p)}} \frac{1}{p!}\left(\frac{t_1^2(1-t_2)}{t_2^2(1-t_1)}\right)^p\label{thesum}.
\end{align}
By the same formula (\ref{NISTPARTIAL}) and using the relation between Pochhammer symbols and $\Gamma$-functions $(a)^{(n)} = \Gamma(a+n)/\Gamma(a)$ to absorb the extra $\Gamma$ factors in $\gamma''^{j,n}_{m_1,m_2}$, the matrix entries for the long intertwining operator has an expression
\begin{align}
	&[A(\lambda)]^{j,n}_{m_1,m_2}(t_1,t_2)=\nonumber\\
	&\frac{(-1)^{n}((2j)!)^2}{\left(\frac{\lambda_1+1}{2}\right)^{(\frac{j+n-1}{2})}\left(\frac{\lambda_1+1}{2}\right)^{(-\frac{j+n-1}{2})}\left(\frac{\lambda_1-\lambda_2+1}{2}\right)^{(j)}\left(\frac{\lambda_1+\lambda_2+1}{2}\right)^{(j)}}\times\nonumber\\
	&\frac{1}{c^j_{m_1}c^j_{m_2}}\frac{\left(\frac{\lambda_1-\lambda_2}{2}\right)^{(\frac{j+m_1-1}{2})}\left(\frac{\lambda_1+\lambda_2}{2}\right)^{(\frac{-j-m_2+1}{2})}}{\left(\frac{\lambda_2+1}{2}\right)^{\frac{m_2-n}{2}}\left(\frac{\lambda_2+1}{2}\right)^{-\frac{m_2-n}{2}}}\frac{\Gamma\left(\frac{j-m_1}{2}\right)}{\Gamma\left(\frac{j-m_2}{2}\right)}\times\nonumber\\
	&(1-t_1)^{\frac{-j+m_1}{2}}(1-t_2)^{\frac{j-m_2-2}{2}}t_2^{1-2j}t_1^{-1}\times\nonumber\\
	&{}_2F_1\left(\begin{smallmatrix}-j+m_1,\frac{\lambda_1-\lambda_2-2j-1}{2}\\-2j\end{smallmatrix};t_1\right){}_2F_1\left(\begin{smallmatrix}-j-m_2,\frac{\lambda_1+\lambda_2-2j-1}{2}\\-2j\end{smallmatrix};t_2\right)\times\nonumber\\
	 &\,
   _5F_4\left(\begin{smallmatrix}-j+1,1,\frac{-j+m_2+2}{2},\frac{-j-n-\lambda_1+2}{2},\frac{-j-m_2+\lambda_1+\lambda_2+1}{2}\\-j+1,\frac{-j+m_1+2}{2},\frac{-j-n+\lambda_1+2}{2},\frac{-j-m_1-\lambda_1+\lambda_2+1}{2}+1\end{smallmatrix};\frac{t_2^2 \left(1-t_1\right)}{   t_1^2\left(1-t_2\right)
}\right)\label{oddcase}.
\end{align}
Therefore, if we set
$
	\epsilon^{j,n}_{\delta} = \left\{\begin{smallmatrix}0&j-n\equiv\delta_i\text{ }\mathrm{mod}\text{ } 2\\1&j-n\not\equiv\delta_i\text{ }\mathrm{mod}\text{ } 2\end{smallmatrix}\right.,
$
we can summarize the above two cases (\ref{evencase}),(\ref{oddcase}) into one single formula:
\begin{align}
	&[A(\lambda)]^{j,n}_{m_1,m_2}(t_1,t_2)=\nonumber\\
	&\frac{(-1)^{n}((2j)!)^2}{\left(\frac{\lambda_1+1}{2}\right)^{(\frac{j+n-\epsilon^{j,n}_{\delta}}{2})}\left(\frac{\lambda_1+1}{2}\right)^{(-\frac{j+n-\epsilon^{j,n}_{\delta}}{2})}\left(\frac{\lambda_1-\lambda_2+1}{2}\right)^{(j)}\left(\frac{\lambda_1+\lambda_2+1}{2}\right)^{(j)}}\times\nonumber\\
	&\frac{1}{c^j_{m_1}c^j_{m_2}}\frac{\left(\frac{\lambda_1-\lambda_2}{2}\right)^{(\frac{j+m_1-\epsilon^{j,n}_{\delta}}{2})}\left(\frac{\lambda_1+\lambda_2}{2}\right)^{(\frac{-j-m_2+\epsilon^{j,n}_{\delta}}{2})}}{\left(\frac{\lambda_2+1}{2}\right)^{\frac{m_2-n}{2}}\left(\frac{\lambda_2+1}{2}\right)^{-\frac{m_2-n}{2}}}\frac{\Gamma\left(\frac{1-\epsilon^{j,n}_{\delta}+j-m_1}{2}\right)}{\Gamma\left(\frac{1-\epsilon^{j,n}_{\delta}+j-m_2}{2}\right)}\times\nonumber\\
	&(1-t_1)^{\frac{-1+\epsilon^{j,n}_{\delta}-j+m_1}{2}}(1-t_2)^{\frac{-1-\epsilon^{j,n}_{\delta}+j-m_2}{2}}t_2^{\epsilon^{j,n}_{\delta}-2j}t_1^{-\epsilon^{j,n}_{\delta}}\times\nonumber\\
	&{}_2F_1\left(\begin{smallmatrix}-j+m_1,\frac{\lambda_1-\lambda_2-2j-1}{2}\\-2j\end{smallmatrix};t_1\right){}_2F_1\left(\begin{smallmatrix}-j-m_2,\frac{\lambda_1+\lambda_2-2j-1}{2}\\-2j\end{smallmatrix};t_2\right)\times\nonumber\\
	 &\,
   _5F_4\left(\begin{smallmatrix}-j+\epsilon^{j,n}_{\delta},1,\frac{-j+m_2+1+\epsilon^{j,n}_{\delta}}{2},\frac{-j-n-\lambda_1+1+\epsilon^{j,n}_{\delta}}{2},\frac{-j-m_2+\lambda_1+\lambda_2+\epsilon^{j,n}_{\delta}}{2}\\-j+\epsilon^{j,n}_{\delta},\frac{-j+m_1+1+\epsilon^{j,n}_{\delta}}{2},\frac{-j-n+\lambda_1+1+\epsilon^{j,n}_{\delta}}{2},\frac{-j-m_1-\lambda_1+\lambda_2+\epsilon^{j,n}_{\delta}}{2}+1\end{smallmatrix};\frac{t_2^2 \left(1-t_1\right)}{t_1^2\left(1-t_2\right)
   }\right).
\end{align}
Thus we have finished the proof of Theorem \ref{LongSp4R}.

\appendix
\section{An Introduction to Wigner $D$-Matrices and Hypergeometric Functions}\label{cptgrp}
This appendix is a summary of \cite[Section~2]{zz2019}. The Lie algebra $\fraksl(2,\CC)$ is generated by the  \emph{Pauli matrices} :
\begin{align*}
	\sigma_0 = \begin{spmatrix}1&0\\0&1\end{spmatrix},	\sigma_1 =\begin{spmatrix}0&1\\1&0\end{spmatrix},
	\sigma_2 = \begin{spmatrix}0&-\ii\\\ii&0\end{spmatrix},\sigma_3 = \begin{spmatrix}1&0\\0&-1\end{spmatrix}
\end{align*}
and its compact real form $\fraku(2)$ is generated by $\gamma_i=\frac{\ii}{2}\sigma_i$:
\begin{align*}
	\gamma_0 = \begin{spmatrix}\frac{\ii}{2}&0\\0&\frac{\ii}{2}\end{spmatrix},	\gamma_1 =\begin{spmatrix}0&\frac{\ii}{2}\\\frac{\ii}{2}&0\end{spmatrix},
	\gamma_2 = \begin{spmatrix}0&\frac{1}{2}\\-\frac{1}{2}&0\end{spmatrix},\gamma_3 = \begin{spmatrix}\frac{\ii}{2}&0\\0&-\frac{\ii}{2}\end{spmatrix}.
\end{align*}
These generators exponentiate to an \emph{Euler angle} parametrization of the compact Lie group $U(2)$:
\begin{align}
	e^{-\zeta\gamma_0}\mathcal{U}(\psi, \theta, \phi)= e^{-\zeta\gamma_0}e^{-\psi\gamma_3}e^{-\theta\gamma_2}e^{-\phi\gamma_3}= \begin{spmatrix}
	e^{\frac{\ii}{2}(-\zeta-\phi-\psi)}\cos\frac{\theta}{2}&-e^{\frac{\ii}{2}(-\zeta+\phi-\psi)}\sin\frac{\theta}{2}\\
	e^{\frac{\ii}{2}(-\zeta-\phi+\psi)}\sin\frac{\theta}{2}&e^{\frac{\ii}{2}(-\zeta+\phi+\psi)}\cos\frac{\theta}{2}
	\end{spmatrix}\label{GenericSU2Matrix}.
\end{align}
An arbitrary irreducible representation $\pi_{j,n}$ of $U(2)$ can be realized on the space $\Sym^{2j}W\otimes \det^{(j+n)}$, on which $U(2)$ acts by right regular action on any degree $2j$ homogeneous polynomial $f\in \Sym^{2j}(W)\otimes \det^{(j+n)}$ in 2 variables $z_1,z_2$:
\[
	\pi_{j,n}(g)f(z)=(\det g)^{j+n}f(g^{-1}z).
\]
Let $m\in\frac{1}{2}\ZZ$ such that $-j\leq m\leq j$ and $j\pm m$ are integers, the weight basis $\{v^j_m\}_{-j\leq m\leq j}$ for $\Sym^{2j}W \otimes \det^{(j+n)}$ is defined as:
\[
	v^j_m = \frac{z_1^{j-m}z_2^{j+m}}{\sqrt{(j-m)!(j+m)!}}.
\]
There is a hermitian inner product $\langle v^j_{m_1}, v^j_{m_2} \rangle = \delta_{m_1,m_2}
$ on $\Sym^{2j} W\otimes \det^{(j+n)}$.  The \emph{Wigner $D$-functions} $W^{(j,n)}_{m_1,m_2}(\zeta,\psi,\theta,\phi)$ are the matrix coefficients of the irreducible representation $\pi_{j,n}$ under this hermitian inner product:
\begin{align}
	W^{(j,n)}_{m_1,m_2}(\zeta,\psi,\theta,\phi)&=\langle v^j_{m_1},e^{-\gamma_0\zeta}\mathcal{U}(\psi,\theta,\phi)v^j_{m_2}\rangle\nonumber\\
	&= c^j_{m_1}c^j_{m_2}e^{\ii n\zeta}e^{\ii(m_1\psi+m_2\phi)}d^{(j,n)}_{m_1,m_2}(\theta)\label{WignerDefinition},
\end{align}
where $c^j_m = \sqrt{(j+m)!(j-m)!}$ is a normalization factor, and the function $d^{(j,n)}_{m_1,m_2}(\theta)$ is given by the trigonometric polynomial
\begin{align}
d^{(j,n)}_{m_1,m_2}(\theta)&=\sum_{p=\max(0,m_1-m_2)}^{\min(j-m_2,j+m_1)}\frac{(-1)^{m_2-m_1+p}}{(j+m_1-p)!p!(m_2-m_1+p)!(j-m_2-p)!}\nonumber\\&\sin^{m_2-m_1+2p}\left(\frac{\theta}{2}\right)\cos^{2 j + m_1 - m_2 - 2 p}\left(\frac{\theta}{2}\right).\label{DPart}
\end{align}
which can also be expressed in terms of Jacobi polynomials:
\begin{align}
	&d^{(j,n)}_{m_1,m_2}(\theta)=\frac{\left(\sin\frac{\theta}{2}\right)^{m_1-m_2}\left(\cos\frac{\theta}{2}\right)^{m_1+m_2}}{(j+m_2)!(j-m_2)!}P^{(m_1-m_2,m_1+m_2)}_{j-m_1}(\cos\theta).\label{WignerJacobi}
\end{align}
For $n\geq 0$ and for $\alpha,\beta\in\RR$, the \emph{Jacobi polynomials} $P_{n}^{{(\alpha ,\beta )}}(z)$ are a class of orthogonal polynomials defined in \cite{AbStegun} as
\begin{align}
	P_{n}^{{(\alpha ,\beta )}}(z)={\frac {\Gamma (\alpha +n+1)}{n!\,\Gamma (\alpha +\beta +n+1)}}\sum _{{m=0}}^{n}{n \choose m}{\frac {\Gamma (\alpha +\beta +n+m+1)}{\Gamma (\alpha +m+1)}}\left({\frac {z-1}{2}}\right)^{m}.\label{Jacobi1}
\end{align}
In \cite{WolframAlpha}, the Jacobi polynomial is also defined as
\begin{align}\label{Jacobi2}
P^{\alpha,\beta}_n(x)&=\binom{n+\alpha}{n}\left(\frac{x+1}{2}\right)^n {}_2F_1\left(\begin{smallmatrix}-n,-n-\beta,\alpha+1\end{smallmatrix};\frac{x-1}{x+1}\right).
\end{align}

There is a left and right regular action by $\gamma_i$ on these Wigner $D$-matrices:
\begin{align}
\begin{matrix}
\dr(\gamma_0)W^{(j,n)}_{m_1,m_2}=-\ii nW^{(j,n)}_{m_1,m_2}& \dl(\gamma_0)W^{(j,n)}_{m_1,m_2}=\ii nW^{(j,n)}_{m_1,m_2}\\
	\dr(\gamma_3)W^{(j,n)}_{m_1,m_2}=-\ii m_2W^{(j,n)}_{m_1,m_2}& \dl(\gamma_3)W^{(j,n)}_{m_1,m_2}=\ii m_1W^{(j,n)}_{m_1,m_2}\\
\end{matrix}\label{CompactAction:1}\\
	\dr(\gamma_1\pm \ii\gamma_2)W^{(j,n)}_{m_1,m_2}=\ii\sqrt{(j\pm m_2)(j\mp m_2+1)}W^{(j,n)}_{m_1,m_2\mp 1}\label{CompactAction:2}\\
	\dl(\gamma_1\pm \ii\gamma_2)W^{(j,n)}_{m_1,m_2}=-\ii\sqrt{(j\mp m_1)(j\pm m_1+1)}W^{(j,n)}_{m_1\pm 1,m_2}\label{CompactAction:3}.
\end{align}

The product of Wigner $D$-functions $W^{(j_1,n_1)}_{m_{11},m_{12}}W^{(j_2,n_2)}_{m_{21},m_{22}}$ can be expressed in terms of other Wigner $D$-matrices in terms of the \emph{Clebsch-Gordan coefficients}:
\begin{align}
	W^{(j_1,n_1)}_{m_{11},m_{12}}W^{(j_2,n_2)}_{m_{21},m_{22}} = \sum_{\substack{|j_1-j_2|\leq J \leq j_1+j_2\\J-|j_1-j_2|\in\ZZ\\M_1=m_{11}+m_{21}\\M_2=m_{12}+m_{22}}}\left( \begin{smallmatrix}J, M_1\\j_1,m_{11},j_2,m_{21}\end{smallmatrix}\right) \left( \begin{smallmatrix}J, M_2\\j_1,m_{12},j_2,m_{22}\end{smallmatrix}\right)W^{(J,n_1+n_2)}_{M_1,M_2}.\label{ProductClebschGordan}
\end{align}

The Clebsch-Gordan coefficients for $V^j\otimes V^1$ are listed in the following table:
\begin{table}[h]
\centering
	\begin{tabular}{cccc}
		$\begin{spmatrix}j+j_0,m_1+m_2\\j,m_1,1,m_2\end{spmatrix}$& $m_2=-1$ & $m_2=0$ & $m_2=1$\\\hline
		$j_0=-1$ & $\sqrt{\frac{(j+m_1)(j+m_1-1)}{2j(2j+1)}}$ & $-\sqrt{\frac{(j-m_1)(j+m_1)}{j(2j+1)}}$ & $\sqrt{\frac{(j-m_1)(j-m_1-1)}{2j(2j+1)}}$\\
		$j_0=0$ & $\sqrt{\frac{(j+m_1)(j-m_1+1)}{2j(j+1)}}$ & $\frac{m_1}{\sqrt{j(j+1)}}$ & $-\sqrt{\frac{(j-m_1)(j+m_1+1)}{2j(j+1)}}$\\
		$j_0 =1$ & $\sqrt{\frac{(j-m_1+1)(j-m_1+2)}{(2j+2)(2j+1)}}$ &  $\sqrt{\frac{(j-m_1+1)(j+m_1+1)}{(j+1)(2j+1)}}$ & $\sqrt{\frac{(j+m_1+1)(j+m_1+2)}{(2j+2)(2j+1)}}$
	\end{tabular}.
	\caption{Table for Clebsch-Gordan coefficients of $V^{j}\otimes V^{1}$}\label{TableCG2}
\end{table}

\bibliographystyle{alpha}
\bibliography{Sp4}
\end{document}